\documentclass{smfart}
\usepackage[all]{xy}

\newtheorem{theorem}{Th{\'e}or{\`e}me}
\newtheorem{lemma}{Lemme}
\newtheorem{proposition}{Proposition}
\newtheorem{corollary}{Corollaire}
\newtheorem{pretheorem}{Th{\'e}or{\`e}me}[section]
\newtheorem{preproposition}{Proposition}[section]

\theoremstyle{definition}

\newtheorem{example}{Exemple}

\newtheorem{example1}{Exemple}

\theoremstyle{remark}
\newtheorem{remark}{Remarque}
\numberwithin{equation}{section}

\newcommand{\internalcomment}[1]{}

\begin{document}

\title[approximation par vari{\'e}t{\'e}s analytiques]{R{\'e}sultat n{\'e}gatif en
  th{\'e}orie d'approximation de compacts fonctionnels par des vari{\'e}t{\'e}s
  analytiques et application {\`a} un probl{\`e}me inverse}
\author{Amadeo Irigoyen}
\address{Universit{\'e} Paris VI - Pierre et Marie Curie, UMR 7586, 175,
  rue du Chevaleret 75013 Paris}
\date{30 janvier 2006}

\email{
\begin{minipage}[t]{5cm}
amadeo@math.jussieu.fr
\end{minipage}
}

\begin{abstract}

Dans la th{\'e}orie d'approximation figurent en particulier les probl{\`e}mes
d'approximation de compacts dans les espaces fonctionnels, par des
familles non lin{\'e}aires : on rappelle le cas de la param{\'e}trisation
polynomiale, puis on s'int{\'e}ressera au cas analytique. On montre un
r{\'e}sultat n{\'e}gatif qui dit
qu'une famille de fonctions param{\'e}tr{\'e}e analytiquement par $N$
variables ne peut pas approximer le compact $\Lambda_l(I^s)$ mieux que
de l'ordre de $(N\log N)^{-\frac{l}{s}}$, lorsque $N$ augmente.

Cette assertion fournit, comme application {\`a} un probl{\`e}me inverse dans
la th{\'e}orie de 
Sturm-Liouville, une r{\'e}ponse {\`a} une question sur la meilleure
reconstruction possible du potentiel $Q$, n{\'e}gatif avec
$m+1$ d{\'e}riv{\'e}es int{\'e}grables, {\`a} partir des valeurs propres et valeurs
caract{\'e}ristiques de l'{\'e}quation $-y''+\omega^2Qy=\lambda y$, lorsque
$\omega$ augmente : on montre l'impossibilit{\'e} d'avoir une formule
analytique d'approximation qui puisse avoir une pr{\'e}cision
meilleure que de l'ordre de $(\omega\log\omega)^{-(m+1)}$. Il
existe en outre dans~\cite{henkin} des formules d'approximation
presque optimales. 

\end{abstract}

\maketitle

\begin{altabstract}

In the theory of approximation there are some problems on approximation of
compacts in functional spaces by nonlinear families : first we deal with the
polynomial case, and then we consider the analytic case. We demonstrate a
negative result in which we claim that an analytic family of functions
with $N$ parameters can not approach the compact $\Lambda_l(I^s)$ closer
than of order $(N\log N)^{-\frac{l}{s}}$, when $N$ increases.

As applied to an inverse problem in
Sturm-Liouville theory, this assertion provides an answer to a
question about the best 
possible reconstruction of the negative potential $Q$
with $m+1$ integrable derivatives, from its eigenvalues and
characteristic values of the equation $-y''+\omega^2Qy=\lambda y$, when
$\omega$ increases : we show that it is impossible to get an analytic
approximating formula with precision better than of order
$(\omega\log\omega)^{-(m+1)}$. Moreover there is
in~\cite{henkin} formulas which are almost optimal.

\end{altabstract}

\maketitle

\tableofcontents

\section{Introduction}
\bigskip

On s'int{\'e}resse ici {\`a} un probl{\`e}me d'approximation de compacts dans
certains espaces fonctionnels par des sous-espaces non
lin{\'e}aires. Cette th{\'e}orie est li{\'e}e {\`a} l'{\'e}tude d'entropie d'un compact
dans un espace 
vectoriel norm{\'e}, concept d{\'e}fini par Kolmogorov comme {\'e}tant le
logarithme du nombre minimal de boules 
de rayon $\varepsilon$ pour recouvrir un compact $K$. Cette notion
donne une id{\'e}e de la possibilit{\'e} ou de la 
difficult{\'e} d'approximer un compact par un ensemble donn{\'e} (cf
~\cite{kolmogorov}, ~\cite{lorentz} et  ~\cite{vitushkin}).

Si $K$ est un compact d'un espace vectoriel norm{\'e} $L$, et $F$ un
sous-ensemble de $L$, on pose 
$$D(K,F)=\sup_{y\in K}\inf_{z\in F}\|y-z\|$$
qui repr{\'e}sente le d{\'e}faut d'approximation de $K$ par $F$. Les cas les
plus {\'e}tudi{\'e}s sont ceux o{\`u} $L$ d{\'e}signe un espace de fonctions
continues ou int{\'e}grables sur $I^s$ et $K=\Lambda_{l,s}$ d{\'e}signe la
boule unit{\'e} des fonctions de  
classe $C^l,\,l>0$ sur $I^s$, ou encore la boule
unit{\'e} des restrictions de fonctions analytiques born{\'e}es sur un
voisinage de $I^s$ dans $\mathbb{C}^s$.

Un des principaux exemples est celui de l'approximation lin{\'e}aire,
c'est-{\`a}-dire par des sous-espaces vectoriels 
de dimension finie. Si $n$ est un entier $\geq1$, on pose
$$D_n(K)=\inf_{F_n}D(K,F_n),$$ o{\`u} $F_n$ parcourt l'ensemble
$\mathcal{L}_n$ des 
sous-espaces vectoriels de dimension $n$ de $L$. $D_n(K)$ est appel{\'e}
$n-width$ 
du compact $K$. Dans la suite des travaux de Vitushkin
(cf~\cite{vitushkin}), Tihomirov (cf~\cite{tihomirov1}) {\'e}tablit
quelques propri{\'e}t{\'e}s g{\'e}n{\'e}rales des $n$-widths ainsi que des
exemples de calcul avec explicitation de $F_n$ r{\'e}alisant
le minimum (cf~\cite{pinkus}).

On en particulier, si $\Lambda_{l,s}$ d{\'e}signe la boule unit{\'e} de
$C^l(I^s)\subset C^0(I^s)$, $I=[0,1]$, $$\frac{A}{n^{\frac{l}{s}}}\leq
D_n(\Lambda_{l,s})\leq\frac{B}{n^{\frac{l}{s}}},$$
ce qui montre entre autres une minoration, donc un r{\'e}sultat n{\'e}gatif
d'approximation.

Dans le cadre non lin{\'e}aire, Vitushkin consid{\`e}re dans~\cite{vitushkin}
le cas suivant : 
si $n$ et $d$ sont des entiers $\geq1$, on se donne
l'ensemble 
$$P_{n,d}=\left\{P(\zeta)=\sum_{|k|\leq d}a_k(\cdot)\zeta^k,\,
\zeta=(\zeta_1,\ldots,\zeta_n)\in\mathbb{R}^n\right\},\;a_k\in C^0(I^s),$$
param{\'e}tr{\'e} polynomialement avec $n$ variables et de degr{\'e} $d$, et
$$D(K,P_{n,d})=\sup_{y\in K}\inf_{\zeta\in\mathbb{R}^n}\|y-P(\zeta)\|.$$
On pose finalement 
$$D_{n,d}(K)=\inf_{P_{n,d}}D(K,P_{n,d}),$$
o{\`u} $P_{n,d}$ parcourt l'ensemble $\mathcal{P}_{n,d}$ des vari{\'e}t{\'e}s de
$L$ param{\'e}tr{\'e}es par $n$ variables ind{\'e}pendantes, polynomialement de
degr{\'e} au 
plus $d$. On a comme r{\'e}sultat principal le th{\'e}or{\`e}me de Vitushkin que
l'on peut trouver dans 
~\cite{ivanov}, ~\cite{lorentz}, ~\cite{vitushkin} et~\cite{warren2} :
$$\frac{A}{(n\log d)^{\frac{l}{s}}}\leq
D_{n,d}(\Lambda_{l,s})\leq\frac{B}{(n\log d)^{\frac{l}{s}}},$$
ce qui est une g{\'e}n{\'e}ralisation du cas lin{\'e}aire, et montre en
particulier qu'augmenter la complexit{\'e} par le degr{\'e} n'apportera pas
d'am{\'e}lioration significative. On souligne en particulier des r{\'e}sultats
n{\'e}gatifs d'approximation, i.e. des minorations de $D_{n,d}$.

Vitushkin a montr{\'e} cette estimation en utilisant la th{\'e}orie de
$variation$ $of$ $sets$, sans en expliciter les constantes.
\bigskip

Plus tard, Warren a montr{\'e} ce r{\'e}sultat dans ~\cite{warren}, sans
utilisation des $variations$ $of$ $sets$, en
explicitant les constantes, mais seulement dans des cas
particuliers : $\Lambda_{\omega}([0,1])$ dans le cas uniforme (o{\`u}
$\omega$ est un module de continuit{\'e}), et
$\Lambda_{\alpha,s}$, $0<\alpha\leq1$, dans le cas
$L^1(I^s)$. La m{\'e}thode est diff{\'e}rente, mais elle utilise comme
Vitushkin une estimation du nombre de composantes connexes d'un
ensemble alg{\'e}brique 
dans $\mathbb{R}^n$, de Oleinik et Petrovskii (cf~\cite{oleinik}). Ce
dernier repose finalement sur le th{\'e}or{\`e}me
de B{\'e}zout : une intersection de $n$ ensembles alg{\'e}briques
$\{\zeta\in\mathbb{R}^n, P_j(\zeta)=0\}$, o{\`u} $\deg P_j=p_j$, ne peut avoir
plus de $\prod_{j=1}^np_j$ points, les ensembles de z{\'e}ros se trouvant
en position g{\'e}n{\'e}rique (i.e. les coefficients de ces polyn{\^o}mes {\'e}tant pris
dans un ouvert de Zariski).

Cela permet {\`a} Warren d'obtenir
l'estimation suivante, {\'e}crite en d{\'e}tail dans~\cite{warren}: si
$P_1,\ldots,P_m$ sont des polyn{\^o}mes sur 
$\mathbb{R}^n$, de degr{\'e} $\leq d$, alors le nombre de composantes
connexes de l'ensemble 
$$\mathbb{R}^n\setminus\bigcup_{j=1}^n\{P_j(\zeta)=0\}$$
ne peut d{\'e}passer $\left(\frac{4edm}{n}\right)^n$.
On en d{\'e}duit ainsi une estimation du nombre de suites
$\varepsilon=(\varepsilon_1,\ldots,\varepsilon_m)$,
$\varepsilon_j=\pm1$, prises par la fonction :
$\zeta\in\mathbb{R}^n\mapsto(sgn\,P_1(\zeta),\ldots,sgn\,P_m(\zeta))$. 
Il en r{\'e}sulte que si $K$ poss{\`e}de des fonctions qui oscillent suffisamment,
il existera des points o{\`u} des {\'e}l{\'e}ments de $K$ diff{\`e}reront de
l'ensemble polynomial, ce 
qui permettra d'obtenir une minoration de $D_{n,d}(K)$.

Notre premier but est de d{\'e}montrer le th{\'e}or{\`e}me de Vitushkin dans la
partie \ref{vitwar}, dans le cas g{\'e}n{\'e}ral $\Lambda_{l,s}$ et en
explicitant la constante $A$ : 
\begin{pretheorem}

Soit $\mathcal{P}_{n,d}$ l'ensemble des familles polynomiales de
$C(I^s)$ (resp. $L^1(I^s)$), i.e. form{\'e}e des {\'e}l{\'e}ments
$$P_{n,d}=\left\{\sum_{|k|\leq d}c_k\,\zeta^k,\,\zeta=
(\zeta_1,\ldots,\zeta_n)\in\mathbb{R}^n\right\},$$
o{\`u} $c_k\in C(I^s)$ (resp. $L^1(I^s)$). On a alors
$$D_{n,d}(\Lambda_{l,s})\geq\frac{C_{l,s}}{(n\log
  d)^{\frac{l}{s}}}\,.$$

\end{pretheorem}
On donne en outre une estimation des constantes
$C_{\infty}(l,s)$ et $C_{L^1}(l,s)$.
\bigskip

Dans le cadre de la th{\'e}orie d'approximation non lin{\'e}aire, figure aussi
l'{\'e}tude d'approximation rationnelle, que l'on peut trouver dans
~\cite{shapiro} et~\cite{vitushkin}. Ici aussi, passer du cas
polynomial au cas rationnel n'am{\'e}liore pas essentiellement la pr{\'e}cision.
\bigskip

En vue d'applications pour des probl{\`e}mes inverses, on g{\'e}n{\'e}ralise
l'approximation non lin{\'e}aire sur une classe assez 
naturelle de fonctions analytiques qui ne sont pas des polyn{\^o}mes, qui est
celle des fonctions enti{\`e}res de type 
exponentiel, i.e. qui v{\'e}rifient : $$\forall z\in\mathbb{C}^n, |f(z)|\leq
a\exp\left(b\|z\|_1^d\right),$$
o{\`u} $\|z\|_1=|z_1|+\ldots+|z_n|$. Ce sera l'objet de la
partie~\ref{appranal}, o{\`u} on consid{\`e}re un ensemble
param{\'e}tr{\'e} analytiquement par $N$ variables ind{\'e}pendantes, {\`a} valeurs
dans $C(I^s)$, soit une fonction de la forme
$$f(x,\zeta)=f(x_1,\ldots,x_s,\zeta_1,\ldots,\zeta_N),\,x\in I^s,\,
\zeta\in\mathbb{R}^N,$$
o{\`u}, pour chaque $\zeta$, $f$ est continue par rapport {\`a} $x$, et pour
chaque $x\in I^s$, $f$ est enti{\`e}re de type exponentiel par rapport {\`a}
$\zeta$. L'hypoth{\`e}se du type exponentiel pour $f$ se r{\'e}{\'e}crit avec
$\|f(\cdot,\zeta)\|_{\infty}$. On doit n{\'e}anmoins adopter une
restriction sur le domaine de 
d{\'e}finition : la variable $\zeta$ reste dans un compact dont la taille
grandit avec $N$. Le r{\'e}sultat est alors le suivant : 

\begin{pretheorem}\label{pth1}

Soit $\mathcal{E}_N$ la classe des familles enti{\`e}res de $C(I^s)$
(resp. $L^1(I^s)$), i.e. de la forme
$$E_N=\left\{f(\cdot,\zeta)=\left(x\in I^s\mapsto f(x_1,\ldots,x_s,\,
\zeta_1,\ldots,\zeta_N)\right),\;\zeta\in\mathbb{R}^N,\,|\zeta_j|
\leq BN^r,1\leq j\leq N\right\},$$
o{\`u} $f$ est une fonction enti{\`e}re de $\zeta$ v{\'e}rifiant
$$\forall\,\zeta\in\mathbb{C}^N,\,
\|f(\cdot,\zeta)\|\leq Ae^{uN^v}e^{bN^t\|\zeta\|_1^d}\,.$$
Alors
$$D_N(\Lambda_{l,s})=\inf_{f\in\mathcal{E}_N}\,\sup_{h\in\Lambda_{l,s}}\,
\inf_{|\zeta_j|\leq BN^r}\|h-f(\cdot,\zeta)\|\geq\frac{C}{\left
(N\log N\right)^{\frac{l}{s}}}\,,$$
o{\`u} $C=C(l,s,A,u,v,b,t,d,B,r)$.

\end{pretheorem}

Ici aussi, on explicite les constantes $C_{\infty}$ et
$C_{L^1}$.
\bigskip

Pour finir dans le cadre d'approximation par des ensembles
analytiques, on a {\'e}galement consid{\'e}r{\'e} une classe particuli{\`e}re de fonctions
enti{\`e}res de type exponentiel appel{\'e}es quasipolyn{\^o}mes, et qui sont de la
forme $$P(\zeta_1,\ldots,\zeta_n,\exp<a_1,\zeta>,\ldots,\exp<a_k,\zeta>),$$
o{\`u} $\zeta\in\mathbb{R}^n$, $a_1,\ldots,a_k\in\mathbb{R}^n$,
$<a_j,\zeta>=a_j^1\zeta_1+\ldots+a_j^n\zeta_n$. En utilisant 
la m{\'e}thode de Warren, on peut obtenir des r{\'e}sultats n{\'e}gatifs
d'approximation de $\Lambda_{l,s}$ par des familles non lin{\'e}aires qui
sont des 
quasipolyn{\^o}mes {\`a} coefficients dans $L$ ($=C(I^s)$ ou $L^1(I^s)$). On a
besoin pour cela 
d'estimations du nombre de racines d'un syst{\`e}me de $n$ quasipolyn{\^o}mes,
donc de type B{\'e}zout. Khovanskii a pour cela {\'e}tabli
dans~\cite{khovanskii} plusieurs r{\'e}sultats sur les {\'e}quations de Pfaff,
o{\`u} il a  donn{\'e} des exemples avec estimations explicites : si on
consid{\`e}re un syst{\`e}me 
$$P_1=\ldots=P_p=0$$ 
de $p$ {\'e}quations quasi-polynomiales en $\zeta$, o{\`u} chaque $P_j$ est de
degr{\'e} $m_j$ en les 
$n+k$ variables $\zeta_1,\ldots,\zeta_n,\,y_1,\ldots,y_k$, avec
$y_k=\exp<a_k,\zeta>$, alors
le nombre de cellules ({\`a} homotopie pr{\`e}s) de l'ensemble solution de
dimension $n-p$ sera major{\'e} par 
$$2^{\frac{k(k-1)}{2}}m_1\ldots m_p\left(\sum_{j=1}^pm_j+n-p+1\right)^
{n-p}\left[(n-p+1)\left(\sum_{j=1}^pm_j+n-p+1\right)-n+p\right]^k,$$
ce qui donnera une estimation du m{\^e}me type pour le nombre de
composantes de l'ensemble
$\mathbb{R}^n\setminus\bigcup_{j=1}^m\{P_j=0\}$, et en particulier
le r{\'e}sultat suivant,
que l'on prouvera en derni{\`e}re partie :

\begin{pretheorem}\label{prekhovanskii}

Soit $\mathcal{P}_{n,k,d}$ l'ensemble des familles de $C(I^s)$
repr{\'e}sent{\'e}es par
$$P_{n,k,d}=\left\{\sum_{|j|\leq d}c_j\,\zeta_1^{j_1}\ldots\zeta_n^{j_n}\,
e^{j_{n+1}<a_1,\zeta>}\ldots e^{j_{n+k}<a_k,\zeta>},\;
\zeta\in\mathbb{R}^n\right\},$$
o{\`u} $P_{n,k,d}$ est un quasi-polyn{\^o}me {\`a}
coefficients $c_j\in C(I^s)$, {\`a} $n$ variables $\zeta_i$, $k$
pseudo-variables $e^{<a_i,\zeta>}$ et de degr{\'e} total $d$.
On a alors
$$D_{n,k,d}(\Lambda_{l,s}):=\inf_{P\in\mathcal{P}_{n,k,d}}\sup_{h\in
\Lambda_{l,s}}\,\inf_{\zeta\in\mathbb{R}^n}\|h-P(\zeta)\|_{\infty}\geq
\frac{C_{l,s}}{(k^2n\log n\log d)^{\frac{l}{s}}}\,.$$

En outre, la constante $C_{l,s}$ peut {\^e}tre calcul{\'e}e et vaut
$$\frac{1}{\sqrt{s}\,2^{l+1}\,38^{\frac{l}{s}}([l]+1)^{[l]+1}
(4(1+e))^{s([l]+1)}}\,.$$

\end{pretheorem}

Ainsi, une famille quasi-polynomiale de
$C(I^s)$ ne peut pas approcher le compact $\Lambda_{l,s}$ mieux qu'{\`a}
l'ordre
$$\frac{1}{\left(k^2n\log n\log d\right)^{\frac{l}{s}}}.$$
Cependant la pr{\'e}sence du terme $k^2$ semble brutale, ce qui
donne un r{\'e}sultat n{\'e}gatif relativement faible. En revanche, on n'a ici
aucune restriction sur les param{\`e}tres $\zeta_1,\ldots,\zeta_n\in\mathbb{R}$.
\bigskip

Dans la partie~\ref{pbinverse} on donne une application du
th{\'e}or{\`e}me~\ref{pth1} {\`a} un probl{\`e}me inverse dans la th{\'e}orie de Sturm-Liouville
sur la demi-droite $\mathbb{R}^+$. On consid{\`e}re l'{\'e}quation
$$-\frac{d^2y}{dx^2}-\omega^2Qy=\lambda y,$$
o{\`u} $-\omega^2Q$ est un potentiel strictement n{\'e}gatif et suffisamment
r{\'e}gulier, $\omega$ un 
param{\`e}tre grand : dans notre cas, $Q$ aura $m+1$ d{\'e}riv{\'e}es localement
int{\'e}grables et {\`a} d{\'e}croissance polynomiale. On sait alors que
l'op{\'e}rateur $-\frac{d^2}{dx^2}-\omega^2Q$ admet $N(\omega)$ valeurs
propres $-\xi_j^2$ n{\'e}gatives, avec pour chacune une fonction propre
$\psi_j$ associ{\'e}e v{\'e}rifiant les conditions aux bords :
$$\psi_j(0)=0\text{ et }\int_0^{\infty}\psi_j^2(x)dx=1\,.$$
On pose alors $C_j=(\psi_j'(0))^2$, qui est la valeur caract{\'e}ristique
associ{\'e}e {\`a} $\xi_j$.

Le probl{\`e}me est le suivant : connaissant les valeurs propres $\xi_j$
et valeurs caract{\'e}ristiques $C_j$ de l'{\'e}quation, il s'agit de
reconstruire le potentiel $-\omega^2Q$ sur 
$\mathbb{R}^+$. On dispose pour cela de formules plus ou
moins explicites (cf~\cite{lax-levermore}, avec convergence dans
$L^2$, et ~\cite{henkin}), provenant en
particulier des travaux de Gelfand,
Levitan, Kohn et Jost (cf ~\cite{gelfand}, ~\cite{levitan}). En outre,
G. Henkin et N. Novikova, motiv{\'e}s en particulier par Lax et Levermore, 
donnent dans~\cite{henkin} des estimations sur la vitesse de
convergence : si $Q$ poss{\`e}de $m+1$ d{\'e}riv{\'e}es int{\'e}grables, il existe des
formules qui l'approximent uniform{\'e}ment sur tout 
intervalle $[0,X]$ lorsque $\omega\rightarrow\infty$, avec une
vitesse (au moins) de l'ordre de $\frac{1}{\omega^m}$.

La question que l'on se pose alors (cf~\cite{henkin}, p. 22) est 
de savoir si on peut am{\'e}liorer cette approximation,
i.e. montrer s'il existe une autre estimation ou bien une autre
formule (et si c'est le cas, l'expliciter) qui donne une convergence plus 
rapide vers le potentiel donn{\'e} $-\omega^2Q$.

Il y a deux cas que nous traitons, auxquels nous donnons une r{\'e}ponse
n{\'e}gative : le premier est celui o{\`u} $m=1$. On utilise
dans~\cite{henkin} une formule d'approximation explicite du type
Gelfand-Levitan pour des potentiels avec propri{\'e}t{\'e}s {\'e}nonc{\'e}es ci-dessus : 
$$Q_{\omega}^0(x)=\frac{2}{\omega^2}\frac{d^2}{dx^2}\ln|\det(W_{s,r})(x)|,$$
o{\`u}
$$W_{s,r}(x)=\frac{2sh(\xi_r+\xi_s)x}{\xi_r+\xi_s}-(1-\delta_{s,r})
\frac{2sh(\xi_r-\xi_s)x}{\xi_r-\xi_s}-\delta_{s,r}\left(2x-
\frac{4\xi_j^2}{C_j}\right).$$
Cette formule approxime en primitives tout $Q$ uniform{\'e}ment sur
tout $[0,X]$, {\`a} la vitesse $\frac{\ln\omega}{\sqrt{\omega}}$ (mais
l'approximation peut {\^e}tre vraisemblablement bien meilleure, de l'ordre
de $\frac{1}{\omega^2}$).
On remarque que, mis {\`a} part le d{\'e}nominateur de $Q_{\omega}^0$
provenant du logarithme, $\det(W_{s,r}(x))$ est une fonction enti{\`e}re
de type exponentiel par rapport {\`a} 
$\left(\xi_1,\ldots,\xi_N,\frac{1}{C_1},\ldots,\frac{1}{C_N}\right)$, 
continue par 
rapport {\`a} $x$. On peut donc interpr{\'e}ter $Q_{\omega}^0$ comme une
famille de fonctions continues param{\'e}tr{\'e}e analytiquement. Pour cela,
on peut appliquer notre r{\'e}sultat n{\'e}gatif et le compl{\'e}ter avec le
r{\'e}sultat positif de~\cite{henkin}. On obtient :

\begin{pretheorem}

Soient $\mathcal{Q}$ la classe des fonctions $Q$ d{\'e}finies sur
$\mathbb{R}^+$, strictement positives, {\`a} d{\'e}croissance polynomiale, avec $2$
d{\'e}riv{\'e}es localement int{\'e}grables,
et leurs op{\'e}rateurs associ{\'e}s $-\frac{d^2}{dx^2}-\omega^2Q$ ; pour tout
$N(\omega)=O(\omega)$, soit {\'e}galement $\psi(x,\zeta)$ une fonction
d{\'e}finie sur  
$\mathbb{R}\times\mathbb{C}^N$, de classe $C^1$ par rapport {\`a} $x$ et
enti{\`e}re de type exponentiel par rapport {\`a} $\zeta\in\mathbb{C}^N$. Alors
l'approximation de
$$\int_0^{\cdot}\mathcal{Q}:=\left\{\left(x\mapsto\int_0^xQ(t)dt
\right),\;Q\in\mathcal{Q}\right\},$$
uniform{\'e}ment sur tout intervalle $[0,X]$, par la famille
$$\left\{\left(x\mapsto\left(\frac{1}{\psi}\frac{\partial\psi}{\partial x}
\right)(x,\zeta)\right), \,\zeta_j=O(\omega^r),
\,\forall\,j=1,\ldots,N\right\},$$
lorsque $\omega\rightarrow\infty$, ne peut pas {\^e}tre meilleure que de
l'ordre de $$\frac{1}{(\omega\ln\omega)^3}\,.$$

En outre, un cas d'approximation au moins {\`a} l'ordre
$\frac{\ln\omega}{\sqrt{\omega}}$,
est donn{\'e} par la formule de type Gelfand-Levitan
$$\Psi(x,\zeta)=\det\widetilde{W}_{s,r}(x,\zeta),$$
avec
$$\widetilde{W}_{s,r}(x,\zeta)=\frac{2sh(\zeta_r+\zeta_s)x}{\zeta_r+\zeta_s}
-(1-\delta_{s,r})\frac{2sh(\zeta_s-\zeta_r)x}{\zeta_s-\zeta_r}
-\delta_{s,r}(2x-\exp(\zeta_{r+N})),$$
$s,\,r=1,\ldots,N(\omega)$,
o{\`u} $N(\omega)$ est le nombre de
valeurs propres $\xi_j$ et caract{\'e}ristiques $C_j$ de l'op{\'e}rateur
$-\frac{d^2}{dx^2}-\omega^2Q$, dont l'{\'e}l{\'e}ment optimisant peut {\^e}tre
ainsi choisi : 
$$\zeta_j(Q)=\xi_j, \text{ et } \zeta_{j+N}(Q)=\ln\frac{4\xi_j^2}{C_j},
\text{ } j=1,\ldots,N(\omega).$$

\end{pretheorem}

On aura besoin pour la d{\'e}monstration d'{\'e}tablir une estimation des
valeurs propres $\xi_j$ et valeurs caract{\'e}ristiques $C_j$,
$j=1,\ldots,N(\omega)$, ainsi que de leurs inverses, {\`a} cause de la
restriction des param{\`e}tres $\zeta_j$, n{\'e}cessaire pour pouvoir appliquer
notre r{\'e}sultat n{\'e}gatif. Plus pr{\'e}cis{\'e}ment, on montrera la

\begin{preproposition}

Si $q=-\omega^2Q$ est un potentiel n{\'e}gatif, int{\'e}grable de classe
$C^1$ avec $q'(0)=0$, et {\`a} d{\'e}croissance polynomiale, alors pour tous
$\omega$ assez grand et $j=1,\ldots,N(\omega)$, on a 
$$\frac{1}{a\omega^b}\leq\xi_j\leq a\omega^b,\text{ et }
\frac{1}{\alpha\exp\left(\beta\omega^{\gamma}\right)}\leq
\frac{4\xi_j^2}{C_j}\leq\alpha \exp\left(\beta\omega^{\gamma}\right).$$

\end{preproposition}

Dans l'autre cas, qui n'est pas coml{\`e}tement fini, il s'agit de formules
d'approximation concernant le 
cas plus g{\'e}n{\'e}ral d'un potentiel avec $m+1$ d{\'e}riv{\'e}es, mais pas
explicit{\'e}es, et qui r{\'e}alisent l'approximation (au moins) {\`a} l'ordre
$$\frac{1}{\omega^m}\,.$$
Ces formules sont cependant, dans le cas o{\`u} les $m$ premi{\`e}res d{\'e}riv{\'e}es
de $Q$ s'annulent en $0$, de la forme
$$Q_{\omega}(x)=\frac{2}{\omega^2}\left(-\frac{d}{dx}K(x,x)+\frac{d^2}{dx^2}
\det T_{j,k}(x)\right)\,,$$
o{\`u}
$$K(x,x)=K(x,x,q(0))\text{ et }T_{j,k}(x)=T_{j,k}(x,\xi,C,q(0))\,.$$
Ce sont donc des formules qui, bien que non explicites, sont des
fonctions analytiques en les variables $\xi,\;C,\;q(0)=-\omega^2Q(0)$.
Ici encore, m{\^e}me si ce n'est pas prouv{\'e}, on est
tr{\`e}s optimiste sur la validit{\'e} d'un r{\'e}sultat n{\'e}gatif affirmant qu'il
n'est pas possible d'approximer $Q$, avec des formules du m{\^e}me type,
mieux qu'{\`a} l'ordre 
$$\frac{1}{(\omega\ln\omega)^{m+1}}\,,$$
ce qui montrerait encore une propri{\'e}t{\'e} de presque optimalit{\'e} (mais
cette fois, sans formule de reconstruction explicite).
\bigskip

On terminera en derni{\`e}re partie par quelques discussions sur d'autres
m{\'e}thodes pour obtenir des r{\'e}sultats analogues : on prouvera en
particulier le theor{\`e}me~\ref{prekhovanskii} sur les familles de
quasi-polyn{\^o}mes, et qui utilise les estimations de Khovanskii ; ainsi
qu'un autre r{\'e}sultat n{\'e}gatif, qui est une application directe du
th{\'e}or{\`e}me de 
Descartes. On constatera alors que le th{\'e}or{\`e}me~\ref{pth1} donne une
meilleure minoration, en plus d'{\^e}tre plus g{\'e}n{\'e}ral, malgr{\'e} la
restriction n{\'e}cessaire sur les param{\`e}tres.
\bigskip

\section{Remerciements}

Je remercie G. Henkin pour les id{\'e}es et discussions enrichissantes sur
ce sujet.

\section{Approximation analytique}\label{appranal}
\bigskip
\subsection{Enonc{\'e} et preuve du th{\'e}or{\`e}me}

On consid{\`e}re le compact $\Lambda_l(I^s)$, $I=[0,1]$,
$s\in{\mathbb{N}^{\ast}}$, $l>0$,
$$\Lambda_{l,s}=\Lambda_l(I^s)=\left\{f\in C^l(I^s),\,\forall\,j,0\leq j\leq
m,\,\left\|f^{(j)}\right\|_{\infty}\leq 1,\,\left\|f^{(m)}\right\|_{\alpha}
\leq 1\right\},$$
o{\`u} $l=m+\alpha$, $m\in\mathbb{N}$, $0<\alpha\leq 1$
($m=-[-l]-1$), et $$\|f\|_{\alpha}=\sup_{x\neq
  y}\frac{|f(x)-f(y)|}{\|x-y\|^{\alpha}},$$ o{\`u} $\|\cdot\|$ est la norme
euclidienne usuelle.

($\Lambda_l(I^s)$ est un compact de $(C^0(I^s),\|\cdot\|_{\infty})$
et de $(L^1(I^s),\|\cdot\|_1)$)
\bigskip

\begin{theorem}\label{th1}
On se donne, pour $N\in\mathbb{N},N\geq 2$,
$$f(x,\zeta)=f(x_1,\ldots,x_s,\zeta_1,\ldots,\zeta_N),\,x\in{I^s},\,
\zeta\in\mathbb{R}^N,$$
o{\`u} $f$ est enti{\`e}re par rapport {\`a} $\zeta\in\mathbb{C}^N$, continue
(resp. int{\'e}grable) par rapport {\`a} $x\in I^s$.

On suppose de plus que,
$$\forall\,\zeta\in\mathbb{C}^N,
\|f(\cdot,\zeta)\|\leq Ae^{uN^v}e^{bN^t\|\zeta\|_1^d},$$
o{\`u} $A,u,v,b,t,d$ sont des r{\'e}els $\geq 1$,
$\|f(\cdot,\zeta)\|=\|f(\cdot,\zeta)\|_{\infty}$
(resp. $\|f(\cdot,\zeta)\|_{L^1}$), et
$\|\zeta\|_1=|\zeta_1|+\ldots+|\zeta_N|$.
\bigskip

Alors on a : $\exists\,h\in\Lambda_{l,s}$,
$\forall\,\zeta\in\mathbb{R}^N,|\zeta_j|\leq BN^r,\,j=1,\ldots,N$, avec
$B,\,r\geq 1$,
$$\|h-f(\cdot,\zeta)\|_{\infty}\geq\frac{C}{\left(N\log_{[2]}N\right)^
{\frac{l}{s}}}\,\text{ (resp. }\|h-f(\cdot,\zeta)\|_{L^1}\geq\frac{C}
{\left(N\log_{[2]}N\right)^{\frac{l}{s}}}\text{ )},$$
o{\`u} $C$ est une constante d{\'e}pendant de $l,s,A,u,v,b,t,d,B,r$ (mais pas
de $N$) .
\bigskip

De fa{\c c}on {\'e}quivalente, si $\mathcal{E}_N$ d{\'e}signe la classe des fonctions
analytiques du type $f$ (donc d{\'e}pendant de $A,u,v,b,t,d$), alors
$$D_N(\Lambda_{l,s})=\inf_{f\in\mathcal{E}_N}\,\sup_{h\in\Lambda_{l,s}}\,
\inf_{|\zeta_j|\leq BN^r}\|h-f(\cdot,\zeta)\|\geq\frac{C}{\left(N\log_{[2]}
N\right)^{\frac{l}{s}}}.$$
\end{theorem}
\bigskip

Dans la suite, on notera simplement $\log$ pour le logarithme binaire
$\log_{[2]}$ ($\ln$ d{\'e}signera toujours le logarithme n{\'e}perien).
\begin{proof}

On commence par montrer le
\begin{lemma}
On consid{\`e}re le d{\'e}veloppement de Taylor de $f$ :
$\forall\,(x,\zeta)\in I^s\times\mathbb{R}^N$,
$$f(x,\zeta)=\sum_{|k|\geq 0}c_k(x)\zeta^k,$$
o{\`u} $k$ est le multi-indice $(k_1,\ldots,k_N)\in\mathbb{N}^N$, 
$|k|=k_1+\ldots+k_N$ et $\zeta^k=\zeta_1^{k_1}\ldots\zeta_N^{k_N}$.
\bigskip

Alors on a : $\forall k\in\mathbb{N}^N$,
$$\|c_k(\cdot)\|\leq Ae^{uN^v}\left(\frac{ebdN^{d+t}}{|k|}\right)^
{\frac{|k|}{d}},$$
o{\`u} la norme d{\'e}signe la norme uniforme ou la norme $L^1$.

\end{lemma}

\begin{proof}

On commence par remarquer que, pour tout $k$, la
fonction $x\mapsto c_k(x)$, qui donne le coefficient de Taylor du
mon{\^o}me $\zeta^k$, est bien d{\'e}finie, que ce soit dans l'espace $C(I^s)$
ou $L^1(I^s)$. En effet, elle est donn{\'e}e par la formule de Cauchy qui
utilise l'holomorphie de $f$ par rapport {\`a} $\zeta$ :
$$c_k(x)=\frac{1}{(2i\pi)^N}\int_{|\zeta_1|=\ldots=|\zeta_N|=R}
\frac{f(x,\zeta)}{\zeta^{k+1}}d\zeta,$$
o{\`u} $R>0$, $\zeta=(\zeta_1,\ldots,\zeta_N)$,
$d\zeta=d\zeta_1\wedge\ldots\wedge d\zeta_N$ et
$(k+1)=(k_1+1,\ldots,k_N+1)$.
\bigskip

Pour le cas continu, il s'agit d'une int{\'e}grale {\`a} param{\`e}tre d'une
fonction continue par rapport {\`a} $x\in{I^s}$ ; tandis que pour le cas $L^1$,
c'est le th{\'e}or{\`e}me de Fubini appliqu{\'e} {\`a} la fonction $f$ sur l'ensemble
$I^s\times (bD(0,R))^N$, qui permet de d{\'e}finir sur $I^s$ (presque
partout, puis 
en prolongeant par $0$) la fonction $c_k$, qui sera alors dans
$L^1(I^s)$.
\bigskip

On a donc : $\forall\,x\in I^s$,
$$|c_k(x)|\leq\frac{1}{(2\pi)^N}\int_{|\zeta_j|=R}\frac{|f(x,\zeta)|}
{R^{|k|+N}}d\zeta,$$
soit
$$\|c_k\|\leq\frac{1}{(2\pi)^N}\int_{|\zeta_j|=R}\frac{Ae^{uN^v}e^{bN^t
(|\zeta_1|+\ldots+|\zeta_N|)^d}}{R^{|k|+N}}d\zeta=Ae^{uN^v}\frac{e^
{bN^{d+t}R^d}}{R^{|k|}}.$$
\bigskip

L'in{\'e}galit{\'e} {\'e}tant valable pour tout $R>0$, on l'a en particulier pour
$R=R_{min}$ qui minimise le membre de droite. Pour d{\'e}terminer $R_{min}$, on
consid{\`e}re la fonction d{\'e}finie sur $]0,+\infty[$ par
$R\mapsto\frac{e^{bN^{d+t}R^d}}{R^{|k|}}$. Cette fonction tend vers
$+\infty$ en $0$ et $+\infty$, elle admet donc (au moins) un minimum,
qui annule sa d{\'e}riv{\'e}e, donn{\'e}e par :
$$R\mapsto\frac{e^{bN^{d+t}R^d}}{R^{|k|}}bdN^{d+t}R^{d-1}-|k|\frac{e^
{bN^{d+t}R^d}}{R^{|k|+1}}=\frac{e^{bN^{d+t}R^d}}{R^{|k|+1}}\left(bdN^{d+t}
R^d-|k|\right),$$
qui s'annule exactement en
$R=R_{min}=\left(\frac{|k|}{bdN^{d+t}}\right)^{\frac{1}{d}}$.
\bigskip

Par suite, il y a un seul minimum, et on a :
$$\|c_k\|\leq\frac{Ae^{uN^v}e^{bN^{d+t}\frac{|k|}{bdN^{d+t}}}}{\left(\frac{|k|}{bdN^{d+t}}\right)^{\frac{|k|}{d}}}=Ae^{uN^v}\left(\frac{ebdN^{d+t}}{|k|}\right)^{\frac{|k|}{d}},$$
ce qui prouve le lemme.
\end{proof}
\bigskip

\bigskip

Nous pouvons commencer la preuve du th{\'e}or{\`e}me. 
\bigskip

On a :
$\forall\,\zeta\in\mathbb{R}^N,\,|\zeta_1|,\ldots,|\zeta_N|\leq BN^r$,
$$f(x,\zeta)=\sum_{|k|\geq 0}c_k(x)\zeta^k=\sum_{|k|\leq K}c_k(x)\zeta^k+\sum_{|k|>K}c_k(x)\zeta^k=P_K(x,\zeta)+R_K(x,\zeta),$$
o{\`u} $K$ est un entier $>0$. Or, d'apr{\`e}s le lemme :

\begin{eqnarray*}
\|R_K(\cdot,\zeta)\| & \leq  & 
\sum_{|k|>K}\|c_k\|\,|\zeta_1|^{k_1}\ldots|\zeta_N|^{k_N}\\
& \leq & Ae^{uN^v}\sum_{|k|>K}
\left(\frac{ebdB^dN^{d(r+1)+t}}{|k|}\right)^{\frac{|k|}{d}}\\
& = & Ae^{uN^v}\sum_{n>K}
\left(\frac{ebdB^dN^{d(r+1)+t}}{n}\right)^{\frac{n}{d}}
card\{\,k_1+\cdots+k_N=n\}\\
& = & Ae^{uN^v}\sum_{n>K}
\left(\frac{ebdB^dN^{d(r+1)+t}}{n}\right)^{\frac{n}{d}}
\frac{(n+N-1)!}{n!\,(N-1)!},
\end{eqnarray*}

la derni{\`e}re {\'e}galit{\'e} provenant du

\begin{lemma}
$$card\{k\in\mathbb{N}^N,\;k_1+\cdots+k_N=n\}=\frac{(n+N-1)!}{(N-1)!\,n!}\,.$$
\end{lemma}

\begin{proof}
Le membre de gauche est le coefficient en $X^n$ de la s{\'e}rie formelle
$$\sum_{k_1,\ldots,k_N\geq 0}X_1^{k_1}\ldots X_N^{k_N},$$ apr{\`e}s l'{\'e}valuation
$X_1=\ldots=X_N=X$. Or, cette s{\'e}rie vaut :
$$\prod_{j=1}^N\left(\sum_{k_j\geq 0}X_j^{k_j}\right)=\prod_{j=1}^N\frac{1}{1-X_j},$$
qui apr{\`e}s {\'e}valutation donne $\frac{1}{(1-X)^N}$. On a d'autre part :
\begin{eqnarray*}
\frac{1}{(1-X)^N} & = &
\frac{1}{(N-1)!}\frac{d^{N-1}}{dX^{N-1}}\left(\frac{1}{1-X}\right)\\
& = & \frac{1}{(N-1)!}\sum_{k\geq 0}k(k-1)...(k-N+2)X^{k-N+1},
\end{eqnarray*}
dont le
coefficient en $X^n$ correspond {\`a} $k=n+N-1$, ce qui donne finalement
$$\frac{(n+N-1)\ldots(n+1)}{(N-1)!},$$ et prouve le lemme.
\end{proof}
\bigskip

Ensuite, en utilisant les majorations suivantes
$$\left(\frac{n}{e}\right)^n\sqrt{2\pi n}\leq n!\leq\left(\frac{n}{e}\right)^n2\pi n,$$
on a, pour tout $n\geq N$,
\begin{eqnarray*}
\frac{(n+N-1)!}{n!\,(N-1)!} & \leq &
\left(\frac{n+N-1}{n}\right)^n\left(\frac{n+N-1}{N-1}\right)^{N-1}\frac{n+N-1}{\sqrt{n(N-1)}}\\
& \leq &
e^{N-1}\,\left(\frac{2n}{N-1}\right)^{N-1}\,2\left(\frac{n}{N-1}\right)^{\frac{1}{2}}\\
& \leq & \left(\frac{2en}{N-1}\right)^N\,.
\end{eqnarray*}

On en d{\'e}duit, {\`a} condition de prendre $K\geq N$,
\begin{eqnarray*}
\|R_K(\cdot,\zeta)\| & \leq & Ae^{uN^v}\sum_{n\geq K+1}\left(\frac{ebdB^dN^{d(r+1)+t}}{n}\right)^{\frac{n}{d}}\left(\frac{2en}{N-1}\right)^N\\
& = & Ae^{uN^v}\sum_{n\geq
  K+1}\left(\frac{ebdB^dN^{d(r+1)+t}}{n}\left(\frac{2en}{N-1}\right)^{\frac{dN}{n}}\right)^{\frac{n}{d}}.
\end{eqnarray*}
Or,
\begin{eqnarray*}
\left(\frac{2en}{N-1}\right)^{\frac{dN}{n}} & = & \left(\left(\frac{2en}{N-1}\right)^{\frac{1}{\frac{2en}{N-1}}}\right)^{\frac{2edN}{N-1}}\\
& \leq & 2^{\frac{2edN}{N-1}}\\
& \leq & 2^{4ed},
\end{eqnarray*}
et donc $$\|R_K(\cdot,\zeta)\|\leq Ae^{uN^v}\sum_{n\geq K+1}\left(\frac{ebd\left(2^{4e}B\right)^dN^{d(r+1)+t}}{n}\right)^{\frac{n}{d}}.$$
Par ailleurs, pour tout $$n\geq K\geq
ebd2^{(4e+1)d}B^dN^{d(r+1)+t}\;(\,\geq N),$$
on a
$$\left(\frac{ebd\left(2^{4e}B\right)^dN^{d(r+1)+t}}{n}\right)^{\frac{1}{d}}\leq\frac{1}{2},$$
ce qui donne $$\|R_K(\cdot,\zeta)\|\leq\frac{Ae^{uN^v}}{2^K}.$$
\bigskip

Remarquons maintenant que, pour $K\geq 4\left(\frac{l}{s}\right)^2$,
\begin{eqnarray*}
K-\frac{l}{s}\log(\log K) & = & K\left(1-\frac{l}{s}
\frac{\log(\log K)}{K}\right)\\
& \geq & K\left(1-\frac{l}{s}\frac{1}{\sqrt{K}}\right)\\
& \geq & \frac{K}{2},
\end{eqnarray*}
la premi{\`e}re in{\'e}galit{\'e} provenant du fait que, pour tout $x\geq 4$,
$\log(\log x)\leq\sqrt{x}.$
\bigskip

On choisit finalement $$K\geq
ebud2^{(4e+1)d}B^d\left(1+\frac{l}{s}\right)^2\left(\log\frac{4A}{C}\right)\,N^{d(r+1)+v+t},$$ 
o{\`u} $C$ est une constante $\leq 1$ que l'on pr{\'e}cisera, et on montre le

\begin{lemma}
Sous les hypoth{\`e}ses pr{\'e}c{\'e}dentes, on a :
$$\|R_K(\cdot,\zeta)\|\leq\frac{C}{4(N\log K)^{\frac{l}{s}}}.$$
\end{lemma}

\begin{proof}
En effet, on v{\'e}rifie d'abord que $$K\geq
ebd2^{(4e+1)d}B^dN^{d(r+1)+t},$$ puisque
$\left(1+\frac{l}{s}\right)^2\log\frac{4A}{C}\,uN^v\geq 1$.

Ensuite, $$K\geq e2^{4e+1}\left(1+\frac{l}{s}\right)^2\geq4\left(\frac{l}{s}\right)^2.$$

Il reste {\`a} s'assurer que
$$\|R_K(\cdot,\zeta)\|\leq\frac{Ae^{uN^v}}{2^K}\leq\frac{C}{4(N\log K)^{\frac{l}{s}}},$$
soit $$\frac{2^K}{(\log
  K)^{\frac{l}{s}}}\geq\frac{4Ae^{uN^v}N^{\frac{l}{s}}}{C},$$
ou encore, puisque $K-\frac{l}{s}\log\log
K\geq\frac{K}{2}$, $$K\geq\ 2\log\frac{4A}{C}+2\frac{l}{s}\log
N+2uN^v\log e.$$

Or,
\begin{eqnarray*}
K & \geq & 2\log\frac{4A}{C}\,\left(1+\frac{l}{s}\right)^2N\,2euN^v\\
& \geq & 2\log\frac{4A}{C}+\left(1+\frac{l}{s}\right)^2N+2euN^v,
\end{eqnarray*}
car $2\log\frac{4A}{C}$, $\left(1+\frac{l}{s}\right)^2N$ et $2euN^v$
sont $\geq 2$.
D'autre part, $\left(1+\frac{l}{s}\right)^2\geq 2\frac{l}{s}$,
$N\geq\log N$, ce qui donne 
$$K\geq 2\log\frac{4A}{C}+2\frac{l}{s}\log N+2uN^v\log e,$$ et prouve
le lemme. 

\end{proof}

On a finalement :$$f(\cdot,\zeta)=P_K(\cdot,\zeta)+R_K(\cdot,\zeta),$$
o{\`u} $P_K$ est un 
polyn{\^o}me en $\zeta$, {\`a} coefficients dans l'espace $C(I^s)$
(resp. $L^1(I^s))$, {\`a} $N$ variables et de degr{\'e} $K$, et
$\|R_K(\cdot,\zeta)\|\leq\frac{C}{4(N\log K)^{\frac{l}{s}}}$, o{\`u} la
norme d{\'e}signe la norme uniforme (resp. $L^1$).
\bigskip

D'apr{\`e}s le th{\'e}or{\`e}me de Vitushkin (dont on redonnera une preuve dans la
partie~\ref{vitwar}), il existe $h\in\Lambda_{l,s}$, tel
que, $\forall\,\zeta\in\mathbb{R}^N,$
$$\|h-P_K(\cdot,\zeta)\|\geq\frac{C}{(N\log K)^{\frac{l}{s}}},$$
o{\`u} $C=C(l,s)$ est une constante variant selon le cas continu ou
$L^1$.
Il en r{\'e}sulte que
\begin{eqnarray*}
\|h-f(\cdot,\zeta)\| & \geq &
\|h-P_K(\cdot,\zeta)\|-\|R_K(\cdot,\zeta)\|\\
& \geq & \frac{3C}{4(N\log K)^{\frac{l}{s}}}.
\end{eqnarray*}

Pour terminer la preuve du th{\'e}or{\`e}me, il reste {\`a} obtenir
$\frac{C'}{(N\log N)^{\frac{l}{s}}}$.
Pour cela, on choisit le plus petit entier $K$ possible, soit
\begin{eqnarray*}
K & \leq & 1+ebud2^{(4e+1)d}B^d\left(1+\frac{l}{s}\right)^2\log\frac{4A}{C}\,N^{d(r+1)+v+t}\\
& \leq & 2ebud2^{(4e+1)d}B^d\left(1+\frac{l}{s}\right)^2\log\frac{4A}{C}\,N^{d(r+1)+v+t},
\end{eqnarray*}
donc
\begin{eqnarray*}
\log K & \leq & (d(r+1)+v+t)\log
N+\log\left(2ebud2^{(4e+1)d}B^d\left(1+\frac{l}{s}\right)^2\log\frac{4A}{C}\right)\\ 
& \leq & (d(r+1)+v+t)\log
N.\log\left(2ebud2^{(4e+1)d}B^d\left(1+\frac{l}{s}\right)^2\log\frac{4A}{C}\right),
\end{eqnarray*}
car les deux termes de la somme sont $\geq2$,
ce qui donne finalement
$$\|h-f(\cdot,\zeta)\|\geq\frac{C'}{(N\log N)^{\frac{l}{s}}},$$
avec
$$C'=\frac{3C}{4\left[(d(r+1)+v+t)\log\left(2ebud2^{(4e+1)d}B^d\left(1+\frac{l}{s}\right)^2\log\frac{4A}{C}\right)\right]^{\frac{l}{s}}},$$
et ach{\`e}ve la preuve du th{\'e}or{\`e}me.

\end{proof}

\begin{remark}\label{remth1}

La fonction $h$ qui intervient pour obtenir la minoration dans le cas
analytique, est la m{\^e}me que dans le cas polynomial. On s'en servira
dans la preuve des corollaires suivants.

D'autre part, comme on le constatera dans la partie suivante, la
constante $C$ ne d{\'e}pend que de $l$ et $s$, ainsi 
$C'$ va d{\'e}pendre de $l,s,A,u,v,b,t,d,B,r$ (mais pas de $N$). On verra
en particulier des exemples d'estimation de $C$, et donc de
calcul explicite de $C'$ (cf corollaire~\ref{constantes}).

\end{remark}

\subsection{Quelques r{\'e}sultats d{\'e}riv{\'e}s}

On peut g{\'e}n{\'e}raliser le th{\'e}or{\`e}me~\ref{th1} de la fa{\c c}on suivante :

\begin{theorem}\label{th1bis}

Soit de m{\^e}me $f(x,\zeta,w)$ une fonction qui v{\'e}rifie les m{\^e}mes
conditions que dans le th{\'e}or{\`e}me~\ref{th1}, et qui est en plus
holomorphe par rapport aux variables $w_1,\ldots,w_m$ sur le domaine
$$\{\Re e\,z>0\}^m\,,$$
et de type exponentiel : 
$\forall\,(\zeta,w)\in\mathbb{C}^N\times\left\{\Re e\,z>0\right\}^m,$ 
$$\|f(\cdot,\zeta,w)\|\leq Ae^{uN^v}e^{bN^t\left(\|\zeta\|_1^d
+\|w\|_1^d\right)}\,.$$

Alors l'approximation de $\Lambda_{l,s}$ par la famille
$\{f(\cdot,\zeta,w)\}$, o{\`u} les $\zeta_j,\,w_i$ sont r{\'e}els,
$\zeta_j=O(N^r)$ et 
$$|w_i-a_i(N)|\leq\frac{1}{2}a_i(N),\;i=1,\ldots,m,$$
avec $a_i(N)\geq1$ polynomial en $N$, est au mieux de l'ordre de
$$\frac{1}{(N\log N)^{\frac{l}{s}}}\,.$$

\end{theorem}

\begin{proof}

La preuve est semblable {\`a} celle du th{\'e}or{\`e}me~\ref{th1} : on {\'e}crit de
m{\^e}me le d{\'e}veloppement de Taylor de $f$ au 
point $(0,\ldots,0,a_1(N),\ldots,a_m(N))$, soit
$$f(\cdot,\zeta,w)=\left(\sum_{|k|\leq K,\,|l|\geq0}+
\sum_{|k|>K,\,|l|\geq0}\right)
\left(c_{k,l}(a)\zeta^k(w-a)^l\right),$$
o{\`u} les coefficients $c_{k,l}$ sont major{\'e}s {\`a} l'aide de la formule de
Cauchy : $\forall\,R,\;\forall\,r_i<a$,
\begin{eqnarray*}
|c_{k,l}| & \leq & \frac{1}{(2\pi)^{N+m}}\int_{|\zeta'_j|=R,\,
|w'_i-a_i|=r_i}\frac{|f(\cdot,\zeta',w')|}{|\zeta'^{k+1}|\,
|(w'-a)^{l+1}|}d\zeta'dw'\\
& \leq &
Ae^{uN^v}\,\frac{e^{bN^{d+t}R^d}}{R^{|k|}}\,
\frac{e^{bN^t\|a+r\|_1^d}}{|r^l|}\,.
\end{eqnarray*}
Puisque $|w_i-a_i|\leq\frac{a_i}{2}$, la deuxi{\`e}me somme ainsi peut
{\^e}tre major{\'e}e :
\begin{eqnarray*}
\sum_{|k|>K,\,|l|\geq0}|c_{k,l}|\,\left|\zeta^k\right|\,\left|(w-a)^l\right|
& \leq &
Ae^{uN^v}\sum_{|k|>K}\frac{e^{bN^{d+t}}R^d}{R^{|k|}}
\left(BN^r\right)^{|k|}\sum_{|l|\geq0}\frac{e^{bN^t2^d\|a\|_1^d}}{2^{|l|}}\\
& \leq &
A2^me^{uN^v+bN^t2^d\|a\|_1^d}\sum_{n>K}\left(\frac{ebdB^dN^{d(r+1)+t}}{n}
\right)^{\frac{n}{d}}\frac{(n+N-1)!}{n!\,(N-1)!}\,,
\end{eqnarray*}
ce qui, pour $K$ de l'ordre de $N^{\alpha}$ (cf preuve du
th{\'e}or{\`e}me~\ref{th1}), est major{\'e} par 
$$\frac{C}{8(N\log K)^{\frac{l}{s}}}\,.$$

Pour l'autres somme, on r{\'e}{\'e}crit de m{\^e}me
$$\sum_{|k|\leq K,\,|l|\geq0}=\sum_{|k|\leq K}\sum_{|l|\leq p(k)}\;+
\sum_{|k|\leq K}\sum_{|l|>p(k)}\,,$$
o{\`u} pour chaque $n=|k|$, on choisit $p(k)=p(n)=p(K),$ ainsi :
\begin{eqnarray*}
\sum_{|l|>p(K)}\frac{1}{2^{|l|}} & = & \sum_{p>p(K)}\frac{1}{2^p}
\frac{(p+m-1)!}{p!\,(m-1)!},\\
& \leq & \sum_{p>p(K)}\frac{(p+1)^{m-1}}{2^p}\\
& \leq & \left(\frac{3}{4}\right)^{p(K)},
\end{eqnarray*}
pour $p(K)\geq p(K,m)$, ce qui donne
\begin{eqnarray*}
\sum_{|k|\leq K}\sum_{|l|>p(k)}|c_{k,l}|\left|\zeta^k\right|
\left|(w-a)^l\right| & \leq &
Ae^{uN^v+bN^t2^d\|a\|_1^d}\sum_{n=0}^K\left(\frac{ebdB^dN^{d(r+1)+t}}
{n}\right)^{\frac{n}{d}}\frac{(n+N-1)!}{n!\,(N-1)!}
\sum_{|l|>p(K)}\frac{1}{2^{|l|}}\\
& \leq &
Ae^{uN^v+bN^t2^d\|a\|_1^d}\sum_{n=0}^K\left(\frac{(ebd)^{\frac{1}{d}}B
N^{r+1+\frac{t}{d}}}{n^{\frac{1}{d}}}\right)^nN^n\left(\frac{3}{4}\right)^{p(K)}\\
& \leq & Ae^{uN^v+bN^t2^d\|a\|_1^d}\,(K+1)\left((ebd)^{\frac{1}{d}}
BN^{r+t+2}\right)^K\left(\frac{3}{4}\right)^{p(K)}\\
& \leq & \frac{1}{2^K}\,,
\end{eqnarray*}
si $p(K)\geq\alpha K$, ce qui donnera encore une majoration par
$$\frac{C}{8(N\log K)^{\frac{l}{s}}}\,.$$

Pour finir, le somme restante est
$$\sum_{|k|\leq K}\sum_{|l|\leq p(K)}c_{k,l}\zeta^k(w-a)^l\,,$$
qui est un polyn{\^o}me en $(\zeta,w)$ (donc {\`a}
$N+m$ variables), de degr{\'e} total (au plus) $K+p(K)\leq\gamma K$, ce
qui donne, par le th{\'e}or{\`e}me~\ref{vitushkin} de Vitushkin :
$\exists\,h\in\Lambda_{l,s},\;\forall\,\zeta,\,w$,
\begin{eqnarray*}
\|f(\cdot,\zeta,w)-h\| & \geq & \frac{3C'}{4((N+m)\log\gamma K)
^{\frac{l}{s}}}\\
& \geq & \frac{C''}{(N\log N)^{\frac{l}{s}}}\,,
\end{eqnarray*}
ce qui prouve l'assertion.

\end{proof}

On en d{\'e}duit respectivement les corollaires suivants des
th{\'e}or{\`e}mes~\ref{th1bis} et~\ref{th1}, qui nous serviront dans la
partie~\ref{pbinverse}.

\begin{corollary}\label{inverse}

Soit $\psi(x,\zeta,w)$ une fonction d{\'e}finie sur
$[0,1]\times\mathbb{C}^N\times\left\{\Re e\,z>0\right\}^m$, de classe
$C^2$ par rapport {\`a} $x$, et telle 
que, pour tout $x\in[0,1]$, les fonctions $\psi(x,\zeta,w)$, 
$\frac{\partial\psi}{\partial x}(x,\zeta,w)$ et
$\frac{\partial^2\psi}{\partial x^2}(x,\zeta,w)$ v{\'e}rifient les
conditions du th{\'e}or{\`e}me~\ref{th1bis} pour le cas continu. Soit d'autre
part $k(x,\zeta,w)$ continue par rapport {\`a} $x$, analytique de type
exponentiel par rapport {\`a} $(\zeta,w)$, dont la restriction sur
$\mathbb{R}^N\times\left(\mathbb{R}^+\right)^m$ est domin{\'e}e polynomialement.

On suppose {\'e}galement que,
$\forall\,(\zeta,w)\in\mathbb{R}^N\times\left(\mathbb{R}^+\right)^m$,
$\zeta_j=O(N^r),\;|\zeta_i-a_i(N)|\leq\frac{a_i(N)}{2}$, 
$$\frac{1}{\psi(0,\zeta,w)}=O\left(e^{\alpha N^{\beta}}\right),$$ et
$$\frac{\partial\psi}{\partial x}(0,\zeta,w)=0.$$
Soit finalement une constante $b(N)>0$, telle que $b(N)$ et
$\frac{1}{b(N)}$ soient polynomiales en $N$.
\bigskip

Alors l'approximation de $\Lambda_l([0,1])$ par la famille 
$$\left\{\left(x\in[0,1]\mapsto\frac{1}{b(N)}\left(k(x,\zeta,w)+
\frac{\partial}{\partial x}
\left(\frac{\frac{\partial\psi}{\partial x}}{\psi}\right)(x,\zeta,w)\right)
\right),\,\zeta_j=O(N^r),\,|w_i-a_i(N)|\leq\frac{a_i(N)}{2}\right\},$$
au sens uniforme sur $[0,1]$, lorsque $N\rightarrow+\infty$, ne peut
    pas {\^e}tre meilleure que $$\frac{\widetilde{C}}{N^l(\ln N)^l}.$$
($\widetilde{C}$ d{\'e}pend des param{\`e}tres qui apparaissent dans la
    constante $C'$ du th{\'e}or{\`e}me~\ref{th1bis}, ainsi que de $\alpha$, $\beta$)

En outre, une fonction $h$ qui r{\'e}alise la minoration peut {\^e}tre prise dans
$\Lambda_l\cap C^{[l]}$, et {\^e}tre identiquement nulle sur
$\left[0,\frac{1}{20}\right]$.

\end{corollary}

\begin{proof}
On peut supposer que, $\forall\,\zeta,\,w$, $\psi(0,\zeta,w)>0$ : en effet,
la condition sur $\psi(0,\zeta,w)$ montre qu'elle ne s'annule pas sur
l'ensemble (connexe) $\mathbb{R}^N\times\left(\mathbb{R}^+\right)^m$,
donc reste de signe constant, et peut {\^e}tre remplac{\'e}e par $-\psi$ (ce
qui ne change pas la famille).

On se donne $\psi,\,k$ et on pose
$$\widetilde{\psi}(x,\zeta,w)=e^{\widetilde{k}(x,\zeta,w)}
\psi(x,\zeta,w)\,,$$ 
o{\`u} $\widetilde{k}(x,\zeta,w)=\int_0^xdt\int_0^tk(s,\zeta,w)ds$. Alors le
th{\'e}or{\`e}me~\ref{th1bis} est encore vrai (m{\^e}me si $\widetilde{\psi}$ peut
ne plus {\^e}tre de type exponentiel) : en effet, on {\'e}crit
$$\psi(x,\zeta,w)=P_K(x,\zeta,w)+R_K(x,\zeta,w),$$
et puisque 
$e^{\widetilde{k}(x,\zeta,w)}=O\left(e^{\alpha N^{\beta}}\right)$, quitte {\`a}
agrandir $K$ (tout en le gardant polynomial en $N$), on a encore
$$\left\|e^{\widetilde{k}(\cdot,\zeta,w)}R_K(\cdot,\zeta,w)\right\|_{\infty}
\leq\frac{C_{\infty}(l,1)}{4((N+m)\log K)^l}\,.$$
D'autre part, il existe $h\in\Lambda_l([0,1])$ pour $P_K$, tel que 
$\forall\,(\zeta,w),\;\exists\,x(\zeta,w)$,
$$\frac{C'_{\infty}(l,1)}{((N+m)\log K)^l}\leq|h(x)|\;
\text{ et }h(x)P_K(x,\zeta,w)\leq0,$$
ce qui donne encore
$$\left\|h-e^{\widetilde{k}(\cdot,\zeta,w)}P_K(\cdot,\zeta,w)
\right\|_{\infty}\geq\frac{C'_{\infty}(l,1)}{((N+m)\log K)^l}\,.$$
\bigskip

On peut maintenant prouver le corollaire. Posons
$$\widetilde{\chi}(x,\zeta_1,\ldots,\zeta_N,\zeta_{N+1},w)=e^{\zeta_{N+1}}
\left(\frac{\partial^2\widetilde{\psi}}{\partial^2x}\widetilde{\psi}-
\left(\frac{\partial\widetilde{\psi}}
{\partial x}\right)^2\right)(x,\zeta,w)\,.$$
On voit alors que
$$\widetilde{\chi}(\cdot,\zeta,w)
=e^{2\widetilde{k}(\cdot,\zeta,w)}\chi(\cdot,\zeta,w),$$
o{\`u} $\chi$ v{\'e}rifie encore les conditions de l'{\'e}nonc{\'e} (contrairement {\`a}
$\widetilde{\chi}$ qui n'est plus n{\'e}cessairement de type
exponentiel). En particulier, celles du th{\'e}or{\`e}me~\ref{th1bis} sont
satisfaites, et ce qui pr{\'e}c{\`e}de prouve qu'on a l'existence de
$h\in\Lambda_l([0,1])$ qui v{\'e}rifie, pour tous $\zeta_j=O(N^r)$,
$j=1,\ldots,N+1,\;|w_i-a_i(N)|\leq\frac{a_i(N)}{2},\,i=1,\ldots,m$,
$$\|h-\widetilde{\chi}(\cdot,\zeta,\zeta_{N+1},w)\|_{\infty}\geq
\frac{C'_{\infty}}{((N+m+1)\log(N+m+1))^l},\,$$
et m{\^e}me plus pr{\'e}cis{\'e}ment, si
$\chi=P_K+R_K$
($K=O(N^{\gamma})$), pour $x=x(\zeta,\zeta_{N+1},w)$ : 
$$(\ast)
\begin{cases}
\frac{C'_{\infty}(l,1)}{((N+m+1)\log K)^l}\leq|h(x)|\leq\|h\|_{\infty}
\leq\frac{2^lC'_{\infty}(l,1)}{((N+m+1)\log K)^l},\\
h(x)e^{2\widetilde{k}(x,\zeta,w)}P_K(x,\zeta,\zeta_{N+1},w)\leq0,\\
\left\|e^{2\widetilde{k}(\cdot,\zeta,w)}R_K(\cdot,\zeta,\zeta_{N+1},w)
\right\|_{\infty}\leq\frac{C'_{\infty}(l,1)}{4((N+m+1)\log K)^l}\,.
\end{cases}
$$

En outre, $h\in C^{[l]}([0,1])$ et est identiquement nulle sur
$\left[0,\frac{1}{20}\right]$ 
(la troisi{\`e}me ligne de $(\ast)$  provient de
la preuve du th{\'e}or{\`e}me~\ref{th1bis} ; les autres hypoth{\`e}ses
r{\'e}sultent de la 
construction de $h$, elles seront justifi{\'e}es dans la
partie~\ref{vitwar} comme le pr{\'e}cisera le corollaire~\ref{1min} qui
suit la preuve du th{\'e}or{\`e}me~\ref{vitushkin} de Vitushkin).

D'autre part, soit $\mu\in[0,1]$
un r{\'e}el. On a
$$h(x)\,\mu\,e^{2\widetilde{k}(x,\zeta,w)}P_K(x,\zeta,\zeta_{N+1},w)
\leq0,$$
et
$$\left\|\mu\,e^{\widetilde{k}(\cdot,\zeta,w)}
R_K(\cdot,\zeta,\zeta_{N+1},w)\right\|_{\infty}
\leq\frac{C'_{\infty}(l,1)}{4((N+m+1)\log K)^l},$$
ce qui donne 
$$\left|h(x)-\mu\widetilde{\chi}(x,\zeta,\zeta_{N+1},w)\right|
\geq\frac{3C'_{\infty}(l,1)}{4((N+m+1)\log K)^l}\,.$$
Par ailleurs, $\psi$ {\'e}tant de type exponentiel, on a en particulier,
pour tous $\zeta_j,\,w_i=O(N^r)$,
$$0<\widetilde{\psi}(0,\zeta,w)\leq Ae^{uN^v+bN^t
\left(((N+1)BN^r)^d+(2ma_i(N))^d\right)}=
O\left(e^{\alpha' N^{\beta'}}\right),$$
et il en est de m{\^e}me pour $\frac{1}{\widetilde{\psi}(0,\zeta,w)}$ et
$e^{b(N)}$ par hypoth{\`e}se, ce qui montre que 
$$\ln\left(\frac{e^{b(N)}}{b(N)\widetilde{\psi}(0,\zeta,w)^2}\right)
=O\left(N^{r'}\right),$$ 
que l'on peut alors choisir comme valeur du param{\`e}tre $\zeta_{N+1}$,
ce qui donnera (quitte {\`a} agrandir $K=O(N^{\gamma})$ ), pour tout
$\mu\in[0,1]$ : 
$$\left|h(x)-\frac{\mu e^{b(N)}}{b(N)}
\frac{\frac{\partial^2\widetilde{\psi}}
{\partial^2x}(x,\zeta,w)\,\widetilde{\psi}(x,\zeta,w)
-\left(\frac{\partial\widetilde{\psi}}{\partial x}\right)^2(x,\zeta,w)}
{\left(\widetilde{\psi}(0,\zeta,w)\right)^2}\right|\geq
\frac{3C'_{\infty}(l,1)}{4((N+m+1)\log K)^l}.$$
\bigskip

Supposons alors
$e^{-b(N)}\left(\frac{\widetilde{\psi}(0,\zeta,w)}{\widetilde{\psi}
(x,\zeta,w)}\right)^2\leq1$.
On peut donc, dans l'in{\'e}galit{\'e}, remplacer $\mu$ par
$e^{-b(N)}\left(\frac{\widetilde{\psi}(0,\zeta,w)}{\widetilde{\psi}
(x,\zeta,w)}\right)^2$, 
pour obtenir : 
\begin{eqnarray*}
\left|h(x)-\frac{1}{b(N)}\frac{\frac{\partial^2\widetilde{\psi}}
{\partial^2x}(x,\zeta,w)\,\widetilde{\psi}(x,\zeta,w)
-\left(\frac{\partial\widetilde{\psi}}
{\partial x}\right)^2(x,\zeta,w)}{\left(\widetilde{\psi}
(x,\zeta,w)\right)^2}\right| 
& \geq & \frac{3C'_{\infty}(l,1)}{4((N+m+1)\log K)^l}\\
& \geq & \frac{C'}{(N\log N)^l}\,.
\end{eqnarray*}

Dans l'autre cas, on a
$\left|\frac{\widetilde{\psi}(x,\zeta,w)}{\widetilde{\psi}
(0,\zeta,w)}\right|<e^{-\frac{b(N)}{2}}$.
Quitte {\`a} rapprocher $x$ de $0$, on peut supposer que la fonction
$t\mapsto\widetilde{\psi}(t,\zeta,w)$ 
ne s'annule pas sur $[0,x]$. En effet, car sinon, soit $x_0$ le
premier z{\'e}ro ($>0$) de $\widetilde{\psi}(\cdot,\zeta,w)$. Comme
$\lim_{t\rightarrow x_0}\widetilde{\psi}(t,\zeta,w)=0$, il suffit de choisir
$x$ suffisamment proche de $x_0$, pour que
$$0<\left|\frac{\widetilde{\psi}(x,\zeta,w)}{\widetilde{\psi}
(0,\zeta,w)}\right|<e^{-\frac{b(N)}{2}}.$$
En particulier, $\widetilde{\psi}(x,\zeta,w)$ et
$\widetilde{\psi}(0,\zeta,w)$ sont de m{\^e}me signe, la fonction
$$t\in[0,x]\mapsto\ln\frac{\widetilde{\psi}(t,\zeta,w)}{\widetilde{\psi}
(0,\zeta,w)},$$
est donc bien d{\'e}finie et r{\'e}guli{\`e}re. Puisque
$$\ln\frac{\widetilde{\psi}(x,\zeta,w)}{\widetilde{\psi}(0,\zeta,w)}
<-\frac{b(N)}{2},$$
alors il existe $x_1\in[0,x]$, tel que
$$\left|\frac{\partial}{\partial t}\ln\frac{\widetilde{\psi}}
{\widetilde{\psi}(0,\zeta,w)}
(x_1,\zeta,w)\right|=\left|\frac{\frac{\partial\widetilde{\psi}}
{\partial t}(x_1,\zeta,w)}{\widetilde{\psi}(x_1,\zeta,w)}\right|>
\frac{b(N)}{2}.$$
En effet, car sinon on aurait
$$\left|\ln\frac{\widetilde{\psi}(x,\zeta,w)}{\widetilde{\psi}
(0,\zeta,w)}\right|=\left|\int_0^x\frac{\frac{\partial\widetilde{\psi}}
{\partial t}(t,\zeta,w)}{\widetilde{\psi}(t,\zeta,w)}dt\right|
\leq\frac{b(N)}{2},$$
ce qui est impossible.

De m{\^e}me, puisque par hypoth{\`e}se, 
$\frac{\partial\widetilde{\psi}}{\partial t}(0,\zeta,w)=0$, il existe
$x_2\in[0,x_1]$, tel que 
$$\left|\frac{\frac{\partial^2\widetilde{\psi}}{\partial^2t}(x_2,\zeta,w)
\widetilde{\psi}(x_2,\zeta,w)-\left(\frac{\partial\widetilde{\psi}}
{\partial t}\right)^2(x_2,\zeta,w)}{\left(\widetilde{\psi}(x_2,\zeta,w)
\right)^2}\right|>\frac{b(N)}{2}\,.$$
On en d{\'e}duit 
$$\left|h(x_2)-\frac{1}{b(N)}\frac{\frac{\partial^2\widetilde{\psi}}
{\partial^2t}(x_2,\zeta,w)\widetilde{\psi}(x_2,\zeta,w)
-\left(\frac{\partial\widetilde{\psi}}{\partial t}\right)^2(x_2,\zeta,w))}
{\left(\widetilde{\psi}(x_2,\zeta,w)\right)^2}\right|\geq
\frac{1}{2}-\|h\|_{\infty}.$$
Or, d'apr{\`e}s $(\ast)$, quitte {\`a} r{\'e}duire $C'_{\infty}(l,1)$, on a
{\'e}galement $\|h\|_{\infty}\leq\frac{1}{4}$, ce qui donne
$$\left|h(x_2)-\frac{1}{b(N)}\frac{\frac{\partial^2\widetilde{\psi}}
{\partial^2t}(x_2,\zeta,w)\,\widetilde{\psi}(x_2,\zeta,w)-
\left(\frac{\partial\widetilde{\psi}}{\partial t}\right)^2(x_2,\zeta,w)}
{\left(\widetilde{\psi}(x_2,\zeta,w)\right)^2}\right|\geq\frac{1}{4}
\geq\frac{\widetilde{C}}{(N\log N)^l},$$
quitte {\`a} r{\'e}duire $\widetilde{C}$, ce qui ach{\`e}ve le second cas.
\bigskip

Puisqu'on a {\'e}galement
$$\frac{\partial^2}{\partial x^2}\ln\frac{\widetilde{\psi}(x,\zeta,w)}
{\widetilde{\psi}(0,\zeta,w)}=k(x,\zeta,w)+\frac{\partial^2}{\partial x^2}
\ln\frac{\psi(x,\zeta,w)}{\psi(0,\zeta,w)}\,,$$
on peut enfin conclure l'existence de $h\in\Lambda_l([0,1])$,
tel que, pour tous $\zeta_j=O(N^r),\,j=1,\ldots,N$ et
$|w_i-a_i(N)|\leq\frac{a_i(N)}{2},\,i=1,\ldots,m$, 
$$\left\|h-\frac{1}{b(N)}\left(k(\cdot,\zeta,w)+\frac{\frac{\partial^2\psi}
{\partial^2x}(\cdot,\zeta,w)\psi(\cdot,\zeta,w)-\left(\frac{\partial\psi}
{\partial x}\right)^2(\cdot,\zeta,w)}{\left(\psi(\cdot,\zeta,w)\right)^2}
\right)\right\|_{\infty}\geq\frac{\widetilde{C}}{(N\ln N)^l}\,.$$

\end{proof}

On prouve {\'e}galement le

\begin{corollary}\label{inverse2}

On consid{\`e}re cette fois $\psi(x,\zeta)$ sur $[0,1]\times\mathbb{C}^N$, de
classe $C^1$ par rapport {\`a} $x$, et tel 
que $\psi(x,\zeta)$ et 
$\frac{\partial\psi}{\partial x}(x,\zeta)$ v{\'e}rifient les
conditions du th{\'e}or{\`e}me~\ref{th1}, ainsi que
$$\frac{1}{\psi(0,\zeta)}=O\left(e^{\alpha N^{\beta}}\right).$$

Alors l'approximation du compact
$$\int_0^{\cdot}\Lambda_l([0,1])=\left\{\left(x\mapsto\int_0^xh(t)dt\right),
\,\,h\in\Lambda_l\right\}$$
par la famille 
$$\left\{\left(x\in[0,1]\mapsto\frac{1}{b(N)}\left(\frac{1}{\psi}
\frac{\partial\psi}{\partial x}\right)(x,\zeta)\right),\,
\zeta_j=O(N^r),\,\forall\,j=1,\ldots,N\right\},$$
sur $[0,1]$, avec $b(N)>0$ et $\frac{1}{b(N)}$ polynomiaux en $N$, ne
peut pas {\^e}tre meilleure que 
$$\frac{\widetilde{C}}{N^{l+1}(\ln N)^{l+1}}.$$

Ici aussi, $h$ a les m{\^e}mes propri{\'e}t{\'e}s que dans le corollaire
pr{\'e}c{\'e}dent.

\end{corollary}

\begin{proof}

On raisonne de m{\^e}me en posant cette fois
$$\chi(x,\zeta_1,\ldots,\zeta_N,\zeta_{N+1})=e^{\zeta_{N+1}}
\frac{\partial\psi}{\partial x}(x,\zeta_1,\ldots,\zeta_N).$$
On a l'existence de $h_1\in\Lambda_{l+1}([0,1])$, tel que
$$\|h_1-\chi(\cdot,\zeta_1,\ldots,\zeta_{N+1})\|_{\infty}\geq
\frac{C'_{\infty}(l+1,1)}{((N+1)\log(N+1))^{l+1}}\,,$$
et m{\^e}me plus pr{\'e}cis{\'e}ment, $\chi=P_K+R_K$, avec :
$$(\ast\ast)
\begin{cases}
\frac{C'_{\infty}(l+1,1)}{((N+1)\log K)^{l+1}}\leq|h_1(x)|
\leq\|h_1\|_{\infty}
\leq\frac{2^{l+1}C'_{\infty}(l+1,1)}{((N+1)\log K)^{l+1}},\\
h_1(x)P_K(x,\zeta,\zeta_{N+1})\leq0,\\
\left\|R_K(\cdot,\zeta,\zeta_{N+1})
\right\|_{\infty}\leq\frac{C'_{\infty}(l+1,1)}{4((N+1)\log K)^{l+1}}\,.
\end{cases}
$$
Si $\mu\in[0,1]$, on a alors
\begin{eqnarray*}
|h_1(x)-\mu\widetilde{\chi}(x,\zeta,\zeta_{N+1})| & \geq & 
|h_1(x)-\mu P_K(x,\zeta)|-\mu\|R_K(\cdot,\zeta,\zeta_{N+1})\|_{\infty}\\
& \geq & \frac{3C'_{\infty}(l+1,1)}{4((N+1)\log K)^{l+1}}\\
& \geq & \frac{C'}{(N\log N)^{l+1}}\,.
\end{eqnarray*}

Prenons alors 
$\zeta_{N+1}=\ln\frac{e^{b(N)}}{b(N)\psi(0,\zeta)}$ (sachant qu'on
peut ici aussi supposer $\psi(0,\zeta)>0$). De
deux choses l'une : soit
$e^{b(N)}\frac{\psi(x,\zeta)}{\psi(0,\zeta)}\geq1$,
donc en posant $\mu=\frac{\psi(0,\zeta)}{e^{b(N)}\psi(x,\zeta)}$, on a
$$\left|h_1(x)-\frac{1}{b(N)}\frac{\frac{\partial\psi}{\partial x}(x,\zeta)}
{\psi(x,\zeta)}\right|\geq\frac{C'}{(N\log N)^{l+1}}\,;$$
soit
$\psi(x,\zeta)<\frac{\psi(0,\zeta)}{e^{b(N)}}$, donc, quitte {\`a}
rapprocher $x$ de $0$, on peut aussi supposer $\psi(x,\zeta)>0$, et
$\exists\;x'=x'(\zeta)\in[0,1]$, tel que
$$\left|\frac{\frac{\partial\psi}{\partial x}
(x',\zeta)}{\psi(x',\zeta)}\right|>b(N)$$
soit, puisque $\|h_1\|_{\infty}\leq\frac{1}{2}$ (quitte {\`a} le r{\'e}duire),
$$\left|h_1(x')-\frac{\frac{\partial\psi}{\partial x}(x',\zeta)}
{\psi(x',\zeta)}\right|\geq\frac{1}{2}\,.$$

Il ne reste plus qu'{\`a} prendre $h=h_1'\in\Lambda_l([0,1])$, on a
$h_1(x)=\int_0^xh(t)dt$ (puisque $h_1(0)=0$), et $h$ s'annule aussi
identiquement sur $\left[0,\frac{1}{20}\right]$.

\end{proof}

Dans la partie~\ref{pbinverse}, on aura plut{\^o}t recours au

\begin{corollary}\label{detail}

Les corollaires ~\ref{inverse} et ~\ref{inverse2} sont encore valables
(pour $l>1$) si on se limite aux fonctions 
$h>0,\,h\in\Lambda_l\cap C^{[l]}([0,1])$, strictement d{\'e}croissantes et
qui v{\'e}rifient 
$$h^{(j)}(0)=0,\,\forall\,j=1,\ldots,[l].$$

\end{corollary}

\begin{proof}

On se donne $\psi$ et on pose
$$\widetilde{\psi}(x,\zeta,w)=\exp\left(b(N)\left(ax^{[l]+3}-bx^2\right)
\right)\,\psi(x,\zeta,w),$$
$a,\,b$ r{\'e}els $>0$. Alors $\widetilde{\psi}$ v{\'e}rifie encore les
hypoth{\`e}ses du corollaire~\ref{inverse} ($b(N)$ {\'e}tant polynomial en
$N$) : elle est encore de classe $C^2$ sur $[0,1]$, 
$\frac{\partial\widetilde{\psi}}{\partial x}(0,\zeta,w)=
\frac{\partial\psi}{\partial x}(0,\zeta,w)=0,\,\forall\,\zeta,\,w$, et
$$\frac{\partial^2}{\partial x^2}\left(\ln\frac{\widetilde{\psi}(x,\zeta,w)}
{\widetilde{\psi}(0,\zeta,w)}\right)=\frac{\partial^2}{\partial x^2}
\left(\ln\frac{\psi(x,\zeta,w)}{\psi(0,\zeta,w)}\right)
+b(N)(a([l]+3)([l]+2)x^{[l]+1}-2b)\,.$$
D'apr{\`e}s le corollaire~\ref{inverse}, il existe 
$h\in\Lambda_l\cap C^{[l]}([0,1])$ qui v{\'e}rifie, $\forall\,\zeta$,
\begin{eqnarray*}
\left\|h-\frac{1}{b(N)}\left(\frac{\partial^2}{\partial x^2}\ln
\frac{\widetilde{\psi}}{\widetilde{\psi}(0,\zeta,w)}\right)(\cdot,\zeta,w)
\right\|_{\infty} & = & \left\|\widetilde{h}
-\frac{1}{b(N)}
\left(\frac{\partial^2}{\partial x^2}\ln\frac{\psi}{\psi(0,\zeta,w)}\right)
(\cdot,\zeta,w)\right\|_{\infty}\\
& \geq & \frac{\widetilde{C}}{(N\log N)^l}\,,
\end{eqnarray*}
o{\`u} $\widetilde{h}(x)=h(x)+2b-a([l]+3)([l]+2)x^{[l]+1}$.
\bigskip

De m{\^e}me, dans le cas du corollaire~\ref{inverse2}, on a
\begin{eqnarray*}
\left|\int_0^xh(t)dt-\frac{1}{b(N)\widetilde{\psi}(x,\zeta)}
\frac{\partial\widetilde{\psi}}{\partial x}(x,\zeta)\right| & =
& \left|\int_0^x\widetilde{h}(t)dt
-\frac{1}{b(N)\psi(x,\zeta)}
\frac{\partial\psi}{\partial x}(x,\zeta)\right|\\
& \geq & \frac{C'}{(N\log N)^{l+1}}\,.
\end{eqnarray*}

Il ne reste plus qu'{\`a} choisir convenablement $a$ et $b$ : on a
d'abord, sur $\left[\frac{1}{20},1\right]$,
\begin{eqnarray*}
\widetilde{h}'(x) & = & h'(x)-a([l]+3)([l]+2)([l]+1)x^{[l]}\\
& \leq & -a([l]+3)([l]+2)([l]+1)\frac{1}{20^{[l]}}+\|h'\|_{\infty}\\
& < & 0\,,
\end{eqnarray*}
pour $a$ assez grand ($h$ {\'e}tant dans $\Lambda_l([0,1])$), ce qui
montre que $\widetilde{h}$ est srtictement d{\'e}croissante ; ainsi que sur
$\left[0,\frac{1}{20}\right]$ o{\`u} $h$ est constante. $a$ {\'e}tant fix{\'e}, il
ne reste plus qu'{\`a} choisir $b$ assez grand pour avoir
$\widetilde{h}>0$, ce qui montre que $\widetilde{h}$ est strictement
d{\'e}croissante et a ses d{\'e}riv{\'e}es qui s'annulent en $0$ jusqu'{\`a} l'ordre
$[l]$.

Le corollaire est donc prouv{\'e}, {\`a} condition de prendre
$\widetilde{h}$ dans un homoth{\'e}tique de $\Lambda_l([0,1])$. Pour le
cas g{\'e}n{\'e}ral, il suffit de prendre
$\frac{\widetilde{h}}{a([l]+3)([l]+2)+2b+1}$, quitte {\`a} changer $b(N)$
et diminuer $\widetilde{C}$.

\end{proof}

\section{Th{\'e}or{\`e}me de Vitushkin par la m{\'e}thode de Warren}\label{vitwar}
\bigskip

\subsection{Enonc{\'e} et quelques rappels}

On se propose dans cette partie de retrouver le r{\'e}sultat suivant :

\begin{theorem}\label{vitushkin}
Soit $\Lambda_{l,s}\subset C(I^s)$ (resp. $L^1(I^s)$), $l>0$,
et soit
$$P_{n,d}=\left\{\sum_{|k|\leq d}c_k(x)\zeta^k,\,\zeta=(\zeta_1,\ldots,\zeta_n)\in\mathbb{R}^n\right\},$$
o{\`u} $n\geq 1,d\geq 2$ sont des entiers, et $c_k\in C(I^s)$ (resp. $L^1(I^s)$).

On a alors : $\exists\,h\in\Lambda_{l,s},\,\forall\,\zeta\in\mathbb{R}^n$,
$$\|h-P(\zeta)\|\geq\frac{C(l,s)}{(n\log d)^{\frac{l}{s}}},$$ o{\`u}
  $\|\cdot\|$ d{\'e}signe la norme uniforme $\|\cdot\|_{\infty}$
  (resp. $\|\cdot\|_{L^1}$), $C(l,s)=C_{\infty}(l,s)$
  (resp. $C_{L^1}(l,s)$).

De fa{\c c}on {\'e}quivalente, si $\mathcal{P}_{n,d}$ d{\'e}signe l'ensemble des
parties de $C(I^s)$ (resp. $L^1(I^s)$) param{\'e}tr{\'e}es avec $n$ variables,
polynomialement de degr{\'e} (au plus) $d$, on a
$$D_{n,d}(\Lambda_{l,s}):=\inf_{P\in\mathcal{P}_{n,d}}\,\sup_{h\in
\Lambda_{l,s}}\,\inf_{\zeta\in\mathbb{R}^n}\|h-P(\zeta)\|\geq
\frac{C(l,s)}{(n\log d)^{\frac{l}{s}}}.$$

On a en outre le calcul des constantes :
$$C_{\infty}(l,s)=\frac{1}{\sqrt{s}\,2^{l+1}8^{\frac{l}{s}}([l]+1)^{[l]+1}
(4(1+e))^{s([l]+1)}}\;,$$
et
$$C_{L^1}(l,s)=\frac{(([l]+1)!)^{2s}}{5\sqrt{s}\,2^{l+2}18^{\frac{l}{s}}
([l]+1)^{[l]+1}((2[l]+3)!)^s(1+e)^{s([l]+1)}}\;.$$

\end{theorem}

La preuve de ce th{\'e}or{\`e}me utilise la m{\^e}me m{\'e}thode que Warren dans
~\cite{warren}. Elle se fonde seulement sur l'estimation du nombre de
composantes connexes d'un ensemble alg{\'e}brique de $\mathbb{R}^n$. Plus
pr{\'e}cis{\'e}ment, on a :
\begin{proposition}
Soient $p_1,\ldots,p_q$ des polyn{\^o}mes r{\'e}els {\`a} $n$ variables de degr{\'e} au
plus $d$. Alors le nombre de composantes connexes de l'ensemble
$$\mathbb{R}^n\setminus\bigcup_{j=1}^{q}\{\zeta\in\mathbb{R}^n,p_j(\zeta)
=0\}$$ est fini, et est major{\'e} par $\left(\frac{4edq}{n}\right)^n$.
\end{proposition}

L'id{\'e}e est d'estimer le nombre de suites de $\pm1$ atteintes par la
fonction $$\zeta\in\mathbb{R}^n\mapsto(sgn\,p_1(\zeta),\ldots,sgn\,
p_q(\zeta))$$
(en convenant que $sgn\,0=0$), en utilisant le principe que la
restriction de chaque $p_j$ sur chaque composante connexe de 
$\mathbb{R}^n\setminus\bigcup_{j=1}^{q}\{\zeta\in\mathbb{R}^n,p_j(\zeta)
=0\}$ donne une 
fonction continue qui ne s'annule pas, donc de signe constant
(et c'est un argument qu'on ne peut pas transposer dans le cas
complexe). La d{\'e}monstration du th{\'e}or{\`e}me~\ref{vitushkin} repose alors
sur le : 
\begin{corollary}\label{warren}
Si $q\geq 8n\log d$, alors il existe une suite
$\varepsilon=(\varepsilon_1,\ldots,\varepsilon_q)$, o{\`u}
$\varepsilon_j=\pm1$, qui n'est jamais atteinte par 
$sgn\,P(\zeta)=(sgn\,p_1(\zeta),\ldots,sgn\,p_q(\zeta)),\,\zeta\in
\mathbb{R}^n$.

De m{\^e}me, si $q\geq 18n\log d$, il existe une suite $\varepsilon$
qui diff{\`e}re de plus de $\frac{q}{10}$ places de toute suite
$sgn\,P(\zeta),\,\zeta\in\mathbb{R}^n$.
\end{corollary}

Ces deux {\'e}nonc{\'e}s sont prouv{\'e}es par H. E. Warren dans
~\cite{warren}. Ils utilisent les r{\'e}sultats de Oleinik et Petrovskii
sur l'estimation du nombre de composantes connexes d'un ensemble
alg{\'e}brique, qui proviennent finalement du th{\'e}or{\`e}me de B{\'e}zout
(cf.~\cite{oleinik}). 

\subsection{Construction d'une famille de $\Lambda_{l,s}$}

Avant de donner la preuve du th{\'e}or{\`e}me, on doit montrer la proposition
suivante (qui est une reformulation du calcul de la capacit{\'e} de
$\Lambda_{l,s}$).

\begin{proposition}\label{propf}

Soient $r$ un entier $\geq1$, $q=r^s$, et $\varepsilon=(\varepsilon_1,\ldots,\varepsilon_q)$ une suite de
$\pm1$. On d{\'e}finit sur $I^s$ la fonction $f_{\varepsilon}$ de la
fa{\c c}on suivante : on consid{\`e}re le pavage (moyennant les bords) de $I^s$
en les $q=r^s$ cubes $K_i$, $i=1,\ldots,q$. Si $x\in
K_i=\prod_{j=1}^s\left[t_{i,j},t_{i,j}+\frac{1}{r}\right]$, alors $x=t_i+y$, o{\`u}
$t_i=(t_{i,1},\ldots,t_{i,s})$, $y\in\left[0,\frac{1}{r}\right]^s$, et on pose
$$f_{\varepsilon}(x)=\varepsilon_i\frac{g_{l,s}(ry)}{2r^lM_{l,s}},$$
o{\`u} $g_{l,s}$ est d{\'e}finie sur $I^s$ par
$$g_{l,s}(x)=\prod_{j=1}^s(x_j(1-x_j))^{[l]+1},$$
$[l]$ {\'e}tant la partie enti{\`e}re de $l$, et $$M_{l,s}=\sqrt{s}(1+e)^{s([l]+1)}([l]+1)^{[l]+1}.$$

Alors $f_{\varepsilon}$ est bien d{\'e}finie et est un {\'e}l{\'e}ment de
$\Lambda_{l,s}$.

\end{proposition}

On commence par prouver le 

\begin{lemma}

Soit la fonction d{\'e}finie sur $[0,1]$ par $$t\mapsto
t^{[l]+1}(1-t)^{[l]+1}.$$
Alors $\forall k,$ $0\leq k\leq[l]+1$,
$$\left\|\frac{d^k}{dt^k}\left(t^{[l]+1}(1-t)^{[l]+1}\right)\right\|_{\infty}\leq(1+e)^{[l]+1}([l]+1)^k.$$

\end{lemma}

\begin{proof}

$\forall\,k,\,0\leq k \leq[l]+1$, $\forall\,t\in[0,1]$,
$$t^{[l]+1}(1-t)^{[l]+1}=\sum_{j=0}^{[l]+1}C_{[l]+1}^j(-1)^jt^{[l]+j+1},$$
donc
$$\left|\frac{d^k}{dt^k}\left(t^{[l]+1}(1-t)^{[l]+1}\right)\right|\leq\sum_{j=0}^{[l]+1}C_{[l]+1}^j([l]+j+1)\ldots([l]+j-k+2)\,t^{[l]+j+1-k}.$$
Or,
\begin{eqnarray*}
([l]+j+1)\ldots([l]+j-k+2) & \leq &
([l]+1)^k\left(1+\frac{j}{[l]+1}\right)^k\\
& \leq & ([l]+1)^k\exp\frac{jk}{[l]+1},
\end{eqnarray*}
ce qui donne
\begin{eqnarray*}
\left|\frac{d^k}{dt^k}\left(t^{[l]+1}(1-t)^{[l]+1}\right)\right| &
\leq &
([l]+1)^k\sum_{j=0}^{[l]+1}C_{[l]+1}^j\left(e^{\frac{k}{[l]+1}}\right)^j\\
& = & ([l]+1)^k\left(1+e^{\frac{k}{[l]+1}}\right)^{[l]+1},
\end{eqnarray*}
et prouve l'assertion.

\end{proof}

On en d{\'e}duit le

\begin{lemma}\label{lemmeg}

Soit $$g_{l,s}(x)=\prod_{j=1}^s(x_j(1-x_j))^{[l]+1},$$
alors $\frac{g_{l,s}}{M_{l,s}}\in\Lambda_l(I^s)$.

\end{lemma}

\begin{proof}

On voit d'abord que $g_{l,s}\in C^{[l]+1}(I^s)\subset C^m(I^s)$, en tant
que polyn{\^o}me (o{\`u} $l=m+\alpha$).

Il s'agit de montrer que, $\forall\,k=(k_1,\ldots,k_s)$,
$0\leq|k|=k_1+\cdots+k_s\leq m$, 
$$\left\|\frac{\partial^{|k|}g_{l,s}}{\partial
    x^k}\right\|_{\infty}=\left\|\frac{\partial^{|k|}g_{l,s}}{\partial
    x_1^{k_1}\ldots\partial x_s^{k_s}}\right\|_{\infty}\leq M_{l,s},$$
et $\forall\,k$, $|k|=m$, $$\left\|\frac{\partial^mg_{l,s}}{\partial
  x^k}\right\|_{\alpha}\leq M_{l,s}.$$

Pour la premi{\`e}re estimation, prenons m{\^e}me $k$, $|k|=k_1+\cdots+k_s\leq[l]+1$, et soit
$x=(x_1,\ldots,x_s)\in I^s$ :
\begin{eqnarray*}
\frac{\partial^{|k|}g_{l,s}}{\partial x^k}(x) & = & \prod_{j=1}^s\frac{\partial^{k_j}}{\partial x_j^{k_j}}\left(x_j^{[l]+1}(1-x_j)^{[l]+1}\right)\\
& = &
\prod_{j=1}^s\left(\frac{d^{k_j}}{dt^{k_j}}t^{[l]+1}(1-t)^{[l]+1}\right)(x_j),
\end{eqnarray*}
ce qui donne, par le lemme pr{\'e}c{\'e}dent puisque $0\leq k_j\leq[l]+1$,
\begin{eqnarray*}
\left\|\frac{\partial^{|k|}g_{l,s}}{\partial
  x^{k}}\right\|_{\infty} & \leq & \prod_{j=1}^s(1+e)^{[l]+1}([l]+1)^{k_j}\\
& = & (1+e)^{s([l]+1)}([l]+1)^{|k|}\\
& \leq & M_{l,s}.
\end{eqnarray*}
Si $l$ n'est pas entier, cela prouve l'estimation. Si $l$ est entier, on a montr{\'e} un peu plus puisque $[l]=m+1$
: $g_{l,s}\in C^{m+1}(I^s)$ et
$\left\|\frac{1}{M_{l,s}}\frac{\partial^{m+1}g_{l,s}}{\partial x^k}\right\|_{\infty}\leq1$.
\bigskip

Pour la seconde majoration, on a par les accroissements
finis et Cauchy-Schwarz : $\forall\,k,\,|k|=m,\,\forall\,x,\,y\in I^s,$
$$\left|\frac{\partial^mg_{l,s}}{\partial
  x^k}(x)-\frac{\partial^mg_{l,s}}{\partial
  x^k}(y)\right|\leq\left\|\overrightarrow{\nabla}\frac{\partial^mg_{l,s}}{\partial
  x^k}\right\|\|x-y\|,$$
soit puisque $\|x-y\|\leq\|x-y\|^{\alpha}$,
$$\left\|\frac{\partial^mg_{l,s}}{\partial
  x^k}\right\|_{\alpha}\leq\left(\sum_{j=1}^s\left\|\frac{\partial}{\partial x_j}\left(\frac{\partial^mg_{l,s}}{\partial
  x^k}\right)\right\|^2_{\infty}\right)^{\frac{1}{2}}.$$
Puisque $m+1\leq[l]+1$, on obtient d'apr{\`e}s la premi{\`e}re estimation,
\begin{eqnarray*}
\left\|\frac{\partial^mg_{l,s}}{\partial
  x^k}\right\|_{\alpha} & \leq & \left(\sum_{j=1}^s\left((1+e)^{s([l]+1)}([l]+1)^{[l]+1}\right)^2\right)^{\frac{1}{2}}\\
& = & \sqrt{s}(1+e)^{s([l]+1)}([l]+1)^{[l]+1},
\end{eqnarray*}
ce qui prouve le lemme.

\end{proof}

On peut donc prouver la proposition~\ref{propf}.

\begin{proof}

$f_{\varepsilon}$ est effectivement bien d{\'e}finie sur $I^s$ car elle
s'annule sur les bords des $K_i$. D'autre part, $g_{l,s}\in C^{[l]}(I^s)$
et toutes ses d{\'e}riv{\'e}es partielles jusqu'{\`a} $[l]$ (qui sont continues sur
$I^s$) s'annulent sur le bord de $I^s$, ce qui montre que
$f_{\varepsilon}$ construite par recollements, est {\'e}galement dans
$C^{[l]}(I^s)\subset C^m(I^s)$.

On v{\'e}rifie ensuite que la fonction d{\'e}finie sur
$\left[0,\frac{1}{r}\right]^s$ par
$$x\mapsto\frac{1}{r^l}\frac{g_{l,s}(rx)}{M_{l,s}},$$ est bien dans
$\Lambda_{l,s}\left(\left[0,\frac{1}{r}\right]^s\right)$. En effet,
d'apr{\`e}s le lemme~\ref{lemmeg}, $\forall |k|\leq m,\forall
x\in\left[0,\frac{1}{r}\right]^s$,
$$\frac{\partial^{k_1+\cdots+k_s}}{\partial^{k_1}x_1\ldots\partial^{k_s}x_s}(g_{l,s}(rx_1,\ldots,rx_s))=r^{k_1+\cdots+k_s}\left(\frac{\partial^{|k|}g_{l,s}}{\partial
    x^k}\right)(rx),$$
et donc $$\frac{1}{r^l}\left|\frac{\partial^{|k|}}{\partial x^k}(g_{l,s}(rx))\right|\leq\frac{1}{r^{l-|k|}}M_{l,s}\leq
M_{l,s}.$$

De m{\^e}me, si $x\neq y$, $|k|=m$, alors $rx,ry\in I^s$, et :
\begin{eqnarray*}
\frac{1}{\|x-y\|^{\alpha}}\left|\frac{\partial^m}{\partial
    x^k}(g_{l,s}(rx))-\frac{\partial^m}{\partial
    x^k}(g_{l,s}(ry))\right| & = & \frac{r^{\alpha}}{\|rx-ry\|^{\alpha}}r^m\left|\left(\frac{\partial^mg_{l,s}}{\partial
    x^k}\right)(rx)-\left(\frac{\partial^mg_{l,s}}{\partial
    x^k}\right)(ry)\right|\\
& \leq & r^lM_{l,s}.
\end{eqnarray*}

Il en r{\'e}sulte que $f_\varepsilon$ a toutes ses d{\'e}riv{\'e}es d'ordre $\leq m$ born{\'e}es par $1$
en norme uniforme sur $I^s$, car $\forall\,K_i,\;\forall\,|k|\leq m$,
$$\left\|\frac{\partial^{|k|}f_{\varepsilon}}{\partial x^k}|_{K_i}\right\|_{\infty}=\frac{1}{2M_{l,s}}\left\|\frac{\partial^{|k|}g_{l,s}}{\partial^kx}\right\|_{\infty}\leq\frac{1}{2}\leq1.$$

Reste {\`a} s'assurer que, $\forall\,|k|=m$,
$$\left\|\frac{\partial^mf_{\varepsilon}}{\partial x^k}\right\|_{\alpha}\leq1.$$ 
C'est imm{\'e}diat si $x$ et $y$ sont dans le m{\^e}me cube $K_i$, et c'est
m{\^e}me major{\'e} par $\frac{1}{2}$. Sinon, on a
$x\in K_{i_x}$, $y\in K_{i_y}$, avec $i_x\neq i_y$. Le segment $[x,y]$
va donc respectivement couper les bords de $K_{i_x}$ et $K_{i_y}$ en
$z_x$ et $z_y$ (si $x$ est sur le bord de $K_{i_x}$, on prend $z_x=x$
; sinon, l'{\'e}l{\'e}ment $z_x$ est bien d{\'e}fini ; de m{\^e}me pour $y$), ce qui donne, puisque
$\frac{\partial^mf_{\varepsilon}}{\partial
  x^k}(z_x)=\frac{\partial^mf_{\varepsilon}}{\partial x^k}(z_y)=0$, et $\|x-z_x\|,\,\|y-z_y\|\leq\|x-y\|$,
$$\frac{\left|\frac{\partial^mf_{\varepsilon}}{\partial^kx}(x)\right|}{\|x-y\|^{\alpha}}\leq\frac{\left|\frac{\partial^mf_{\varepsilon}}{\partial^kx}(x)-\frac{\partial^mf_{\varepsilon}}{\partial^kx}(z_x)\right|}{\|x-z_x\|^{\alpha}}\leq\frac{1}{2},$$
dans le cas o{\`u} $x\neq z_x$. Sinon, $x$ est sur le bord de $K_{i_x}$,
donc $\frac{\partial^mf_{\varepsilon}}{\partial x^k}(x)=0$, et la
  majoration est triviale.

De m{\^e}me,
$$\frac{\left|\frac{\partial^mf_{\varepsilon}}{\partial^kx}(x)\right|}{\|x-y\|^{\alpha}}\leq\frac{1}{2},$$
donc $$\frac{\left|\frac{\partial^mf_{\varepsilon}}{\partial^kx}(x)-\frac{\partial^mf_{\varepsilon}}{\partial^kx}(y)\right|}{\|x-y\|^{\alpha}}\leq\frac{\left|\frac{\partial^mf_{\varepsilon}}{\partial^kx}(x)\right|}{\|x-y\|^{\alpha}}+\frac{\left|\frac{\partial^mf_{\varepsilon}}{\partial^kx}(y)\right|}{\|x-y\|^{\alpha}}\leq1,$$
ce qui donne la majoration cherch{\'e}e pour tout $|k|\leq m$, et prouve l'assertion.

\end{proof}

\begin{remark}\label{classe}

Comme on l'a vu auparavant, pour $l$ entier $\geq1$, $g_{l,s}$ et
$f_{\varepsilon}$ sont toujours dans l'espace $C^l(I^s)$ usuel :
$g_{l,s}$ comme restriction d'un polyn{\^o}me, et $f_{\varepsilon}$
puisque ses d{\'e}riv{\'e}es partielles sont continues sur chacun des $K_i$ et
qu'elles s'annulent sur les bords des $K_i$. En particulier,
$f_{\varepsilon}\in\Lambda_{l,s}$, car par d{\'e}finition,
$m=l-1,\,\alpha=1$, et $f_{\varepsilon}$ admet des d{\'e}riv{\'e}es
partielles continues d'ordre $l-1$, qui sont lipschitziennes. Ainsi,
le th{\'e}or{\`e}me~\ref{vitushkin} pourra {\'e}galement se formuler avec le
compact $\Lambda_{l,s}\cap C^l(I^s)$.

\end{remark} 

\subsection{Preuve du th{\'e}or{\`e}me de Vitushkin et quelques cons{\'e}quences}

On est maintenant en mesure de donner la preuve du th{\'e}or{\`e}me~\ref{vitushkin}.

\begin{proof}

Traitons dans un premier temps le cas continu, qui utilise la premi{\`e}re assertion du
corollaire~\ref{warren}. $q$ {\'e}tant donn{\'e}, il s'agit de montrer qu'il existe
$\eta(q)$ tel que, pour toute
suite $(\varepsilon_1,\ldots,\varepsilon_q)$, il existe
$f\in\Lambda_{l,s}$ et $q$ formes lin{\'e}aires $\lambda_j$ de norme $\leq
1$ sur $C(I^s)$
qui v{\'e}rifient $$\varepsilon_j\lambda_j(f)\geq\eta(q),\,j=1,\ldots,q$$ (cette
condition traduit la mesure d'oscillation arbitraire de $\Lambda_{l,s}$ autour
de $0$ ; c'est une id{\'e}e qui rejoint la d{\'e}finition
d' $\eta$-capacit{\'e}, cf~\cite{kolmogorov}). Cela entrainera
$$D_{n,d}(\Lambda_{l,s})\geq\eta(q).$$

En effet, {\'e}tant donn{\'e} $P\in\mathcal{P}_{n,d}$, en prenant $q\geq 8n\log d$ et $\varepsilon$ comme dans le corollaire~\ref{warren},
et $\eta(q)$, $f$, $\lambda_1,\ldots,\lambda_q$ associ{\'e}s, on
aura, puisque $\|\lambda_j\|\leq 1$,
$$\sup_{h\in\Lambda_{l,s}}\inf_{\zeta\in\mathbb{R}^n}\|h-P(\zeta)\|_{\infty}\geq\inf_{\zeta\in\mathbb{R}^n}\|f-P(\zeta)\|_{\infty}\geq\inf_{\zeta\in\mathbb{R}^n}\sup_{1\leq
  j\leq q}|\lambda_j(f)-\lambda_j(P(\zeta))|.$$

Pour tout $j$, $p_j(\zeta):=\lambda_j(P(\zeta))$ est un polyn{\^o}me r{\'e}el {\`a} $n$ variables
de degr{\'e} au plus $d$. D'apr{\`e}s le corollaire~\ref{warren}, pour tout
$\zeta\in\mathbb{R}^n$, la suite $sgn\,P(\zeta)=(sgn\,p_1(\zeta),\ldots,sgn\,p_q(\zeta))$
diff{\`e}re de $\varepsilon$ d'au moins une place : $\exists\,k,\,\varepsilon_k\neq sgn
\,p_k(\zeta)$ ($sgn\,p_k(\zeta)$ peut {\^e}tre nul). Dans ce cas,
$$\sup_{1\leq j\leq q}|\lambda_j(f)-\lambda_j(P(\zeta))|\geq
|\lambda_k(f)-p_k(\zeta)|\geq|\lambda_k(f)|\geq\eta(q),$$ ce
qui d{\'e}montre l'in{\'e}galit{\'e} cherch{\'e}e pour $P$ fix{\'e}. L'arbitraire sur $P\in\mathcal{P}_{n,d}$
nous donne finalement $$D_{n,d}(\Lambda_{l,s})\geq\eta(q).$$
\bigskip

Choisissons alors un entier
$r\geq1$, tel que $q=r^s\geq8n\log d$, donnons-nous la suite $\varepsilon$
gr{\^a}ce au corollaire~\ref{warren}, ainsi que $f_{\varepsilon}$
construite {\`a} partir du pavage de $I^s$. On sait d'apr{\`e}s la
proposition~\ref{propf} que $f_{\varepsilon}\in\Lambda_{l,s}$.

Prenons alors pour les $\lambda_i$ les morphismes
d'{\'e}valuation en les centres des $K_i$. On a
$$\lambda_i(f_{\varepsilon})=\varepsilon_i\|f_{\varepsilon}\|_{\infty}=\varepsilon_i\frac{g_{l,s}\left(\frac{1}{2},\ldots,\frac{1}{2}\right)}{2r^lM_{l,s}}=\varepsilon_i\frac{1}{2M_{l,s}4^{s([l]+1)}r^l}.$$ 
En particulier, si $r$ est le plus petit entier ($\geq2$) tel que $r^s\geq
8n\log d$, on a
$$\frac{1}{r}\geq\frac{1}{2(r-1)}\geq\frac{1}{2(8n\log
  d)^{\frac{1}{s}}},$$
soit $$\varepsilon_i\lambda_i(f_{\varepsilon})\geq\frac{1}{2M_{l,s}4^{s([l]+1)}2^l8^{\frac{l}{s}}}\frac{1}{(n\log
  d)^{\frac{l}{s}}},$$ ce qui permet de conclure pour le cas continu :
$$D_{n,d}(\Lambda_{l,s})\geq\frac{C_{\infty}(l,s)}{(n\log
  d)^{\frac{l}{s}}},$$ avec
$$C_{\infty}(l,s)=\frac{1}{2^{l+1}4^{s([l]+1)}8^{\frac{l}{s}}M_{l,s}}\,.$$
\bigskip

Traitons maintenant le cas int{\'e}grable (qui utilise la deuxi{\`e}me
assertion du corollaire~\ref{warren}). 
Comme pour le cas continu, choisissons $q=r^s$, $r\geq2$, tel que
$r^s$ soit le plus petit entier $\geq18n\log d$.

Soient $\varepsilon$ une suite, $f_{\varepsilon}\in\Lambda_{l,s}$ la
fonction associ{\'e}e, et pour les $\lambda_i$,
posons $$\lambda_i(h)=\int_{K_i}h(x)dx.$$
On a alors
$$\sum_{i=1}^m|\lambda_i(h)|\leq\sum_{i=1}^m\int_{K_i}|h(x)|dx=\|h\|_{L^1},$$
donc $$\sum_{i=1}^m\|\lambda_i\|\leq1.$$
D'autre part,
\begin{eqnarray*}
\lambda_i(f_{\varepsilon}) & = & \varepsilon_i\frac{1}{2r^lM_{l,s}}\int_{\left[0,\frac{1}{r}\right]^s}g_{l,s}(ry)dy\\
& = & \varepsilon_i\frac{1}{2r^lM_{l,s}r^s}\int_{[0,1]^s}\prod_{j=1}^s(x_j(1-x_j))^{[l]+1}dx_1\ldots
dx_s\\
& = & \varepsilon_i\frac{1}{2r^lM_{l,s}r^s}\left(\int_0^1(t(1-t))^{[l]+1}dt\right)^s.
\end{eqnarray*}

On peut calculer cette int{\'e}grale par
int{\'e}grations par parties successives, ce qui nous donne
$$\int_0^1(t(1-t))^{[l]+1}dt=\frac{(([l]+1)!)^2}{(2[l]+3)!},$$
et donc
$$\varepsilon_i\lambda_i(f_{\varepsilon})\geq\frac{(([l]+1)!)^{2s}}{2M_{l,s}((2[l]+3)!)^sr^lr^s},$$
ce qui montre une propri{\'e}t{\'e} analogue d'oscillation arbitraire de
$\Lambda_{l,s}$ autour de $0$.
\bigskip

Soit alors $P\in\mathcal{P}_{n,d}$ ; les $p_i(\zeta)=\lambda_i(P(\zeta))$ sont
des polyn{\^o}mes r{\'e}els {\`a} $n$ variables de degr{\'e} $\leq d$. D'apr{\`e}s le
corollaire et par hypoth{\`e}se sur $q=r^s$, il existe une suite
$\varepsilon$ qui diff{\`e}re d'au moins $\frac{q}{10}$ places de
$sgn\,P(\zeta)=(sgn\,p_1(\zeta),\ldots,sgn\,p_q(\zeta)),\forall\,\zeta\in\mathbb{R}^n$,
ce qui entra{\^\i}ne, puisque $\sum_{i=1}^q\|\lambda_i\|\leq1$,
$$\sup_{h\in\Lambda_{l,s}}\inf_{\zeta\in\mathbb{R}^n}\|h-P(\zeta)\|_{L^1}\geq\inf_{\zeta\in\mathbb{R}^n}\|f_{\varepsilon}-P(\zeta)\|_{L^1}\geq\inf_{\zeta\in\mathbb{R}^n}\sum_{i=1}^q|\lambda_i(f_{\varepsilon})-\lambda_i(P(\zeta))|.$$
Comme, pour au moins $\frac{q}{10}$ indices $i$, on a
$$|\lambda_i(f_{\varepsilon})-p_i(\zeta)|\geq\varepsilon_i\lambda_i(f_{\varepsilon}),$$
il vient
$$\inf_{\zeta\in\mathbb{R}^n}\|f_{\varepsilon}-P(\zeta)\|\geq\frac{q}{10}\frac{(([l]+1)!)^{2s}}{r^s2M_{l,s}((2[l]+3)!)^sr^l}=\frac{(([l]+1)!)^{2s}}{20M_{l,s}((2[l]+3)!)^sr^l}.$$
Enfin, le fait d'avoir ici encore
$$\frac{1}{r}\geq\frac{1}{2(18n\log d)^{\frac{1}{s}}},$$ et
l'arbitraire sur $P\in\mathcal{P}_{n,d}$ aboutissent {\`a}
$$D_{n,d}(\Lambda_{l,s})\geq\frac{C_{L^1}(l,s)}{(n\log
  d)^{\frac{l}{s}}},$$
avec
  $$C_{L^1}(l,s)=\frac{(([l]+1)!)^{2s}}{20M_{l,s}2^l18^{\frac{l}{s}}((2[l]+3)!)^s},$$ 
ce qui ach{\`e}ve la preuve dans le cas $L^1(I^s)$, et prouve le th{\'e}or{\`e}me.

\end{proof}

On peut en particulier expliciter les constantes de minoration.

\begin{corollary}\label{constantes}

On a : pour le cas continu,
$$D_{n,d}(\Lambda_{l,s})\geq\frac{C_{\infty}(l,s)}
{(n\log d)^{\frac{l}{s}}},$$
o{\`u}
$$C_{\infty}(l,s)=\frac{1}{\sqrt{s}\,2^{l+1}8^{\frac{l}{s}}([l]+1)^{[l]+1}
(4(1+e))^{s([l]+1)}}\;;$$
et pour le cas $L^1$,
$$D_{n,d}(\Lambda_{l,s})\geq\frac{C_{L^1}(l,s)}
{(n\log d)^{\frac{l}{s}}},$$
o{\`u}
$$C_{L^1}(l,s)=\frac{(([l]+1)!)^{2s}}{5\sqrt{s}\,2^{l+2}18^{\frac{l}{s}}
([l]+1)^{[l]+1}((2[l]+3)!)^s(1+e)^{s([l]+1)}}.$$
\bigskip

En particulier, on trouve pour la constante $C'$ du th{\'e}or{\`e}me~\ref{th1} :
\begin{eqnarray*}
\lefteqn{C'_{\infty}=\frac{3}{\sqrt{s}\,2^{l+3}8^{\frac{l}{s}}
([l]+1)^{[l]+1}(4(1+e))^{s([l]+1)}[d(r+1)+v+t]^{\frac{l}{s}}}}\\
& &
\times\,\frac{1}{\left[\log\left[2ebud2^{(4e+1)d}B^d\left(1+\frac{l}{s}
\right)^2\log\left(A\sqrt{s}\,2^{l+3}8^{\frac{l}{s}}([l]+1)^{[l]+1}
(4(1+e))^{s([l]+1)}\right)\right]\right]^{\frac{l}{s}}}\;,
\end{eqnarray*}
et
\begin{eqnarray*}
\lefteqn{C'_{L^1}=\frac{3(([l]+1)!)^{2s}}{5\sqrt{s}\,2^{l+4}18^{\frac{l}{s}}
([l]+1)^{[l]+1}((2[l]+3)!)^s(1+e)^{s([l]+1)}\left[(d(r+1)+v+t)\right]
^{\frac{l}{s}}}}\\
& &
\times\,\frac{1}{\left[\log\left[2ebud2^{(4e+1)d}B^d\left(1+\frac{l}{s}
\right)^2\log\left(\frac{5A\sqrt{s}\,2^{l+4}18^{\frac{l}{s}}([l]+1)^{[l]+1}
((2[l]+3)!)^s(1+e)^{s([l]+1)}}{(([l]+1)!)^{2s}}\right)\right]\right]
^{\frac{l}{s}}}\;.
\end{eqnarray*}

\end{corollary}
\bigskip

Pour justifier la preuve du corollaire~\ref{inverse}, on a besoin d'un
r{\'e}sultat donnant quelques restrictions sur la 
fonction $h=f_{\varepsilon}$, quitte {\`a} diminuer la constante
$C_{\infty}(l,1)$.

\begin{corollary}\label{1min}

Dans l'{\'e}nonc{\'e} du th{\'e}or{\`e}me~\ref{vitushkin}, on peut de plus prendre $h$
nulle sur un sous-intervalle $[0,\delta]$ (ici
$\left[0,\frac{1}{20}\right]$).

En particulier, cela justifie ce qu'on a dit dans la
remarque~\ref{remth1}, et permet surtout de prendre $h$ (intervenant 
dans la preuve des th{\'e}or{\`e}mes~\ref{th1} et~\ref{th1bis}) qui
v{\'e}rifie pour le cas continu : 
$$
\begin{cases}
\frac{C'_{\infty}(l,s)}{(n\log K)^{\frac{l}{s}}}\leq\|h\|_{\infty}=
|h(x(\zeta))|\leq\frac{2^lC'_{\infty}(l,s)}{(n\log K)^{\frac{l}{s}}},\\
h(x(\zeta),\zeta)P_K(x(\zeta),\zeta)\leq0\\
h(x)=0,\;\forall\,x\in\left[0,\frac{1}{20}\right].
\end{cases}
$$
(avec une autre constante $C'_{\infty}(l,s)\leq C_{\infty}(l,s)$)

Cela justifie finalement les hypoth{\`e}ses suppl{\'e}mentaires utilis{\'e}es dans
la preuve des corollaires~\ref{inverse} et~\ref{inverse2}.

\end{corollary}

\begin{proof}

Il suffit juste de s'assurer que l'{\'e}l{\'e}ment $x(\zeta)$ peut {\^e}tre choisi
dans l'intervalle $\left[\frac{1}{20},1\right]$, ce qui d{\'e}coule de la
seconde assertion du corollaire~\ref{warren} : si $q\geq18n\log d$, il
existe une suite $\varepsilon$ qui diff{\`e}re de plus de $\frac{q}{10}$
places de toute suite $sgn\,P(\zeta),\;\zeta\in\mathbb{R}^n$. Il
existe donc au moins une place pour $x(\zeta)$ se trouvant en dehors de
$$\left[0,\frac{1}{q}\right]\cup\cdots\cup\left[\frac{k-1}{q},
\frac{k}{q}\right]\supset\left[0,\frac{1}{20}\right]\,.$$
$k$ {\'e}tant le plus grand entier $\leq\frac{q}{10}$. Il ne reste plus
qu'{\`a} changer $h$ en $0$ sur $\left[0,\frac{k}{q}\right]$.

\end{proof}

\section{Application {\`a} un probl{\`e}me inverse}\label{pbinverse}
\bigskip

\subsection{Quelques rappels et motivation du probl{\`e}me}

On se donne l'{\'e}quation de Sturm-Liouville sur la demi-droite
$\mathbb{R}^{+}$
$$-y''(x)-\omega^2Q(x)=\lambda y(x),$$ avec les hypoth{\`e}ses suivantes :
$-\omega^2Q$ est un potentiel strictement n{\'e}gatif, int{\'e}grable avec
$m+1$ d{\'e}riv{\'e}es {\`a} d{\'e}croissance polynomiale en l'infini.

On sait que pour tout $\omega$ assez grand, l'op{\'e}rateur
$-\frac{d^2}{dx^2}-\omega^2Q$ admet $N(\omega)$ valeurs propres
discr{\`e}tes n{\'e}gatives 
$\lambda_j=-\xi_j^2$, $0<\xi_1<\ldots<\xi_N$, et $N(\omega)$ fonctions
propres 
$\varphi_j$ qui v{\'e}rifient les conditions aux bords :
$$\varphi_j(0)=0, \text{ et } \int_{0}^{\infty}|\varphi_j(x)|^2dx=1.$$
En outre, $N(\omega)$ est du m{\^e}me ordre de grandeur que
$\omega$ ; plus pr{\'e}cis{\'e}ment, on dispose de la borne de Calogero
$N(\omega)\leq\omega\frac{2}{\pi}\int_0^{\infty}\sqrt{Q(x)}\,dx$
(cf. ~\cite{lax-levermore}, ~\cite{reed-simon}), et m{\^e}me d'un
encadrement : pour tout $\omega$ assez grand, 
$a\omega\leq N(\omega)\leq b\omega$.
\bigskip

On consid{\`e}re ici un probl{\`e}me inverse : {\'e}tant donn{\'e}es les
$N(\omega)$ valeurs propres 
$$\lambda_j=-\xi_j^2$$ 
et valeurs caract{\'e}ristiques 
$$C_j=\left(\varphi'_j(0)\right)^2,$$
$j=1,\ldots,N(\omega)$, il s'agit de reconstruire le potentiel
$-\omega^2Q$ sur $\mathbb{R}^+$. C'est un cas particulier d'un probl{\`e}me plus
g{\'e}n{\'e}ral qui consiste, 
{\`a} partir de donn{\'e}es sur les fonctions propres, {\`a} retrouver $Q$. On
s'int{\'e}resse essentiellement au cas o{\`u} les informations sont donn{\'e}es
par la fonction de Weyl, d{\'e}finie par
$$j(k)=\frac{\varphi'(0,k)}{\varphi(0,k)},\;\Im m\,k>0,$$
o{\`u} $\varphi$ est une solution $L^2$-int{\'e}grable de l'{\'e}quation avec
$\lambda=k^2$. C'est une fonction m{\'e}romorphe dont les p{\^o}les sont
exactement les valeurs propres $\xi_j$ de l'{\'e}quation, et dont les
r{\'e}sidus sont ({\`a} $2i\xi_j$ pr{\`e}s) les valeurs caract{\'e}ristiques $C_j$. La
connaissance de cette fonction permet de d{\'e}terminer la mesure
spectrale $\sigma(d\tau)$ de l'{\'e}quation d{\'e}finie par
$$\sigma(d\tau)=
\begin{cases}
\sigma_+(d\tau),\;\tau\geq0,\\
\sum_{j=1}^{N(\omega)}C_j\delta(\tau+\xi_j^2),\;\tau<0,
\end{cases}
$$
o{\`u} $\sigma_+$ est une mesure positive {\`a} densit{\'e}, et $\delta$ la mesure
de Dirac. Gr{\^a}ce aux travaux de Gelfand et Levitan (cf ~\cite{gelfand},
~\cite{levitan}), on est alors capable de reconstruire $Q$.

On va ici consid{\'e}rer le cas o{\`u} seuls sont connus les param{\`e}tres
$\xi_j$ et $C_j$, ainsi que $Q$ et ses premi{\`e}res d{\'e}riv{\'e}es en
$0$. On a pour cela les r{\'e}sultats de G. H. Henkin et N. N. Novikova dans
~\cite{henkin} (th{\'e}or{\`e}me 1, p. 21) : sous certaines hypoth{\`e}ses de
r{\'e}gularit{\'e} de 
$Q$ (que l'on suppose avec $m+1$ d{\'e}riv{\'e}es int{\'e}grables), on peut
l'approcher, uniform{\'e}ment sur tout intervalle $[0,X]$, par une fonction
$Q_{\omega}$ {\`a} l'ordre $\frac{1}{\omega^m}$, d{\'e}duite du potentiel 
$$q_{\omega}(x)=-\omega^2Q_{\omega}(x),$$
associ{\'e} {\`a} la mesure spectrale explicite $\sigma_{\omega}(d\tau)$ 
construite {\`a} partir des $\xi_j$, $C_j$ et des $Q^{(s)}(0)$,
$s=0,\ldots,m$.

Pour $m=1$, on dispose d'un potentiel explicite
$q^0_{\omega}=-\omega^2Q_{\omega}^0$, d{\'e}fini par (cf~\cite{henkin}, p. 23)
$$Q_{\omega}^0(x)=\frac{2}{\omega^2}\frac{d^2}{dx^2}\ln|\det W_{s,r}(x)|,$$
o{\`u}
$$W_{s,r}(x)=\frac{2sh(\xi_s+\xi_r)x}{\xi_s+\xi_r}-(1-\delta_{s,r})
\frac{2sh(\xi_s-\xi_r)x}{\xi_s-\xi_r}-\delta_{s,r}\left(2x-\frac{4\xi_r^2}
{C_r}\right),$$
$s,\,r=1,\ldots,N(\omega)$, et qui v{\'e}rifie, uniform{\'e}ment sur tout
$[0,X]$,
$$\left|\int_0^xQ(y)dy-\int_0^xQ_{\omega}^0(x)\right|=
O\left(\frac{\ln\omega}{\sqrt{\omega}}\right).$$
\bigskip

Une question alors assez naturelle concerne le probl{\`e}me de l'optimalit{\'e}
d'une telle formule approximation : G. Henkin a conjectur{\'e}
dans~\cite{henkin} 
p. 22, qu'une telle formule utilisant $2N$ param{\`e}tres ({\`a} savoir les
$\xi_j$, $C_j$) ne pouvait pas approcher uniform{\'e}ment une fonction en
position g{\'e}n{\'e}rale avec $m$ d{\'e}riv{\'e}es born{\'e}es (soit le compact
$\Lambda_m$), essentiellement mieux qu'{\`a} l'ordre
$\frac{1}{\omega^m}$. En d'autres 
termes, peut-on trouver, du moins prouver l'existence d'une formule
qui r{\'e}aliserait une approximation meilleure ?

En supposant qu'une telle formule donn{\'e}e s'{\'e}crit comme une fonction
analytique en les param{\`e}tres $\xi_j$, $C_j$ (et $Q^{(s)}(0)$), et en
faisant le lien avec les r{\'e}sultats n{\'e}gatifs prouv{\'e}s dans la
partie~\ref{appranal}, on va
pouvoir donner des bornes inf{\'e}rieures qui ne peuvent pas {\^e}tre franchies.

\subsection{Une estimation des valeurs propres et valeurs caract{\'e}ristiques}

On a besoin d'{\'e}tablir les encadrements suivants pour pouvoir appliquer
les r{\'e}sultats n{\'e}gatifs.

\begin{proposition}\label{estimation}

Supposons que $q=-\omega^2Q$ soit un potentiel n{\'e}gatif, o{\`u} $Q>0$ est
int{\'e}grable de classe 
$C^1$ avec $Q'(0)=0$, strictement d{\'e}croissant sur $\mathbb{R}^+$ et {\`a}
d{\'e}croissance polynomiale {\`a} l'infini. Alors pour tous
$\omega$ assez grand et $j=1,\ldots,N(\omega)$, on a  
$$\frac{a}{\omega^b}\leq\xi_j\leq c\omega, \text{ et }
\frac{1}{\alpha\exp\left(\beta\omega^{\gamma}\right)}\leq
\frac{4\xi_j^2}{C_j}\leq\alpha \exp\left(\beta\omega^{\gamma}\right),$$
o{\`u} $a,b,c,\alpha,\beta,\gamma$ sont des constantes ne d{\'e}pendant que de $Q$.

\end{proposition}

Cette proposition {\'e}tant un corollaire de la th{\'e}orie WKB qui n'est pas
pr{\'e}cis{\'e}ment formul{\'e} dans les r{\'e}f{\'e}rences, nous allons en donner une
preuve abr{\'e}g{\'e}e en utilisant principalement la m{\'e}thode WKB comme dans
~\cite{lax-levermore}. 
\begin{proof}

D'abord, l'estimation $\xi_j\leq c\omega$ est imm{\'e}diate : car si
$\phi_j$ est une fonction propre norm{\'e}e associ{\'e}e, on a
\begin{eqnarray*}
-\xi_j^2\,=\,\lambda_j & = & -\int_0^{\infty}\phi_j''\phi_jdx-\omega^2\int_0^{\infty}Q\,\phi_j^2dx\\
& \geq &
-\left[\phi_j'\phi_j\right]_0^{\infty}+\int_0^{\infty}\phi_j'^2dx-\omega^2\inf_{[0,+\infty[}Q\\
& \geq & -\omega^2\sup_{[0,+\infty[}Q,
\end{eqnarray*}
soit
$$\xi_j\leq\omega\sup_{[0,+\infty[}\sqrt{Q(x)}\,.$$
\bigskip

Passons maintenant {\`a} la minoration de $\xi_j$ : d'abord, quitte {\`a}
diviser par $Q(0)$, on peut supposer $Q(0)=1$. De plus, $Q$ est par
hypoth{\`e}se strictement d{\'e}croissant sur 
$\mathbb{R}^+$, donc admet $0$ comme unique maximum. Pour
adapter la m{\'e}thode de WKB par Lax et Levermore dans
~\cite{lax-levermore}, on doit prolonger l'equation dans $\mathbb{R}$ tout
entier de la mani{\`e}re suivante : on prolonge $Q$ en $\widetilde{Q}$ par
parit{\'e}, ce qui fait de $\widetilde{q}=-\omega^2\widetilde{Q}$ un potentiel 
n{\'e}gatif, int{\'e}grable {\`a} d{\'e}croissance polynomiale et de classe $C^1$
(puisque $Q'(0)=0$), et qui est monotone sur $\mathbb{R}^+$ 
et $\mathbb{R}^-$, avec $0$ comme unique minimum sur $\mathbb{R}$ (on
continuera de noter $Q$ le prolongement).

Quant {\`a} chaque fonction propre $\phi_j$, on la prolonge par imparit{\'e},
ce qui est possible car $\phi_j(0)=0$, et nous donne une fonction
$\widetilde{\phi}_j$ de classe $C^2$ sur $\mathbb{R}$ (puisque
$\phi_j''(0)=0$).

On v{\'e}rifie imm{\'e}diatement que chaque $\widetilde{\phi}_j$ est fonction
propre de 
l'op{\'e}rateur $-\frac{d^2}{dx^2}-\omega^2Q$, avec la m{\^e}me valeur
propre $\lambda_j=-\xi_j^2$ ; inversement, par restriction sur
$\mathbb{R}^+$, toute fonction propre sur $\mathbb{R}$ est l'une des
$\phi_j$, ce qui montre qu'on les a toutes.

On se place finalement
dans le cadre quasi-classique en ramenant l'{\'e}quation {\`a}
$$-\varepsilon^2y''-Qy=-\eta^2y,$$
o{\`u} $\varepsilon=\frac{1}{\omega}$ et
$\eta_j=\frac{\xi_{N-j+1}}{\omega}$, $j=1,\ldots,N$. On a donc
$$0\leq\eta_N\leq\ldots\leq\eta_1\leq 1.$$
En appliquant la m{\'e}thode WKB, on a :
$$\Phi(\eta_j)=\left(j-\frac{1}{2}\right)\varepsilon\pi,$$
o{\`u}
$$\Phi(\eta)=\int_{x_-(\eta)}^{x_+(\eta)}\left(Q(y)
-\eta^2\right)^{\frac{1}{2}}dy,$$
avec $x_+(\eta)>0$, $x_-(\eta)<0$ v{\'e}rifiant
$-Q(x_+(\eta))=-Q(x_-(\eta))=-\eta^2$ (tout est int{\'e}grable car $Q$ est
suffisamment r{\'e}guli{\`e}re {\`a} l'infini).

On a aussi
$$N(\omega)=N(\varepsilon)=\left[\frac{1}{\varepsilon\pi}\Phi(0)\right]=\left[\frac{1}{\varepsilon\pi}\int_{-\infty}^{\infty}(Q(y))^{\frac{1}{2}}dy\right].$$
Et pour finir,
$$s_j=s(\eta_j)=\exp\left(\frac{\theta_+(\eta_j)}{\varepsilon}\right),$$
o{\`u} 
$$\theta_+(\eta)=\eta x_+(\eta)+\int_{x_+(\eta)}^{\infty}\eta
-\left(\eta^2-Q(y)\right)^{\frac{1}{2}}dy,$$
et $s_j$ est le coefficient normalisateur de la solution de Jost :
$$\phi_j(x)\sim s_j\exp\left(-\frac{\eta_j x}{\varepsilon}\right).$$
\bigskip

Ainsi, pour montrer l'estimation de $\xi_j$, il suffit de s'assurer que
$$\eta_N\geq\frac{a}{\omega^b}$$
(pour plus de clart{\'e}, on d{\'e}signera par $a$ plusieurs constantes
diff{\'e}rentes). Puisque par hypoth{\`e}se, au voisinage de l'infini, 
$$\frac{1}{a}\frac{1}{|x|^k}\leq Q(x)\leq\frac{a}{|x|^k},$$
(avec $k\geq4$) et que, pour $\omega$ grand, $\eta_{N}\rightarrow0$,
on a (par exemple pour $x_+(\eta_N)$, le cas $x_-(\eta_N)$ {\'e}tant analogue)
$$\frac{1}{a}\frac{1}{x_+(\eta_N)^k}\leq
Q(\eta_N)=\eta_N^2\leq\frac{a}{(x_+(\eta_N))^k},$$
soit
$$\frac{1}{a}\frac{1}{(\eta_N)^{\frac{2}{k}}}\leq
x_+(\eta_N)\leq\frac{a}{(\eta_N)^{\frac{2}{k}}}.$$
On a alors
$$\Phi(0)-\Phi(\eta_N)=\left(\int_{-\infty}^{x_-(\eta_N)}+\int_{x_+(\eta_N)}^{\infty}\right)(Q(y))^{\frac{1}{2}}dy\;+\;\int_{x_-(\eta_N)}^{x_+(\eta_N)}(Q(y))^{\frac{1}{2}}-\left(Q(y)-\eta_N^2\right)^{\frac{1}{2}}dy.$$
Or,
$$\int_{x_+(\eta_N)}^{\infty}(Q(y))^{\frac{1}{2}}dy=\int_0^{\frac{1}{x_+(\eta_N)}}\left(\frac{Q\left(\frac{1}{y}\right)}{y^4}\right)^{\frac{1}{2}}dy\leq\frac{a}{x_+(\eta_N)^{\frac{k-2}{2}}},$$
de m{\^e}me pour
$\int_{-\infty}^{x_-(\eta_N)}(Q(y))^{\frac{1}{2}}dy$.
Par ailleurs,
\begin{eqnarray*}
\int_{x_-(\eta_N)}^{x_+(\eta_N)}(Q(y))^{\frac{1}{2}}-\left(Q(y)-\eta_N^2\right)^{\frac{1}{2}}dy
& \leq &
\int_{x_-(\eta_N)}^{x_+(\eta_N)}\left(Q(y)-\left(Q(y)-\eta_N^2\right)\right)^{\frac{1}{2}}dy\\
& = & \eta_N(x_+(\eta_N)-x_-(\eta_N)),
\end{eqnarray*}
ce qui donne $$\Phi(0)-\Phi(\eta_N)\leq
a(\eta_N)^{1-\frac{2}{k}}+a(\eta_N)^{1-\frac{2}{k}}+a(\eta_N)^{1-\frac{2}{k}}=a(\eta_N)^{1-\frac{2}{k}}.$$
On a d'autre part
$N=\left[\frac{1}{\varepsilon\pi}\Phi(0)\right]\leq\frac{\Phi(0)}{\varepsilon\pi}$,
donc $$\Phi(0)-\Phi(\eta_N)\geq\frac{\varepsilon\pi}{2},$$
et finalement
$$a(\eta_N)^{\frac{k-2}{k}}\geq\frac{\varepsilon\pi}{2},$$
soit $$\eta_N\geq
a\varepsilon^{\frac{k}{k-2}}=\frac{a}{\omega^{\frac{k}{k-2}}}\geq\frac{a}{\omega^2},$$
qui est l'estimation cherch{\'e}e.
\bigskip

\bigskip

On passe maintenant {\`a} l'estimation de $\frac{4\xi_j^2}{C_j}$. On
{\'e}tablit pour cela une formule explicite de $C_j$ en reprenant
~\cite{chadan} (chap. II).

On consid{\`e}re l'{\'e}quation
$$-y''-\omega^2Qy=k^2y,$$
de param{\`e}tre $k$ dans le demi-plan complexe sup{\'e}rieur. Pour tout
$k,\;\Im m\,k>0$, il existe toujours les solutions
$\phi(k,\cdot),\,\psi(k,\cdot)$, et $f(k,\cdot)$ qui v{\'e}rifient
respectivement
$$
\begin{cases}
\phi(k,\cdot)\in L^2([0,\infty[),\text{ et
}\int_0^{\infty}|\phi(k,x)|^2dx=1.\\
\psi(k,0)=0\;\text{ et }\;\psi'(k,0)=1.\\
f(k,x)\sim\exp(ikx),\;x\rightarrow\infty.\\
\end{cases}
$$
$\phi$ correspond {\`a} la solution physique et $f$ est appel{\'e}e solution
de Jost. Lorsque $k$ est une valeur propre $k_j=i\xi_j,\;\xi_j>0$ (et
seulement dans ce 
cas), ces trois solutions sont proportionnelles (et peuvent {\^e}tre
prises r{\'e}elles), et on a
$$
\begin{cases}
f(k_j,\cdot)=f'(k_j,0)\psi(k_j,\cdot),\\
\phi(k_j,\cdot)=C_j^{\frac{1}{2}}\psi(k_j,\cdot),\\
\phi(k_j,\cdot)=s_jf(k_j,\cdot)
\end{cases}
$$
(les solutions $\phi(k,\cdot)$ et $f(k,\cdot)$ sont toujours
proportionnelles puisque $\Im mk>0$), ce qui donne
$$f'(k_j,0)=\frac{C_j^{\frac{1}{2}}}{s_j}\,.$$
On d{\'e}finit par ailleurs la fonction de Jost $F$ par $F(k)=f(k,0)$. Ce
qui pr{\'e}c{\`e}de montre qu'elle admet exactement les valeurs propres comme
z{\'e}ros. On verra par la suite que c'est une fonction holomorphe sur le
demi-plan sup{\'e}rieur et que ses z{\'e}ros sont simples. On montre alors que
$$\frac{4\xi_j^2}{C_j}=-s_j^2\left(\dot{F}(i\xi_j)\right)^2.$$
\bigskip

En effet, en d{\'e}rivant par rapport {\`a} la variable $k$ l'{\'e}quation
$$-f''(k,x)-\omega^2Q(x)f(k,x)=k^2f(k,x),$$
la d{\'e}riv{\'e}e (par rapport {\`a} $x$) du wronskien de $f(k,x)$ et
$\dot{f}(k,x)=\frac{df}{dk}(k,x)$ vaut 
\begin{eqnarray*}
W'\left(f,\dot{f}\right) & = &
\left(f(k,x)\dot{f}'(k,x)-f'(k,x)\dot{f}(k,x)\right)'\\
& = & f(k,x)\dot{f}''(k,x)-f''(k,x)\dot{f}(k,x)\\
& = & -2kf^2(k,x),
\end{eqnarray*}
ce qui donne, pour $k=k_j$,
\begin{eqnarray*}
\dot{F}(k_j)f'(k_j,0) & = &
-2k_j\int_0^{\infty}f^2(k_j,x)dx=-\frac{2k_j}{s_j^2},\\
& = & \dot{F}(k_j)\frac{C_j^{\frac{1}{2}}}{s_j},
\end{eqnarray*}
et prouve l'expression annonc{\'e}e.
\bigskip

On a par ailleurs pour la majoration de $s_j$ :
$0<\eta_N<\ldots<\eta_1\leq1$, donc $0\leq x_+(\eta_1)\leq\ldots\leq
x_+(\eta_N)$, et $$\theta_+(\eta_j)\leq
x_+(\eta_N)+\int_0^{\infty}(Q(y))^{\frac{1}{2}}dy,$$
qui, en vertu des estimations pr{\'e}c{\'e}dentes, donne
$$\theta_+(\eta_j)\leq a\omega^{\frac{2}{k-2}},$$
et aboutit {\`a}
$$s_j\leq\exp\left(a\omega^{\frac{k}{k-2}}\right)\leq\exp
\left(a\omega^2\right).$$
Et puisque $s_j\geq1$, il suffit pour terminer la preuve de
montrer un encadrement du m{\^e}me type pour $\dot{F}(i\xi_j)$, ce
qui fait l'objet du lemme suivant.

\end{proof}

\begin{lemma}

On a, pour tout $j=1,\ldots,N(\omega)$,
$$\frac{1}{a\exp(b\omega^2)}\leq\left|\dot{F}(i\xi_j)\right|\leq
a\exp\left(b\omega^2\right).$$
(les constantes $a,\,b$ ne d{\'e}pendent que de $Q$)

\end{lemma}

\begin{remark}

En particulier, on trouve, pour tout $j=1,\ldots,N(\omega)$,
$$\frac{1}{a\omega}\leq\xi_j\leq a\omega,$$
et
$$\frac{1}{A\exp\left(b\omega^2\right)}\leq\frac{4\xi_j^2}{C_j}\leq
A\exp\left(b\omega^2\right),$$
ce qui donne une majoration des exposants $c$ et $\gamma$ dans
l'{\'e}nonc{\'e} de la proposition~\ref{estimation}.

\end{remark}

\begin{proof}

On commence pour cela par consid{\'e}rer l'{\'e}quation $$y''-Vy=-k^2y,$$
o{\`u} $x\geq0$, $V(x)=-\omega^2Q(x)$ et $\Im m\,k\geq0$. On peut
construire la solution de Jost par approximations successives en
posant
$$
\begin{cases}
f_0(k,x)=e^{ikx},\\
f_n(k,x)=\int_x^{\infty}\frac{\sin k(x-t)}{k}V(t)f_{n-1}(t)dt,
\end{cases}
$$
ce qui donne
$$f(k,\cdot)=\sum_{n=0}^{\infty}f_n(k,\cdot),\;\text{et}\;F(k)=f(k,0)$$
pour la fonction de Jost. On sait (cf.~\cite{chadan}) que $F$ est
holomorphe sur le demi-plan $\{\Im m\,k>0\}$, continue sur $\{\Im
m\,k\geq0\}$, s'annule exactement en les valeurs propres $i\xi_j$, et
tend vers $1$ {\`a} l'infini. Elle poss{\`e}de en outre
l'estimation suivante : $\forall\,x,\,k$,
$$|f_n(k,x)|\leq\exp(-\Im
m\,kx)\,\frac{1}{n!}\left(\int_x^{\infty}\frac{2\sqrt{2}\,t}{1+|k|t}|V(t)|dt\right)^n,$$
ce qui donne, pour tout $k$,
\begin{eqnarray*}
|F(k)|\leq\sup_{x\in\mathbb{R}}|f(k,x)| & \leq & 
\exp\left(2\sqrt{2}\,\omega^2\int_0^{\infty}t|Q(t)|dt\right)\\
& = & \exp\left(a\omega^2\right),
\end{eqnarray*}
o{\`u} $a=a(Q)$ ne d{\'e}pend pas de $\omega$. Pour tout $i\xi_j$, en
appliquant la formule de Cauchy sur un petit disque
$D(i\xi_j,\varepsilon)$ contenu dans le domaine d'holomorphie de $F$,
on a alors
$$\left|\dot{F}(i\xi_j)\right|=\left|\frac{1}{2\pi i}
\int_{|k-i\xi_j|=\varepsilon}\frac{F(k)}{\left(k-i\xi_j\right)^2}dk\right|
\leq\frac{\|F\|_{\infty}}{\varepsilon}.$$
On peut prendre $\varepsilon=\frac{\xi_1}{2}$, $\xi_1=\omega\eta_N$ 
{\'e}tant la plus petite valeur propre et qui est $\geq\frac{b}{\omega}$ 
d'apr{\`e}s ce qui pr{\'e}c{\`e}de, ce qui donne, pour tout $j=1,\ldots,N$,
$$\left|\dot{F}(i\xi_j)\right|\leq b\omega\exp\left(a\omega^2\right)\leq
b'\exp\left(a'\omega^2\right),$$
et prouve la majoration cherch{\'e}e.
\bigskip

Pour la minoration, on commence par poser
$$\widetilde{F}(k)=\frac{F(k)}{\prod_{l=1}^N\frac{k-i\xi_l}{k+i\xi_l}}.$$
Alors $\widetilde{F}$ est aussi continue sur le demi-plan ferm{\'e},
holomorphe {\`a} l'int{\'e}rieur, tend vers $1$ {\`a} l'infini, et est sans
z{\'e}ro. En particulier,
\begin{eqnarray*}
\dot{F}(i\xi_j) & = &
\frac{d}{dk}\left(\prod_{l=1}^N\frac{k-i\xi_l}{k+i\xi_l}\right)(i\xi_j)\;\widetilde{F}(i\xi_j)\\
& = & \frac{1}{2i\xi_j}\left(\prod_{l\neq j}\frac{\xi_j-\xi_l}{\xi_j+\xi_l}\right)\widetilde{F}(i\xi_j).
\end{eqnarray*}
On a d'abord $\xi_j\leq\xi_N\leq b\omega$.
Ensuite, pour minorer l'{\'e}cart $|\xi_j-\xi_l|$, bien qu'on puisse
supposer la r{\'e}partition
asymptotique des valeurs propres approximativement uniforme, on va montrer une
estimation plus faible, mais suffisante : $$\forall\,j,\,l,\,1\leq j<l\leq N,\,\eta_j-\eta_l\geq\frac{a}{\omega^2}.$$
En effet, $$\eta_j-\eta_l\geq\min_{1\leq l\leq N-1}(\eta_l-\eta_{l+1}),$$ et
$$\frac{\pi}{\omega}=\Phi(\eta_{l+1})-\Phi(\eta_l)=\int_{x_-(\eta_{l+1})}^{x_+(\eta_{l+1})}\left(Q(y)-\eta_{l+1}^2\right)^{\frac{1}{2}}dy-\int_{x_-(\eta_l)}^{x_+(\eta_l)}\left(Q(y)-\eta_l^2\right)^{\frac{1}{2}}dy.$$
Sur $[x_+(\eta_l),x_+(\eta_{l+1})]$
(resp. $[x_-(\eta_{l+1}),x_-(\eta_l)]$), on a $-\eta_l^2\leq
-Q(y)\leq-\eta_{l+1}^2$, donc
$$Q(y)-\eta_{l+1}^2\leq\eta_l^2-\eta_{l+1}^2,$$
et sur $[x_-(\eta_l),x_+(\eta_l)]$,
$$\left(Q(y)-\eta_{l+1}^2\right)^{\frac{1}{2}}-\left(Q(y)-\eta_l^2\right)^{\frac{1}{2}}\leq
\left(\eta_l^2-\eta_{l+1}^2\right)^{\frac{1}{2}},$$
ce qui donne
\begin{eqnarray*}
\frac{\pi}{\omega} & \leq & (\eta_l-\eta_{l+1})^{\frac{1}{2}}
(\eta_l+\eta_{l+1})^{\frac{1}{2}}(x_+(\eta_{l+1})-x_-(\eta_{l+1}))\\
& \leq & \sqrt{2}\,(\eta_l-\eta_{l+1})^{\frac{1}{2}}(x_+(\eta_{l+1})
-x_-(\eta_{l+1}))
\end{eqnarray*}
Or, $$x_+(\eta_{l+1})-x_-(\eta_{l+1})\leq x_+(\eta_N)-x_-(\eta_N)\leq
a\omega^{\frac{2}{k-2}},$$
et donc $$\eta_l-\eta_{l+1}\geq\frac{a}{\omega^{\frac{k}{k-2}}}\geq\frac{a}{\omega^2}.$$
On en d{\'e}duit
$$\prod_{l\neq
  j}\left|\frac{\xi_j-\xi_l}{\xi_j+\xi_l}\right|\geq\left(\frac{a}{2\omega^2}\right)^{N-1}\geq\frac{a''}{\exp(b''\omega^2)}\,.$$
\bigskip

Il ne reste donc plus qu'{\`a} minorer
$\left|\widetilde{F}(i\xi_j)\right|$ : $\widetilde{F}$ 
{\'e}tant holomorphe sans z{\'e}ro, on a pour tout $R$ assez grand :
$$\frac{1}{4i\xi_j\,\widetilde{F}(i\xi_j)}=\frac{1}{2\pi i}\left(\int_{-R}^R\frac{k}{\widetilde{F}(k)(k-i\xi_j)(k+i\xi_j)^2}dk+\int_0^{\pi}\frac{iR^2e^{2i\theta}}{\widetilde{F}(Re^{i\theta})(Re^{i\theta}-i\xi_j)(Re^{i\theta}+i\xi_j)^2}d\theta\right)\,,$$
qui donne {\`a} la limite (puisque $\widetilde{F}(k)\rightarrow1$)
$$\frac{1}{\widetilde{F}(i\xi_j)}=\frac{4i\xi_j}{2\pi
  i}\int_{-\infty}^{+\infty}\frac{k}{\widetilde{F}(k)(k-i\xi_j)(k+i\xi_j)^2}dk,$$
soit
$$\frac{1}{\left|\widetilde{F}(i\xi_j)\right|}\leq\frac{4\omega}{\pi}\int_0^{\infty}\frac{k}{|F(k)|\,|k-i\xi_j|\,|k+i\xi_j|^2}dk,$$
car, $\forall\,k\in\mathbb{R},\;\left|\widetilde{F}(k)\right|=|F(k)|$,
et $F(-k)=\overline{F(k)}$.

En reprenant les estimations de la solution de Jost, on a pour tout $k>0$ :
\begin{eqnarray*}
|F(k)-1| & \leq & \sum_{n\geq1}|f_n(k,0)|\\
& \leq &
\sum_{n\geq1}\frac{1}{n!}\left(\int_0^{\infty}\frac{2\sqrt{2}}{k}|V(t)|dt\right)^n\\
& \leq & \frac{\omega^2}{k}\int_0^{\infty}2\sqrt{2}\,|Q(t)|dt\;\exp\left(\omega^2\int_0^{\infty}\frac{2\sqrt{2}}{k}|Q(t)|dt\right),
\end{eqnarray*}
ce qui montre que, pour $k\geq C\omega^2$ (o{\`u} $C$
ne d{\'e}pend que de $Q$), $|F(k)|\geq\frac{1}{2}$. Il en r{\'e}sulte que
$$\int_{C\omega^2}^{\infty}\frac{k}{|F(k)|\,|k-i\xi_j|\,|k+i\xi_j|^2}dk\leq2\int_{C\omega^2}^{\infty}\frac{1}{k^2}dk=\frac{2}{C\omega^2}\,.$$

Et pour minorer $\frac{1}{|F(k)|}$ sur $\left]0,C\omega^2\right]$, on
utilise le 
wronskien de $f(k,\cdot)$ et $f(-k,\cdot)$, qui est constant et qui
vaut $-2ik$ (on le calcule en prenant
$\lim_{x\rightarrow\infty}$, sachant que $f'(k,x)\sim ike^{ikx}$). Il
vaut par ailleurs, pour tout $k>0$,
$$f(k,0)f'(-k,0)-f(-k,0)f'(k,0)=2i\Im
m\left(F(k)\overline{f'(-k,0)}\right).$$
A l'aide de l'{\'e}quation int{\'e}grale v{\'e}rifi{\'e}e par $f(k,\cdot)$,
$$f(k,x)=e^{ikx}+\int_x^{\infty}\frac{\sin k(x-t)}{k}V(t)f(k,t)dt,$$
on en d{\'e}duit, $\forall\,k,\;0<k\leq C\omega^2$,
\begin{eqnarray*}
|f'(k,0)| & \leq & 
k+\int_0^{\infty}|\cos kt|\,|V(t)|\,|f(k,t)|dt\\ 
& \leq & C\omega^2+\omega^2\exp\left(a\omega^2\right)\int_0^{\infty}|Q(t)|dt\\
& \leq & C'\exp\left(a'\omega^2\right),
\end{eqnarray*}
ainsi que pour $f'(-k,0)=\overline{f'(k,0)}$. On a alors
$$2k\leq2|F(k)|C'\exp\left(a'\omega^2\right),$$
ce qui permet d'obtenir
\begin{eqnarray*}
\int_0^{C\omega^2}\frac{k}{|F(k)|\,|k-i\xi_j|\,|k+i\xi_j|^2}dk & \leq
& C'\exp\left(a'\omega^2\right)\int_0^{C\omega^2}\frac{1}{\xi_j^3}dk\\
& \leq & C'\exp\left(a'\omega^2\right)(b\omega)^3C\omega^2\\
& \leq & C''\exp\left(a''\omega^2\right),
\end{eqnarray*}
et aboutit finalement {\`a}
$$\frac{1}{\left|\widetilde{F}(i\xi_j)\right|}\leq
\frac{4\omega}{\pi}\left(C''\exp\left(a''\omega^2\right)
+\frac{2}{C\omega^2}\right),$$
ce qui termine la preuve. 

\end{proof} 
\bigskip 

Avant de traiter le th{\'e}or{\`e}me de cette partie, on consid{\`e}re deux 
exemples o{\`u} on peut donner des estimations explicites. 

\begin{example} 

Soit $$Q_1(x)=\frac{1}{\left(1+x^2\right)^2}.$$
Alors $Q_1$ est paire, et v{\'e}rifie les conditions de la
proposition~\ref{estimation}. En appliquant la m{\'e}thode WKB, on a (pour
$\omega\geq1$) : 
$$N(\omega)=\left[\frac{1}{\varepsilon}\right]=[\omega].$$
On trouve $x_+(\eta)=-x_-(\eta)=\sqrt{\frac{1}{\eta}-1}$, et
$$\frac{\pi^2}{256}\,\frac{1}{\omega^2}\leq\eta_N\leq\frac{\pi^2}{16}
\frac{1}{\omega^2},$$
soit, $\forall\,n,\,1\leq n\leq N$, $$\frac{\pi^2}{256}\frac{1}{\omega}\leq\xi_n\leq\omega.$$
Pour l'{\'e}cart des valeurs propres, on trouve
$$\xi_n-\xi_{n-1}\geq\frac{1}{5\omega}.$$
On a finalement, $\forall\,n,\,1\leq n\leq N$,
$$1\leq s_n\leq\exp\left(\frac{4}{\pi}\omega^2+\frac{\pi}{2}\omega\right),$$
et
$$\frac{1}{\exp\left(2\omega\ln4\omega\right)}\leq\prod_{j\neq
  n}\left(\frac{\xi_n-\xi_j}{\xi_n+\xi_j}\right)^2\leq1,$$
puis
$$\left|\dot{F}(i\xi_j)\right|\leq\frac{512}{\pi^2}\omega\exp\left(\pi\sqrt{2}\,\omega^2\right),$$
et
$$\frac{1}{\left|\widetilde{F}(i\xi_j)\right|}\leq2^{22}\omega^7\exp\left(\pi\sqrt{2}\,\omega^2\right),$$
donc
$$\frac{1}{45\,\omega^9\exp\left(26\,\omega^2\right)}\leq
\frac{4\xi_n^2}{C_n}
\leq22\,\omega^2\exp\left(15\,\omega^2\right)\,.$$

\end{example}
\bigskip

\begin{example}

Soit
$$Q_2(x)=
\begin{cases}
1,\text{ si } x\in[0,1],\\
0,\text{ si } x>1.
\end{cases}
$$
Le potentiel associ{\'e} s'av{\`e}re {\^e}tre un cas limite vu sa discontinuit{\'e} en
$x=1$. L'exemple se traite ici directement.

L'equation $-y''-\omega^2y=\lambda y,\;\lambda=-\xi^2,\,0<\xi<\omega$,
admet pour fonctions propres les
$$
y_{\xi}(x)=
\begin{cases}
\sqrt{\frac{2\xi}{1+\xi}}\sin\left(\sqrt{\omega^2-\xi^2}\,x\right),\text{ si
}x\in[0,1],\\
\sqrt{\frac{2\xi}{1+\xi}}\,e^{\xi}\sin\sqrt{\omega^2-\xi^2}\,e^{-\xi
  x},\text{ si }x\geq1,
\end{cases}
$$
o{\`u} $\xi$ satisfait l'{\'e}quation
$$\xi\sin\sqrt{\omega^2-\xi^2}+\sqrt{\omega^2-\xi^2}\cos\sqrt{\omega^2-\xi^2}=0.\pod\diamond$$
Dans ce cas, les $y_{\xi}\in C^1(\mathbb{R}^+)$, v{\'e}rifient les
conditions aux bords en $0$ et $+\infty$, et sont norm{\'e}es :
$$\int_0^{+\infty}y_{\xi}^2(x)dx=1.$$
On sait d'abord que $\xi=O(\omega)$. Et puisque
$$C_{\xi}=\left(y_{\xi}'(0)\right)^2=\frac{2\xi}{1+\xi}\left(\omega^2-\xi^2\right),$$
alors
$$\frac{4\xi^2}{C_{\xi}}=\frac{2\xi(\xi+1)}{\omega^2-\xi^2}.$$
Pour la premi{\`e}re estimation, on a m{\^e}me
$0<\frac{4\xi^2}{C_{\xi}}\leq220\,\omega^2$ : en effet, si
$\omega\geq10,\,\varepsilon=\sqrt{\omega^2-\xi^2}\leq\frac{1}{10}$,
l'{\'e}quation$\pod\diamond$ devient
$$\sqrt{\omega^2-\varepsilon^2}\sin\varepsilon+\varepsilon\cos\varepsilon\geq
10\,\varepsilon>0,$$
ce qui impose
$\varepsilon\geq\frac{1}{10}$, soit
$$\xi\leq\sqrt{\frac{99}{100}}\,\omega.$$

En revanche, la minoration de $\frac{4\xi^2}{C_{\xi}}$ n{\'e}cessite une
minoration de $\xi$ du m{\^e}me type, qui est fausse en g{\'e}n{\'e}ral, car la
premi{\`e}re valeur propre $\xi_1$ ne peut pas {\^e}tre minor{\'e}e par
$\frac{1}{\omega^k}$ ou $\exp(-a\omega^b)$, ce qui est n{\'e}cessaire pour
avoir $\ln\frac{4\xi^2}{C_{\xi}}=O(\omega^b)$.

En effet, fixons par exemple $\omega_0$ tr{\`e}s grand,
$\equiv\frac{\pi}{2}\pmod{2\pi}$, et consid{\'e}rons l'{\'e}quation
$\pod\diamond$ en $(\xi,\omega)$, qui est donn{\'e}e par une fonction $g$
suffisamment r{\'e}guli{\`e}re au voisinage de $(0,\omega_0)$. Puisque
$\frac{\partial g}{\partial\xi}(0,\omega_0)=\sin\omega_0=1$, et
$\frac{\partial g}{\partial\omega}(0,\omega_0)=-\omega_0$, une
application imm{\'e}diate du th{\'e}or{\`e}me des fonctions implicites permet de
consid{\'e}rer la fonction $\xi(\omega)$, qui admet au voisinage de
$\omega_0$ le d{\'e}veloppement
$$\xi(\omega)=\omega_0(\omega-\omega_0)+O\left((\omega-\omega_0)^2\right).$$
On voit alors que,
$\forall\,\omega,\,0<\omega-\omega_0\leq\eta(\omega_0),\;\eta(\omega_0)$
assez petit, on a $$0<\xi(\omega)\leq2\exp\left(-\frac{1}{2}\exp\omega\right),$$
ce qui montre, pour tout $\omega$ arbitrairement grand,
l'impossibilit{\'e} de minorer $\xi$ comme voulu. Ce probl{\`e}me provient du
fait que $Q_2$ n'est m{\^e}me pas continue en $x=1$.
\bigskip

Par contre, si on impose une restriction sur $\omega\geq10$ du type
$$\left|\omega-\frac{\pi}{2}+\pi\mathbb{Z}\right|\geq\frac{1}{5},$$
on voit que, pour tout $\xi,\,0<\xi\leq\frac{1}{10}$, on a 
$\left|\sqrt{\omega^2-\xi^2}-\frac{\pi}{2}+\pi\mathbb{Z}\right|\geq\frac{1}{10}$, et
$$\left|\xi\sin\sqrt{\omega^2-\xi^2}+\sqrt{\omega^2-\xi^2}\cos\sqrt{\omega^2-\xi^2}\right|\geq\frac{1}{2}.$$
Ainsi, on doit avoir $\xi\geq\frac{1}{10}$, ce qui donne
$$\frac{4\xi^2}{C_{\xi}}\geq\frac{1}{5\omega^2},$$
et aboutit {\`a} un encadrement m{\^e}me plus fin de $\frac{4\xi^2}{C_{\xi}}$
en $O(\omega^2)$.

\end{example}

\subsection{Applications au probl{\`e}me inverse}

On peut finalement {\'e}noncer les cons{\'e}quences des r{\'e}sultats n{\'e}gatifs {\'e}tablis
dans la partie~\ref{appranal} et des r{\'e}sultats positifs {\'e}nonc{\'e}s
pr{\'e}c{\'e}demment. 

\begin{theorem}\label{optim1}

Soit $\mathcal{Q}$ la classe des fonctions $Q$ d{\'e}finies sur
$\mathbb{R}^+$, qui sont strictement positives, strictement
d{\'e}croissantes avec d{\'e}croissance polynomiale {\`a} l'infini, et qui ont $2$
d{\'e}riv{\'e}es localement int{\'e}grables tendant polynomialement vers $0$,
avec $Q'(0)=0$, et soit $\Lambda_{\mathcal{Q}}$ un 
compact de $\mathcal{Q}$, du type $\Lambda_l$. Consid{\'e}rons, pour tout
$\omega$ assez grand, la classe des op{\'e}rateurs de Sturm-Liouville
$-\frac{d^2}{dx^2}-\omega^2Q$, o{\`u} $Q\in\Lambda_{\mathcal{Q}}$.

Donnons-nous {\'e}galement, pour tous $N$, 
$a_1\omega\leq N\leq a_2\omega$, et $b(\omega)>0$, 
$\frac{1}{b(\omega)}$ polynomiaux en $\omega$, une fonction 
$\psi(x,\zeta)$ d{\'e}finie sur 
$\mathbb{R}\times\mathbb{C}^N$, de classe $C^1$ par rapport {\`a} $x$ et
v{\'e}rifiant les conditions du corollaire~\ref{inverse2} sur tout
$[0,X]$. Alors l'approximation de
$$\int_0^{\cdot}\Lambda_{\mathcal{Q}}:=\left\{\left(x\mapsto\int_0^xQ(t)dt
\right),\;Q\in\Lambda_{\mathcal{Q}}\right\},$$
au sens uniforme sur tout intervalle $[0,X]$, $X\geq1$, par la famille
$$\left\{\left(x\mapsto\frac{1}{b(\omega)}\left(\frac{1}{\psi}
\frac{\partial\psi}{\partial x}
\right)(x,\zeta)\right), \,\zeta_j=O(\omega^r),
\,\forall\,j=1,\ldots,N\right\},$$
lorsque $\omega\rightarrow\infty$, ne peut pas {\^e}tre meilleure que de
l'ordre de $$\frac{1}{(\omega\ln\omega)^3}\,.$$

On dispose d'autre part d'une formule d'approximation telle que, si
$N(\omega)$ est le nombre de 
valeurs propres $\xi_j$ et caract{\'e}ristiques $C_j$ de l'op{\'e}rateur
$-\frac{d^2}{dx^2}-\omega^2Q$, et
$$\Psi(x,\zeta)=\det\widetilde{W}_{s,r}(x,\zeta),$$
avec
$$\widetilde{W}_{s,r}(x,\zeta)=\frac{2sh(\zeta_r+\zeta_s)x}{\zeta_r+\zeta_s}
-(1-\delta_{s,r})\frac{2sh(\zeta_s-\zeta_r)x}{\zeta_s-\zeta_r}
-\delta_{s,r}(2x-\exp(\zeta_{r+N})),$$
$s,\,r=1,\ldots,N(\omega)$,
alors la famille 
$\left\{\frac{2}{\omega^2}\frac{1}{\Psi}\frac{\partial\Psi}{\partial x}
\right\}$
approche le compact $\int_0^{\cdot}\Lambda_{\mathcal{Q}}$ au moins {\`a} l'ordre
$\frac{\ln\omega}{\sqrt{\omega}}$.

En outre, si $Q$ est donn{\'e}e, un {\'e}l{\'e}ment
$\zeta(Q)$ optimisant peut {\^e}tre ainsi choisi :
$$\zeta_j(Q)=\xi_j, \text{ et } \zeta_{j+N}(Q)=\ln\frac{4\xi_j^2}{C_j},
\text{ } j=1,\ldots,N(\omega).$$

\end{theorem}

Avant de prouver ce th{\'e}or{\`e}me, on tient {\`a} en pr{\'e}ciser l'interpr{\'e}tation
en terme de probl{\`e}me inverse : il n'existe pas de formule
qui puisse donner analytiquement une approximation meilleure que de
l'ordre de 
$\frac{1}{(\omega\ln\omega)^3}$, de tout potentiel (avec deux
d{\'e}riv{\'e}es) {\`a} partir de ses valeurs propres et valeurs caract{\'e}ristiques.
Ainsi, la formule d'approximation (explicit{\'e}e) de Gelfand-Levitan-Jost-Kohn
donne un r{\'e}sultat positif avec une vitesse au minimum de l'ordre de
$\frac{\ln\omega}{\sqrt{\omega}}$, 
ce qui remplit notre objectif pour le cas de $2$ d{\'e}riv{\'e}es et r{\'e}pond {\`a} la
question pos{\'e}e dans~\cite{henkin} p. 22,
sur le probl{\`e}me d'approximation de $C^{m+1}$ par une famille
{\`a} $\omega$ param{\`e}tres.

Il est d'autre part int{\'e}ressant de constater que le choix de cette fonction
optimisante n'a pas {\'e}t{\'e} construite sp{\'e}cialement
dans le cadre de la th{\'e}orie d'approximation, puisqu'elle provient de
la th{\'e}orie physique math{\'e}matique.

\begin{proof}

Il est d'abord {\`a} noter que la restriction sur $[0,1]$ d'une
telle fonction $Q$ est bien dans 
$\Lambda_2([0,1])$ (du moins dans un homoth{\'e}tique), et
qu'inversement si on se donne $h\in\Lambda_2\cap C^2([0,1])$,
$h$ positif, strictement d{\'e}croissant avec $h'(0)=0$, on 
peut le prolonger sur $\mathbb{R}^+$ en 
une fonction $Q_h$ (soit un potentiel $-\omega^2Q_h$) de classe $C^2$,
positive avec d{\'e}croissance polynomiale (ainsi que ses deux d{\'e}riv{\'e}es).

Cela permet, {\'e}tant donn{\'e}s $\psi$ et le compact $\Lambda_{\mathcal{Q}}$,
d'appliquer le corollaire~\ref{detail} (avec $m=2$), $N$ et $\omega$
ayant le m{\^e}me ordre de grandeur. Il existe donc un potentiel
$-\omega^2Q_h$, $Q_h>0$, tel que $\int_0^{\cdot}Q_h$ soit 
distant de
$\left\{\frac{1}{b(\omega)}\frac{1}{\psi}\frac{\partial\psi}
{\partial x}(\cdot,\zeta)\right\}$ sur
$[0,1]$ d'au moins 
$\frac{C}{(2N\ln(2N))^3}$, pour tout $\zeta$ de taille polynomiale en
$\omega$, soit de l'ordre de $\frac{1}{(\omega\ln\omega)^3}$ sur tout
$[0,X]$. 
\bigskip

Quant {\`a} la fonction $\Psi$ ainsi d{\'e}finie, elle v{\'e}rifie bien les
conditions du corollaire~\ref{inverse}. En effet, c'est d'abord une
fonction {\`a} $2N(\omega)$ param{\`e}tres ($N(\omega)$ et $\omega$ ayant le
m{\^e}me ordre de grandeur), enti{\`e}re de type exponentiel :

c'est le cas pour
$$\left\|\frac{sh(\zeta_s\pm\zeta_r)x}{\zeta_s\pm\zeta_r}\right\|_{\infty}\leq\sum_{n\geq
  0}\frac{1}{(2n+1)!}|\zeta_s\pm\zeta_r|^{2n}X^{2n+1}\leq X\exp[(|\zeta_r|+|\zeta_s|)X],$$
donc
\begin{eqnarray*}
\left\|\widetilde{W}_{s,r}\right\|_{\infty} & \leq & 4X\exp[(|\zeta_r|+|\zeta_s|)X]+2X+\exp(|\zeta_{s+N}|+|\zeta_{r+N}|),\\
& \leq & 7X\exp[(|\zeta_r|+|\zeta_{r+N}|+|\zeta_s|+|\zeta_{s+N}|)X],
\end{eqnarray*}
ainsi que pour chaque produit du d{\'e}terminant :
\begin{eqnarray*}
\left\|\prod_{j=1}^N\widetilde{W}_{j,\tau(j)}\right\|_{\infty} & \leq
&
(7X)^N\exp\left[\sum_{j=1}^N(|\zeta_j|+|\zeta_{j+N})X+\sum_{j=1}^N(|\zeta_{\tau(j)}|+|\zeta_{\tau(j)+N}|)X\right]\\
& = & (7X)^N\exp(2X\|\zeta\|_1).
\end{eqnarray*}
Comme ils sont au nombre de
$N!=O\left(e^{N^2}\right)$, on obtient une estimation de $\Psi$ en
$O\left(e^{2N^2}e^{4X\|\zeta\|_1}\right)$.

Pour les estimations de $\frac{\partial\Psi}{\partial x}$, il suffit
de remarquer que $\Psi$ est 
{\'e}galement enti{\`e}re (et m{\^e}me de type exponentiel) par rapport {\`a} la variable
$x$, comme on le voit pour chaque $\widetilde{W}_{s,r}$ (et donc pour
le d{\'e}terminant). La formule de Cauchy appliqu{\'e}e sur le disque
$D(0,X+1)$ et l'estimation pr{\'e}c{\'e}dente nous donnent des majorations
analogues pour $\frac{\partial\Psi}{\partial x}$ sur $[0,X]$,
$\forall\zeta\in\mathbb{C}^{2N}$.
\bigskip

Ensuite,
$$\det\left(\widetilde{W}_{s,r}\right)(0)^{-1}=
\frac{1}{2}\prod_{j=1}^N\exp(-\zeta_{N+j})\leq\exp(\|\zeta\|_1)=
O\left(e^{\alpha\omega^{\beta}}\right).$$
\bigskip

Enfin, le choix des {\'e}valuations
$$\zeta_j(Q)=\xi_j\left(-\omega^2Q\right),\text{ et }\zeta_{j+N}(Q)=
\ln\frac{4\xi_j^2}{C_j}\left(-\omega^2Q\right),\;j=1,\ldots,N,$$
est possible, car chaque $Q\in\mathcal{Q}$ v{\'e}rifiant les conditions de la
proposition~\ref{estimation}, on a
$\xi_j\left(-\omega^2Q\right)=O(\omega)$, et
$$\frac{1}{\alpha\exp(\beta\omega^{\gamma})}\leq
\frac{4\xi_j^2}{C_j}\left(-\omega^2Q\right)\leq
\alpha\exp(\beta\omega^{\gamma}),$$
donc
$\ln\frac{4\xi_j^2}{C_j}\left(-\omega^2Q\right)=
O\left(\omega^{\gamma}\right)$.

Et ce choix de $\zeta(Q)$ donne l'approximation voulue de
$\int_0^{\cdot}Q$ sur tout $[0,X]$ {\`a} l'ordre de 
$\frac{\ln\omega}{\sqrt{\omega}}$ :
cela provient en effet du th{\'e}or{\`e}me 2 donn{\'e} 
dans~\cite{henkin} p. 22, qui nous dit que, uniform{\'e}ment sur tout
$[0,X]$,
\begin{eqnarray*}
\left|\int_0^xQ(y)dy-\frac{1}{\omega^2}\frac{1}{\Psi(x)}\frac{\partial\Psi}
{\partial x}(x,\zeta(Q))\right| & = &
\left|\int_0^xQ(y)dy-\frac{2}{\omega^2}\int_0^x 
\frac{\partial^2}{\partial x^2}\ln\left|\Psi(y,\zeta(Q))\right|dy\right|\\
& = & \left|\int_0^xQ(y)dy-2\int_0^xQ_{\omega}^0(y)dy\right|\\ 
& = & O\left(\frac{\ln\omega}{\sqrt{\omega}}\right)\,.
\end{eqnarray*}
Les conditions de r{\'e}gularit{\'e} pour $Q$ sont v{\'e}rifi{\'e}es : $Q$ est
strictement positif (m{\^e}me minor{\'e} par une borne ne d{\'e}pendant que du
compact $\Lambda_{\mathcal{Q}}$ ) et dans (un homoth{\'e}tique de)
$\Lambda_2$, hypoth{\`e}se qui peut remplacer 
celle du nombre fini d'intervalles de monotonie de ses d{\'e}riv{\'e}es.

\end{proof}

\begin{remark}

La fonction $\Psi$ v{\'e}rifie {\'e}galement,
$\forall\,\zeta\in\mathbb{C}^N$,
$$\left(\frac{\partial}{\partial x}\det\widetilde{W}_{s,r}\right)
(0,\zeta)=0.$$
En effet, $\forall\,r,\,s=1,\ldots,N$
$$\frac{\partial\widetilde{W}_{s,r}}{\partial x}(x,\zeta)=
2ch(\zeta_r+\zeta_s)x-(1-\delta_{s,r})
2ch(\zeta_s-\zeta_r)x-2\delta_{s,r}\,,$$
donc
$$\frac{\partial\widetilde{W}_{s,r}}{\partial x}(0,\zeta)=
2-2(1-\delta_{s,r})-2\delta_{s,r}=0.$$
C'est une hypoth{\`e}se inutile pour ce th{\'e}or{\`e}me, mais
qui peut nous servir si on veut utiliser $\Psi$ dans le cadre du
corollaire~\ref{inverse}. 

\end{remark}
\bigskip

Si on consid{\`e}re maintenant le cas plus g{\'e}n{\'e}ral d'un potentiel
$-\omega^2Q$ avec
$m+1$ d{\'e}riv{\'e}es localement int{\'e}grables et qui s'annulent en $0$, on est
au moins capable d'expliciter une mesure spectrale voisine de celle de
l'op{\'e}rateur $-\frac{d^2}{dx^2}-\omega^2Q$, qui est
$$\sigma_{\omega}(d\tau)=
\begin{cases}
\frac{1}{\pi}\sqrt{\tau+\omega^2Q(0)},\;\tau\geq0,\\
\sum_{j=1}^{N(\omega)}C_j\delta(\tau+\xi_j^2),\;\tau<0\,.
\end{cases}
$$
Posons alors, pour $0\leq y\leq x\leq X$,
$$\Phi(x,y)=\frac{1}{\pi}\int_0^{\infty}\frac{\sin(x\sqrt{\tau})}
{\sqrt{\tau}}\frac{\sin(y\sqrt{\tau})}{\sqrt{\tau}}\frac{\omega^2Q(0)
\;d\tau}{\sqrt{\tau+\omega^2Q(0)}+\sqrt{\tau}}\,,$$
(int{\'e}grale absolument convergente pour tous $x,\,y$) et consid{\'e}rons le
noyau $K(x,y)$ solution de l'{\'e}quation int{\'e}grale
$$K(x,y)+\int_0^xK(x,s)\Phi(s,y)ds+\Phi(x,y)\equiv0\,.$$
Construisons alors le potentiel $q_{\omega}$ de la fa{\c c}on suivante :
$$q_{\omega}(x)=2\frac{d}{dx}K(x,x)-2\frac{d^2}{dx^2}\ln|\det T(x)|\,,$$
o{\`u} $T(x)$ est la matrice d'ordre $N(\omega)$ d{\'e}finie par (cf
~\cite{levitan}) 
\begin{eqnarray*}
\lefteqn{T_{j,k}(x)=\frac{4\xi_j^2}{C_j}\delta_{j,k}\,+}\\
& & +\,4\int_0^x\left(sh(\xi_jt)+\int_0^tK(t,s)sh(\xi_js)ds\right)
\left(sh(\xi_kt)+\int_0^tK(t,s)sh(\xi_ks)ds\right)\,dt\,.
\end{eqnarray*}
Alors (~\cite{henkin}, th{\'e}or{\`e}me 1 p. 21) la fonction
$Q_{\omega}=-\frac{q_{\omega}}{\omega^2}$ r{\'e}alise, uniform{\'e}ment sur tout 
compact $[0,X]$, une approximation de $Q$ (au moins) {\`a} l'ordre
$\frac{1}{\omega^m}$.
\bigskip

Comme on le voit, $Q_{\omega}$ s'{\'e}crit comme une fonction 
analytique $\widetilde{Q}_{\omega}$ en les variables $\zeta,\,w$, o{\`u}
$\zeta\in\mathbb{C}^{2N},\,\Re e\,w>0$ ($w$ remplace la variable
$\omega^2Q(0)$ ), soit
$$\widetilde{Q}_{\omega}(x,\zeta,w)=\frac{2}{\omega^2}\left(
-\frac{\partial\widetilde{K}}{\partial x}(x,x,w)+
\frac{\partial^2}{\partial x^2}\ln\left|\det\widetilde{T}(x,\zeta,w)
\right|\right)\,,$$
o{\`u} $\widetilde{K}$ est la solution, pour tous 
$0\leq y\leq x\leq X,\,\Re e\,w>0$, de l'{\'e}quation
$$\widetilde{K}(x,y,w)+\int_0^x\widetilde{K}(x,s,w)
\widetilde{\Phi}(s,y,w)ds+\widetilde{\Phi}(x,y,w)\equiv0\,,$$
avec
$$\widetilde{\Phi}(x,y,w)=\frac{1}{\pi}\int_0^{\infty}
\frac{\sin(x\sqrt{\tau})}{\sqrt{\tau}}\frac{\sin(y\sqrt{\tau})}{\sqrt{\tau}}
\frac{w\;d\tau}{\sqrt{\tau+w}+\sqrt{\tau}}\,,$$
($\widetilde{\Phi}$ est bien d{\'e}finie puisque $\Re e\,w>0$) et
\begin{eqnarray*}
\lefteqn{\widetilde{T}_{j,k}(x,w,\zeta)=\exp(\zeta_{N+j})\delta_{j,k}
\,+}\\
& & +\,4\int_0^x\left(sh(\zeta_jt)+\int_0^t\widetilde{K}(t,s,w)
sh(\zeta_js)ds\right)\left(sh(\zeta_kt)+\int_0^t\widetilde{K}(t,s,w)
sh(\zeta_ks)ds\right)\,dt\,.
\end{eqnarray*}

Nous esp{\'e}rons avoir prouv{\'e} que les fonctions
$\frac{\partial\widetilde{K}}{\partial x}$ et
$\frac{\partial\widetilde{K}}{\partial y}$ soient holomorphes
de type exponentiel par rapport {\`a} $w\in\left\{\Re e\,z>0\right\}$,
et de restriction sur $\mathbb{R}^+$ de type polynomial. 
On pourrait en d{\'e}duire un r{\'e}sultat de presque optimalit{\'e} analogue au
th{\'e}or{\`e}me~\ref{optim1}, qui se d{\'e}duirait des corollaires~\ref{inverse}
et~\ref{detail} : 
si $\mathcal{Q}$ est la classe des fonctions $Q$ d{\'e}finies sur
$\mathbb{R}^+$, strictement positives, {\`a} d{\'e}croissance polynomiale,
avec $m+1$ 
d{\'e}riv{\'e}es localement int{\'e}grables qui s'annulent en $0$, et
$\psi(x,\zeta,w)$ une fonction d{\'e}finie sur 
$\mathbb{R}\times\mathbb{C}^N\times\left\{\Re e\,w>0\right\}$,
continue par rapport {\`a} $x$ et 
v{\'e}rifiant les conditions du corollaire~\ref{inverse}, alors
l'approximation de tout $\Lambda_{\mathcal{Q}}$ sur tout
$[0,X]$, par la famille 
$$\left\{\left(x\mapsto\frac{1}{b(\omega)}\frac{\partial}{\partial x}
\left(\frac{1}{\psi}\frac
{\partial\psi}{\partial x}\right)(x,\zeta,w)\right),\,
\zeta_j=O(\omega^r),\,|w-a(\omega)|\leq\frac{a(\omega)}{2}\right\},$$
o{\`u} $b(\omega)>0$, $\frac{1}{b(\omega)}$ et $a(\omega)\geq1$ sont
polynomiaux en $\omega$, ne peut pas {\^e}tre 
meilleure que de l'ordre de $$\frac{1}{(\omega\ln\omega)^{m+1}}.$$
En outre, la fonction d{\'e}finie pr{\'e}c{\'e}demment
$$\widetilde{\Psi}(x,\zeta,w)=\exp\left(
-\int_0^x\widetilde{K}(t,t,w)dt\right)\det\widetilde{T}
(x,\zeta,w)\,,$$
$j,\,k=1,\ldots,N(\omega)$,
r{\'e}aliserait un cas d'approximation presque optimale {\`a} l'ordre
$\frac{1}{\omega^m}$, avec comme choix optimisant l'{\'e}l{\'e}ment
$$\zeta_j(Q)=\xi_j,\;\zeta_{j+N}(Q)=\ln\frac{4\xi_j^2}{C_j},
\;j=1,\ldots,N(\omega),$$
et
$$w(Q)=a(\omega)=\omega^2Q(0).$$

\begin{remark}

Lorsque $Q$ s'annule en $0$ avec d{\'e}riv{\'e}es, la mesure spectrale devient
$$\sigma^0_{\omega}(d\tau)=
\begin{cases}
\frac{1}{\pi}\sqrt{\tau},\;\tau\geq0,\\
\sum_{j=1}^{N(\omega)}C_j\delta(\tau+\xi_j^2),\;\tau<0\,,
\end{cases}
$$
qui redonne la formule d'approximation de type Gelfand-Levitan
$Q_{\omega}^0$, 
qui est compl{\`e}tement explicit{\'e}e. Le seul probl{\`e}me est qu'on ne sait
pas si la vitesse d'approximation est toujours de l'ordre de
$\frac{1}{\omega^m}$, car pour utiliser le th{\'e}or{\`e}me 1
dans~\cite{henkin}, il faut que la fonction $Q$ soit strictement
positive sur $\mathbb{R}^+$ (en particulier en $0$). On est cependant
assez optimiste, car en pratique, m{\^e}mes pour des potentiels pas tr{\`e}s
r{\'e}gulier, on utilise $Q_{\omega}^0$, qui approxime $Q$ avec une bonne
vitesse : les exemples d'applications num{\'e}riques, m{\^e}me s'ils ne sont
pas des preuves, donnent malgr{\'e} tout un bon pronostic
(cf~\cite{henkin}, partie 4).

\end{remark}

\subsection{D'autres applications aux probl{\`e}mes inverses}

On formule quelques exemples de probl{\`e}mes inverses o{\`u} nous esp{\'e}rons pouvoir
appliquer nos r{\'e}sultats.

\begin{example1}

Le th{\'e}or{\`e}me 1.2 p. 260 dans~\cite{lax-levermore} nous donne un
r{\'e}sultat original dans le cadre d'approximation au sens $L^2$ : si $u$
est un potentiel n{\'e}gatif de classe $C^1$, alors
$\lim_{\varepsilon\rightarrow0}u(\cdot,\varepsilon)=u$, o{\`u}
$$u(x,\varepsilon)=-2\varepsilon^2\frac{d^2}{dx^2}\ln\det(I+
G(x,\varepsilon)),$$
avec
$$G(x,\varepsilon)=\varepsilon\left(\frac{\exp\left(-
\frac{\eta_j+\eta_k}{\varepsilon}x\right)}{\eta_j+\eta_k}C_jC_k\right)\,,$$
$1,\leq j,\,k\leq
N(\varepsilon),\;\varepsilon=\frac{1}{\omega}$. Ici aussi on a affaire
{\`a} une fonction analytique en les valeurs propres $\eta_j$ et $\ln
C_j$, de type exponentiel. Plus pr{\'e}cis{\'e}ment, elle s'{\'e}crit
$$\widetilde{G}(x,\zeta,w)=\varepsilon\frac{\exp
(w_j+w_k)x}{w_j+w_k}\exp(\zeta_j+\zeta_k),$$
o{\`u} on choisira $w_j=\frac{\eta_j}{\varepsilon^r}$ et $\zeta_k=\ln C_k,\;
j=1,\ldots,N(\varepsilon)$, et $r$ assez grand pour minorer
$\frac{\eta_j}{\varepsilon^r}$ par une constante positive. On
esp{\`e}re en d{\'e}duire un r{\'e}sultat analogue au th{\'e}or{\`e}me~\ref{optim1},
avec r{\'e}sultats n{\'e}gatif et positif. Cependant, il faudrait appliquer un
r{\'e}sultat du type corollaire~\ref{inverse} en se passant de l'hypoth{\`e}se que 
$\frac{d}{dx}\det(I+G)(0,\varepsilon)=0$.

\end{example1}

\begin{example1}

On a dans le m{\^e}me esprit un r{\'e}sultat d'approximation de
V. A. Marchenko (cf.~\cite{marchenko}) qui nous dit que, pour un
potentiel $q$ qui a $m+1$ d{\'e}riv{\'e}es sur $\mathbb{R}^+$, la
connaissance de la fonction de Weyl $j(k)$ sur l'intervalle $[-A,A]$
permet de reconstruire $q$ avec une vitesse (au moins) de l'ordre de
$\frac{1}{A^m}$. On peut ici aussi pr{\'e}voir un r{\'e}sultat n{\'e}gatif qui
donnerait une minoration de l'ordre de $\frac{1}{(A\ln A)^{m+1}}$ pour
l'approximation de $q$, malgr{\'e} la difficult{\'e} provenant du passage du
discret au continu pour la connaissance de $j(k)$.

\end{example1}

\begin{example1}

Si on s'int{\'e}resse cette fois {\`a} 
des espaces fonctionnels {\`a} plusieurs variables (et plus seulement
$\mathbb{R}^+$), on peut consid{\'e}rer un th{\'e}or{\`e}me de
R. Novikov (cf~\cite{novikov}), qui est un r{\'e}sultat d'approximation,
dans le cas n=$3$, de fonctions avec $l$ d{\'e}riv{\'e}es int{\'e}grables, {\`a}
partir de la connaissance de l'amplitude de diffusion {\`a} une {\'e}nergie
$E$ donn{\'e}e. On a ici une vitesse de l'ordre de
$$O\left(\frac{1}{E^{\frac{l-3-\varepsilon}{2}}}\right),$$
au sens uniforme quand $E\rightarrow+\infty$ ($\varepsilon$
arbitrairement petit), soit un comportement voisin de
$\left(\frac{1}{\sqrt{E}}\right)^{l-3},\;\sqrt{E}$ {\'e}tant homog{\`e}ne au
param{\`e}tre $\omega$. On peut de m{\^e}me pr{\'e}voir un r{\'e}sultat n{\'e}gatif qui
montrerait une impossibilit{\'e} d'approximer tout potentiel, mieux qu'{\`a} l'ordre
$\left(\frac{1}{\sqrt{E}\,\ln E}\right)^l$.

\end{example1}

\section{Autres m{\'e}thodes}
\bigskip

Pour terminer, avec le th{\'e}or{\`e}me~\ref{th1}, on a montr{\'e} une impossibilit{\'e} de
bien approximer analytiquement des compacts d'espaces
fonctionnels. Comme on l'a vu, ce
r{\'e}sultat utilise le th{\'e}or{\`e}me de Vitushkin pour le cas polynomial, et
le fait qu'une fonction enti{\`e}re de type exponentiel est tr{\`e}s bien
approch{\'e}e par des polyn{\^o}mes. Le probl{\`e}me est que tout doit converger :
ainsi, en plus des param{\`e}tres qui d{\'e}finissent la classe de telles
fonctions, il faut borner les variables ; ce qui a pour cons{\'e}quence,
dans le cadre du probl{\`e}me inverse, la n{\'e}cessit{\'e} d'estimer les valeurs
propres et valeurs caract{\'e}ristiques d'un op{\'e}rateur de Sturm-Liouville.
\bigskip

Comme on l'a signal{\'e} dans l'introduction on pouvait s'inspirer de la
m{\'e}thode de Warren en essayant d'estimer le nombre de composantes
connexes de l'ensemble des z{\'e}ros d'une fonction analytique. De m{\^e}me
que l'estimation de Warren donn{\'e}e dans~\cite{warren} remontait au
th{\'e}or{\`e}me de B{\'e}zout, il s'agirait ici d'{\'e}tablir des estimations du
nombres de solutions non d{\'e}g{\'e}n{\'e}r{\'e}es d'un syst{\`e}me de fonctions
analytiques. On aboutirait ainsi {\`a} des r{\'e}sultats n{\'e}gatifs analogues,
mais cette fois-ci sans avoir {\`a} borner les variables.

En s'inspirant de la formule de type Gelfand-Levitan, on pouvait d{\'e}j{\`a}
s'interroger sur une classe relativement simple de fonctions
analytiques qui est celle des pseudo-polyn{\^o}mes, sous-classe des
fonctions enti{\`e}res de type exponentiel, et d{\'e}finis ainsi, 
$$P(\zeta_1,\ldots,\zeta_n,\exp<a_1,\zeta>,\ldots,\exp<a_k,\zeta>),$$
o{\`u} $\zeta\in\mathbb{R}^n$, $a_1,\ldots,a_k\in\mathbb{R}^n$,
$<a_j,\zeta>=a_j^1\zeta_1+\ldots+a_j^n\zeta_n$. Il existe effectivement des
estimations explicites donn{\'e}es par A. G. Khovanskii
dans~\cite{khovanskii}, dont on peut d{\'e}duire par exemple le r{\'e}sultat
suivant, signal{\'e} dans l'introduction :

\begin{theorem}\label{khovanskii}

Soit, pour $n,\,d\geq2,\;k\geq1$,
$$P_{n,k,d}=\left\{\sum_{|j|\leq d}c_j\,\zeta_1^{j_1}\ldots\zeta_n^{j_n}\,
e^{j_{n+1}<a_1,\zeta>}\ldots e^{j_{n+k}<a_k,\zeta>},\;
\zeta\in\mathbb{R}^n\right\},$$
une famille d'{\'e}l{\'e}ments de $C(I^s)$ param{\'e}tr{\'e}e par
un quasi-polyn{\^o}me {\`a}
coefficients $c_j\in C(I^s)$, {\`a} $n$ variables $\zeta_i$, $k$
pseudo-variables $e^{<a_i,\zeta>}$ et de degr{\'e} total $d$.

Alors on a : $\exists\,h\in\Lambda_{l,s},\;\forall\,\zeta\in\mathbb{R}^n$,
$$\|h-P_{n,k,d}\|_{\infty}\geq\frac{C(l,s)}{(k^2n\log n\log
  d)^{\frac{l}{s}}}.$$

Ce fa{\c c}on {\'e}quivalente, si $\mathcal{P}_{n,k,d}$ d{\'e}signe l'ensemble des
familles param{\'e}tr{\'e}es par des quasi-polyn{\^o}mes {\`a} $n$ variables, $k$
pseudo-variables et de degr{\'e} total $d$, on a
$$D_{n,k,d}(\Lambda_{l,s}):=\inf_{P\in\mathcal{P}_{n,k,d}}
\sup_{h\in\Lambda_{l,s}}\,\inf_{\zeta\in\mathbb{R}^n}\|h-P(\zeta)\|_{\infty}
\geq\frac{C(l,s)}{(k^2n\log n\log d)^{\frac{l}{s}}}\,.$$

En outre, la constante $C(l,s)$ peut {\^e}tre calcul{\'e}e et vaut
$$\frac{1}{\sqrt{s}\,2^{l+1}\,38^{\frac{l}{s}}([l]+1)^{[l]+1}(4(1+e))^{s([l]+1)}}\,.$$

\end{theorem}

\begin{remark}

Comme pour le cas polynomial, le r{\'e}sultat est aussi valable si on
consid{\`e}re l'espace $L^1(I^s)$ muni de la norme $\|\;\|_{L^1}$, avec
les familles quasi-polynomiales {\`a} coefficients $c_j\in L^1(I^s)$ (avec
une autre constante $C_{L^1}(l,s)$).

D'autre part, l'{\'e}nonc{\'e} reste valable pour $n=1$ ou $d=1$ quitte {\`a}
remplacer $n$ et $d$ par $n+1$ et $d+1$.

\end{remark}

\begin{proof}

La preuve est du m{\^e}me esprit que pour le th{\'e}or{\`e}me~\ref{vitushkin} : il
s'agit d'estimer de fa{\c c}on analogue le nombre de composantes connexes
de l'ensemble des z{\'e}ros d'un quasi-polyn{\^o}me en fonction de
$n,\,k,\,d$. Le nombre de cellules d'un tel ensemble vaut
au plus, {\`a} partir de ~\cite{khovanskii} (pour $p=1$),
$$2^{\frac{k(k-1)}{2}}d(n+d)^{n-1}(n(n+d)-n+1)^k\leq2^{\frac{k(k-1)}{2}}
d^{n+k}n^{n+2k}.$$
En reprenant la m{\^e}me m{\'e}thode que Warren dans
~\cite{warren}, on obtient une estimation du nombre de composantes de
l'ensemble $\mathbb{R}^n\setminus\bigcup_{j=1}^m\{\zeta,\;P_j(\zeta)=0\}$,
qui est
$$\sum_{j=0}^n2^{\frac{k(k-1)}{2}}(2d)^{n+k}n^{n+2k}2^jC_m^j<
2^{\frac{k(k-1)}{2}}(4emdn)^{n+2k},$$
o{\`u} $d=\max d_j$.

Pour $m\geq38k^2n\log n\log d$, ce nombre est inf{\'e}rieur {\`a} $2^m$, ce
qui donne, si $r$ est le plus petit entier $\geq2$ tel que
$r^s\geq38k^2n\log n\log d$ :
$\exists\,h\in\Lambda_{l,s},\;\forall\,\zeta\in\mathbb{R}^n$,
\begin{eqnarray*}
\|h-P_{n,k,d}(\zeta)\|_{\infty} & \geq &
\frac{1}{2M_{l,s}4^{s([l]+1)}r^l}\\
& \geq &
\frac{1}{\sqrt{s}\,2^{l+1}38^{\frac{l}{s}}([l]+1)^{[l]+1}
(4(1+e))^{s([l]+1)}}
\,\frac{1}{(k^2n\log n\log d)^{\frac{l}{s}}},
\end{eqnarray*}
ce qui termine la preuve.

\end{proof}

Ce r{\'e}sultat pr{\'e}sente cependant deux inconv{\'e}nients :
d'abord comme 
on l'a d{\'e}j{\`a} dit dans l'introduction, la minoration est relativement
faible par rapport {\`a} $k$ ; ensuite, c'est une classe 
trop restreinte car, dans le cadre de la th{\'e}orie d'approximation, on
ne pourra consid{\'e}rer que des familles du type
$$\sum_{j=(j_1,\ldots,j_{n+k})}c_j(x)\,\zeta_1^{j_1}\ldots\zeta_n^{j_n}
(\exp(<a_1,\zeta>))^{j_{n+1}}\ldots(\exp(<a_k,\zeta>))^{j_{n+k}},$$
o{\`u} $c_j\in C(I^s)$ (ou $L^1(I^s)$). En particulier, il n'y a pas
moyen de l'appliquer aux fonctions du genre $\exp(<\zeta,x>)$ o{\`u} les
variables $x\in I^s$ et $\zeta\in\mathbb{R}^n$ sont "m{\'e}lang{\'e}es" (comme
celles qui interviennent dans les formules du type Gelfand-Levitan).
\bigskip

On pouvait dans le m{\^e}me esprit avoir recours {\`a} la conjecture 
de Kouchnirenko, qui est une g{\'e}n{\'e}ralisation {\`a} plusieurs variables du
th{\'e}or{\`e}me de Descartes, et qui dit qu' un syst{\`e}me $P_1=\cdots=P_n=0$
d'{\'e}quations polynomiales {\`a} n variable, o{\`u} $m_i$ est le nombre de
termes de $P_i$, 
ne peut avoir plus de $(m_1-1)\ldots(m_n-1)$ racines positives non
d{\'e}g{\'e}n{\'e}r{\'e}es. Elle aurait pu nous {\^e}tre utile, mais elle a {\'e}t{\'e} r{\'e}cemment
infirm{\'e}e par B. Haas qui nous donne un contre-exemple
dans~\cite{haas}.
\bigskip

Dans notre cas o{\`u} on s'int{\'e}resse aux fonctions du type
$\exp(<\zeta,x>)$, on montre le r{\'e}sultat suivant, qui est une
application du th{\'e}or{\`e}me de Descartes :

\begin{proposition}\label{descartes}

On consid{\`e}re la famille $\Psi\subset C(I)$ d{\'e}finie par
$$\left\{\left(t\in[0,1]\mapsto\psi(t,\zeta)=\sum_{j=1}^nP_j(t)
\exp(\zeta_jt)\right),\;
\zeta=(\zeta_1,\ldots,\zeta_n)\in\mathbb{R}^n\right\},$$
o{\`u} $P_j$ est un polyn{\^o}me en $t$ de degr{\'e} $\leq p_j$.

On consid{\`e}re {\'e}galement la fonction d{\'e}finie sur 
$[0,1]=\bigcup_{i=1}^n\left[\frac{i-1}{n},\frac{i}{n}\right]$, par 
$$f_{\varepsilon}(x)=\varepsilon_i\frac{g_l(nx-i+1)}{M_l\,n^l},\;x\in
\left[\frac{i-1}{n},\frac{i}{n}\right],$$
o{\`u} $\varepsilon_i=1$, si $i$ pair, $-1$ sinon ($g_l$ d{\'e}finie comme
dans la partie~\ref{vitwar}).

On a alors
$$\inf_{\zeta\in\mathbb{R}^n}\|f_{\varepsilon}-\psi(\cdot,\zeta)\|_{\infty}
\geq\frac{C_l}{\left(n+\sum_{j=1}^np_j\right)^l}.$$

\end{proposition}

\begin{proof}

Supposons d'abord les $\zeta_j\in\mathbb{N}$ et les $P_j=c_j$
constants. Par le th{\'e}or{\`e}me de Descartes, le polyn{\^o}me
$\sum_{j=1}^nc_jX^{\zeta_j}$ a au plus $n-1$ racines $>0$. Si
$\zeta\in\mathbb{Z}^n$, on factorise par $X^{-k}$, $k$ assez
grand. Enfin pour $\zeta\in\mathbb{Q}^n$, on pose
$X=\exp\frac{t}{p},\,p$ d{\'e}nominateur commun des $\zeta_j$, ce qui
montre que la fonction
$$t\in\mathbb{Q}^n\mapsto\sum_{j=1}^nc_j\exp(\zeta_j t),$$
a au plus $n-1$ z{\'e}ros dans $\mathbb{R}$.

Or $f_{\varepsilon}$
s'annulant au moins une fois de plus que $\psi(\cdot,\zeta)$, il existe
au moins un sous-intervalle o{\`u} $f_{\varepsilon}$ et $\psi(,\cdot)$
sont de signes contraires, ce qui implique, $\forall\,\zeta\in\mathbb{Q}^n$,
$$\|f_{\varepsilon}-\psi(\cdot,\zeta)\|_{\infty}\geq\frac{C_l}{n^l}.$$
Enfin, l'assertion est encore valable pour $\zeta\in\mathbb{R}^n$ par
densit{\'e}, la convergence {\'e}tant uniforme sur $[0,1]$ (c'est encore une
application de l'id{\'e}e simple mais fondamentale,  qui dit qu'une
fonction qui ne s'annule 
pas beaucoup ne peut pas beaucoup osciller autour de $0$).
\bigskip

Consid{\'e}rons maintenant le cas g{\'e}n{\'e}ral
$$\psi(t,\zeta)=\sum_{j=1}^nP_j(t)\exp(\zeta_jt),$$
o{\`u} $P_j$ est un polyn{\^o}me de degr{\'e} $\leq p_j$. On commence par
remarquer que la 
fonction $t$ est limite uniforme sur $[0,1]$ de $\frac{\exp\eta
  t-1}{\eta}$, pour $\eta\rightarrow0_+$. On a aussit{\^o}t, pour tout $m\geq0$,
$$t^m=\lim_{\eta\rightarrow0_+}\left(\frac{\exp\eta t-1}{\eta}\right)^m=
\lim_{\eta\rightarrow0_+}\sum_{s=0}^m
\frac{(-1)^{m-s}C_m^s}{\eta^m}\exp(s\eta t).$$
Ainsi, chaque $P_j(t)\exp(\zeta_jt)$ va {\^e}tre limite uniforme sur
$[0,1]$ d'une famille de fonctions de la forme
$$\sum_{s=0}^{p_j}a_{j,s}(\eta)\exp(s\eta+\zeta_j)t,$$
qui pour tout $j,\,1\leq j\leq n$, poss{\`e}de au plus $p_j+1$ termes ; ce
qui pour $\psi(\cdot,\zeta)$ donnera $\sum_{j=1}^np_j+n$ termes et
aboutira, par limite uniforme, {\`a}
$$\|f_{\varepsilon}-\psi(\cdot,\zeta)\|_{\infty}\geq\frac{C_l}
{\left(n+\sum_jp_j\right)^l}\,.$$

\end{proof}

\begin{remark}

On peut bien s{\^u}r prolonger le r{\'e}sultat sur tout intervalle $[a,\,b]$
de $\mathbb{R}$ ({\`a} condition qu'il soit compact pour assurer la
convergence uniforme de $\frac{\exp(\eta t)-1}{\eta}$ vers $t$).

\end{remark}

On a alors un r{\'e}sultat qui pourrait s'appliquer dans le cadre de
notre probl{\`e}me inverse. L'inconv{\'e}nient est qu'ici aussi la minoration
est trop faible : en effet, la formule de type Gelfand-Levitan poss{\`e}de en
tant que d{\'e}terminant tous les 
$$\exp2(b_1\zeta_1+\cdots b_N\zeta_N),\;b_j\in\{-1,\,0,\,1\},$$ 
qui sont au nombre de $3^N$, sans compter ceux qui ont une partie
polynomiale, ce qui ne donnera pas mieux que $\frac{C_l}{3^{lN}}$,
qui est d{\'e}j{\`a} insuffisant.

En revanche, l'avantage ici est qu'on n'a pas besoin de
borner les variables $(\zeta_1,\ldots,\zeta_n)$. Ainsi de deux choses
l'une : soit la proposition~\ref{descartes} peut {\^e}tre nettement
am{\'e}lior{\'e}e afin d'obtenir une estimation voisine de celle donn{\'e}e dans le
th{\'e}or{\`e}me~\ref{th1}. Soit au contraire elle n'est pas loin d'{\^e}tre
optimale, et la comparaison avec l'estimation en $\frac{1}{(N\log
  N)^l}$ montre que c'est alors pour des grandes valeurs des
param{\`e}tres $\zeta_1,\ldots,\zeta_n$ (donc au-del{\`a} de la taille
polynomiale) que l'approximation est meilleure. 

De plus, la fonction
$f_{\varepsilon}$ est compl{\`e}tement explicit{\'e}e, et ne d{\'e}pend pas de la
famille exponentielle impos{\'e}e au d{\'e}part. Elle se trouve ainsi
uniform{\'e}ment distante de toutes les familles exponentielles {\`a} $n$
termes. Ceci mis {\`a} part, le r{\'e}sultat donn{\'e} par le th{\'e}or{\`e}me~\ref{th1}
est pr{\'e}f{\'e}rable pour la pr{\'e}cision de l'ordre de $\frac{1}{(N\log
  N)^{\frac{l}{s}}}$, et pour sa forme plus g{\'e}n{\'e}rale qui ne concerne
pas seulement les familles exponentielles, mais toutes les fonctions
enti{\`e}res de type exponentiel.

\end{document}